%% file: hypeq_bd.tex
\documentclass[11pt,oneside,english]{article}

\include{defs_hypeq_bd}

\author{Roberto \titsc{Frigerio} 
\and Carlo \titsc{Petronio}}

\title{Construction and Recognition of\\
Hyperbolic 3-Manifolds with Geodesic Boundary}

\begin{document}

\maketitle

\begin{abstract} We extend to the context of hyperbolic 3-manifolds with
geodesic boundary Thurston's approach to hyperbolization by means of
geometric triangulations.  In particular, we introduce moduli for
(partially) truncated hyperbolic tetrahedra, and we discuss consistency
and completeness equations.  Moreover,
building on previous work of Ushijima, we extend Weeks' tilt formula
algorithm, which computes the Epstein-Penner canonical triangulation, to
an algorithm that computes the Kojima triangulation.  The theory is
particularly interesting in the case of complete finite-volume manifolds
with geodesic boundary in which the boundary is non-compact. We include
this case using a suitable adjustment of the notion of ideal
triangulation, and we show that the case naturally arises within the
theory of knots and links.

\vspace{4pt}

\noindent MSC(2000): 57M50 (primary), 57M25 (secondary).

\end{abstract}

\tableofcontents

\section*{Introduction}
The aim of this paper is to lay down the theoretical background for a census
of orientable hyperbolic 3-manifolds with geodesic boundary. Our starting
point is the idea
of turning the construction of the hyperbolic structure on a manifold 
into an algebraic problem. This idea is 
originally due to Thurston~\cite{thurston:notes} for the case
of cusped manifolds, and has been systematically
exploited by the software SnapPea~\cite{snappea}.
In the setting of cusped manifolds one employs ideal tetrahedra,
which are parameterized by complex numbers, and tries to solve
the consistency and completeness equations. In the bounded case 
one has to consider truncated tetrahedra, and moduli 
get more complicated, but basically the whole scheme extends.
The two phenomena of non-compactness and
presence of geodesic boundary can actually occur simultaneously, 
and, following Kojima~\cite{kojima,kojima-bis}, 
we introduce the notion of partially truncated tetrahedron to deal with 
this fact. One interesting point emerges when the boundary of a 
finite-volume hyperbolic manifold is itself non-compact.
Namely, we show that in this case the combinatorial datum to start from
to build the structure is not an ideal triangulation of the original manifold,
but of a certain quotient of the original manifold.

Working with moduli and equations one can \emph{construct} hyperbolic
manifolds with boundary, but, after a list of manifolds has been put together, one
has to remove duplicates to get the genuine list, so one
is naturally faced with the issue of \emph{recognizing}
the manifolds. It turns out that for both cusped and bounded manifolds a 
certain \emph{canonical} decomposition exists, due to
Epstein and Penner~\cite{epstein-penner} in the former case and to
Kojima~\cite{kojima,kojima-bis} in the latter. One natural strategy to recognize
a manifold decomposed into geometric pieces is then to modify the decomposition
until the canonical one is reached. This method was described 
by Weeks in~\cite{weeks:tilt} in terms of a so-called
``tilt formula'', and it was used in~\cite{census} in the cusped case.
The tilt formula itself was already discussed for the bounded case
in~\cite{ushijima}, 
and we describe in this paper the whole strategy to turn
an arbitrary triangulation of a manifold with boundary into its
Kojima decomposition. 
We warn the reader that, both in the cusped and in the bounded case,
the algorithm to transform a decomposition into the canonical one is
not proved to converge in general, but, at least in the cusped case, 
it usually does in practice.

Various differences arise between the cusped and the bounded case, and it
is maybe worth mentioning here at least the most subtle one, which requires
quite some effort to deal with. Just as the Epstein-Penner decomposition in the
cusped case, the Kojima decomposition for bounded manifolds is obtained
by projecting to hyperbolic 3-space the faces of a certain polyhedron in 
Minkowski 4-space. In both cases the polyhedron is the convex hull of 
certain points that represent, in a suitable sense, liftings of
cusps and of boundary components. When there are cusps only, the height of the
liftings is intrinsically determined \emph{a priori} (up to global rescaling), 
and the basic idea to modify
a triangulation into the canonical one is to lift the ideal tetrahedra
with vertices at the lifted cusps, and make sure the lifted tetrahedra 
bound a convex set. Essentially the same happens when 
the boundary is non-empty but there are no cusps at all.
In the mixed case, however, only boundary components have a prescribed
height to be lifted at, while the height for cusps is a lot 
harder to determine. This matter is discussed in 
Sections~\ref{height:section} and~\ref{computing:section}.

We believe that the issue of understanding and enumerating hyperbolic
3-ma\-ni\-folds with geodesic boundary is a very natural one, and we are
planning to exploit the theory developed in this paper in the close future,
building an analogue ``with boundary'' of the cusped census of~\cite{census}.
Here are three specific reasons for caring about manifolds with boundary:
\begin{itemize}
\item  These manifolds still satisfy the rigidity theorem, so every
geometric invariant, such as the volume or the length spectrum, 
is actually a topological invariant;
\item Thurston's hyperbolization theorem for Haken
manifolds~\cite{thurston:bams} 
implies that all manifolds with boundary
satisfying some very natural and fairly general topological properties
actually are hyperbolic, so one expects to find that ``most'' manifolds
with boundary are hyperbolic;
\item If $L$ is a link in $\matS^3$ and $\Sigma$ is a minimal-genus
Seifert surface for $L$, then the manifold obtained by cutting $\matS^3$ along
$\Sigma$ is a natural candidate for a finite-volume hyperbolic structure
with boundary. In addition, 
$L$ has a well-defined ``length'' with respect to this structure (if any).
\end{itemize}
Given the length and comparative variety of topics touched in the paper, we have 
included at the beginning of each section a couple of explanatory paragraphs,
where we outline the contents of the section and we list the statements,
definitions and notations used later in our work.
The reader willing to reach the core of our arguments may at first 
concentrate on this material only.

\section[Triangulations of hyperbolic 3-manifolds with geodesic boundary]
{Triangulations of hyperbolic 3-manifolds with\\ geodesic boundary}\label{trunc:tria:section}

In this section we 
prove some preliminary facts about the topology and geometry at infinity of 
a finite-volume orientable hyperbolic 3-manifold with geodesic boundary.
We also explain what do we mean by a triangulation of such a manifold,
showing in particular that this notion must be understood with some care when
the boundary of the manifold is non-compact.
The essential points of this section are 
Proposition~\ref{compactification:prop},
Definitions~\ref{partial:trunc:def},~\ref{top:real:def},~\ref{geom:real:def},
~\ref{top:tria:def},~\ref{geom:tria:def},~\ref{Nprimo:def}, and~\ref{rel:tria:def},
and Proposition ~\ref{tria:of:N':prop}.
However, Proposition~\ref{atoroidal:prop} 
and the discussion following it are also quite important as a motivation.

\paragraph{Natural compactification}
Let $N$ be a complete 
finite-volume orientable hyperbolic 3-manifold with geodesic boundary.
(In the rest of the paper we will
summarize all this information saying just that \emph{$N$ is hyperbolic}.)
We denote by $D(N)$ the double of $N$, 
\emph{i.e.}~the manifold obtained by mirroring $N$ in its boundary.
Now $D(N)$ is an orientable finite-volume hyperbolic 3-manifold without boundary, so it 
consists of a compact portion together with several cusps of the form
$T\times[0,\infty)$, where $T$ is the torus ---see \emph{e.g.}~\cite{BePe:libro}. 
Within $D(N)$ we have the surface
$\partial N$ which cuts $D(N)$ into two isometric copies of $N$, and to
understand the geometry of the ends of $N$ we must investigate how $\partial N$
can intersect a cusp $T\times[0,\infty)$. Using the geometry of $T\times[0,\infty)$
one sees that, up to resizing the cusp, either $\partial N$ is disjoint from
$T\times[0,\infty)$ or it is given by $\gamma\times[0,\infty)$,
where $\gamma$ is the union of a finite number of parallel geodesic loops on $T$.
In the first case the cusp $T\times[0,\infty)$
is contained in one of the two isometric copies of $N$. In the second case,
knowing that $\partial N$ is separating in $D(N)$, we see that $\gamma$ contains
at least two loops, and $N$ has an end of the form $A\times[0,\infty)$ where
$A\subset T$ is an annulus bounded by these two loops. Since the double of $A$ already is
a torus, we also see that $\gamma$ consists of precisely two loops.

The previous discussion shows that $N$ consists of a compact portion together
with some cusps based either on tori or on annuli, which implies the following:

\begin{prop}\label{compactification:prop}
If $N$ is hyperbolic (\emph{i.e.}~$N$ is a complete 
finite-volume orientable hyperbolic $3$-manifold with geodesic boundary)
then it has a natural compactification $\Nbar$
obtained by adding some tori and annuli.
\end{prop}

In particular, $\partial N$,
which we know~\cite{kojima} 
to be a finite-area orientable hyperbolic surface, can be non-compact.
Moreover the ends of $\partial N$ naturally come into pairs
$\{\pm1\}\times S^1\times[0,\infty)=\partial( [-1,1]\times S^1\times[0,\infty))$.
For later purpose we denote by $\calA\subset\Nbar$ the family of annuli
added to compactify $N$. No specific notation for the tori is needed.

\begin{rem}\label{rectangular:annuli:rem}
\emph{If a cusp of $N$ is based on a torus, it is well-known that
this torus has a Euclidean structure well-defined up to rescaling.
Now, if a cusp is based on an annulus, its double is a Euclidean torus, so
the annulus is itself Euclidean with geodesic boundary, up to rescaling.
In particular, the annulus is obtained from a Euclidean rectangle
by identifying two opposite edges.
So, if we normalize the width of the annulus to unity, we can assign the
annulus a well-defined \emph{length}.}
\end{rem}

\paragraph{Topological restrictions}
We have shown so far that a
hyperbolic $N$ is obtained from a compact $\Nbar$ by removing from $\partial\Nbar$ some toric 
components and a family $\calA$ of 
closed embedded annuli. We also know that the components
of $\partial N$ are hyperbolic surfaces, whence:

\begin{prop}\label{neg:chi:prop}
The components of $\partial N$ have 
negative Euler characteristic. 
\end{prop}

\begin{cor}\label{no:spere:cor}
$\partial\Nbar$ does not contain spheres, 
and no annulus of $\calA$ can lie
on a toric component of $\partial\Nbar$.
\end{cor}

\dimostraz
There cannot be a sphere because
an innermost annulus on a sphere bounds an open disc, having $\chi=1$.
For the same reason on a toric
component there cannot be trivial annuli, so there are some parallel annuli,
and the complement also consists of annuli, having $\chi=0$.
\finedimo

This lemma shows that from the pair $(\Nbar,\calA)$ determined by $N$
we can get back $N$ in a non-ambiguous way by removing 
from $\Nbar$ both $\calA$ and all the toric components of $\partial\Nbar$.
We also have the following additional topological restrictions, 
stated separately because harder to check directly when
an \emph{a priori} non-hyperbolic $N$ is given.

\begin{prop}\label{atoroidal:prop}
The compact manifold $\Nbar$ is irreducible and geometrically atoroidal. 
Moreover $\Nbar\setminus\calA$ is boundary-incompressible and the only proper essential
annuli it contains are parallel in $\Nbar$ to the 
annuli in $\calA$.
\end{prop}

\dimostraz
Of course $\Nbar$ is irreducible, because its double $\overline{D(N)}$ is. An embedded 
incompressible torus must be boundary parallel 
in $\overline{D(N)}$, whence also in $\Nbar$. 

The toric boundary components of $\Nbar\setminus\calA$ 
are incompressible because they are  
in $\overline{D(N)}$. Let $\Delta$ be a disc that compresses a loop $\gamma$ 
contained in a component $\Sigma$ of $\partial N$. Then $\Delta$
lifts to the universal cover of $N$, and the lifting of $\gamma$ lies
on a hyperbolic plane that covers $\Sigma$. It readily follows that $\gamma$
must be trivial in $\Sigma$.

An essential annulus cannot join two toric components of $\partial(\Nbar\setminus\calA)$,
otherwise it would in $\overline{D(N)}$.
It also cannot join a toric component with a non-toric one, otherwise its
double would join two tori in
$\partial\overline{D(N)}$. If an essential annulus joins
two non-toric components of $\partial (\Nbar\setminus\calA)$ 
then it lies in $N$, and 
its double is an essential torus in $D(N)$. This torus must be
boundary-parallel, which easily implies that 
the annulus is parallel to $\calA$ in $\Nbar$.
\finedimo

In the previous statement one should notice that irreducibility 
holds for $\Nbar$ if and only if it holds for $N$, and similarly for
atoroidality, whereas boundary-incompressibility for
$\Nbar\setminus\calA$ does not imply the same property for $\Nbar$. It is also not
possible to deduce from the statement that $\Nbar$ is anannular.
Note however that $\calA=\emptyset$ when in $N$ there are no annular cusps.

\paragraph{Links and Seifert surfaces}
We show in this paragraph that manifolds satisfying (most of) the topological
restrictions of Propositions~\ref{neg:chi:prop} 
and~\ref{atoroidal:prop} naturally arise in the context of
the theory of knots and links. Namely, let $L\subset\matS^3$ be a link, and let $\Sigma$
be an orientable Seifert surface for $L$. Thicken $\Sigma$ to a product
$\Sigma\times[-1,1]\subset\matS^3$ so that $\Sigma=\Sigma\times\{0\}$,
and define $N$ as $\matS^3\setminus(\Sigma\times(-1,1))$. Note that $N$
compactifies to a manifold $\Nbar$ by adding the annuli
$\calA=L\times[-1,1]$ that define the null framing on   the
components of $L$. Moreover:
\begin{itemize}
\item If $\partial\Nbar$ contains a sphere then $L$
has a trivial component unlinked from the rest;
\item If $\partial\Nbar$ contains a torus then $L$ contains
two parallel components;
\item If $\Nbar$ is not irreducible then $L$ is a split link;
\item If $\Nbar$ is not atoroidal then $L$ is a satellite of a non-trivial
knot $K$ and $L$ is homologically trivial in the neighbourhood of $K$;
\item If $\partial N$ is compressible then $\Sigma$ is the result
of a stabilization of another Seifert surface; in particular, $\Sigma$
cannot have minimal genus.
\end{itemize}
These remarks provide rather flexible sufficient conditions for $\Nbar$
to satisfy most of the topological requirements for
hyperbolicity. The restriction that essential annuli in $N$ should
be parallel to $\calA$ in $\Nbar$ is more involved, and it is not
addressed here.

\paragraph{Partially truncated tetrahedra}
Recall that, when $N$ is finite-volume non-com\-pact hyperbolic
and $\partial N=\emptyset$, it is typically
possible to decompose $N$ into pieces isometric to geodesic ideal tetrahedra in
$\matH^3$, and in practice the hyperbolic structure of $N$ is 
constructed by first taking a topological ideal triangulation and then
choosing the geometric shape of the tetrahedra so that their structures match
under the gluings giving a complete structure on $N$. 
Our wish in the rest this of section is to extend the notion of
ideal triangulation to the case of hyperbolic manifolds with geodesic boundary.
We begin by describing the pieces into which manifolds will be decomposed,
first topologically and then geometrically.

\begin{defn}\label{partial:trunc:def}
\emph{We call \emph{partially truncated tetrahedron} a triple $(\Delta,\calI,\calZ)$
where $\Delta$ is a tetrahedron, $\calI$ is a set of vertices of $\Delta$,
and $\calZ$ is a set of edges of $\Delta$ such that neither of the two endpoints of 
an edge in $\calZ$ belongs to $\calI$. The elements of $\calI$ and $\calZ$ 
will be called \emph{ideal vertices} and \emph{length-$0$ edges} respectively, for a reason
to be explained soon. In the sequel we will always refer to $\Delta$ itself as a
partially truncated tetrahedron, tacitly implying that certain $\calI$ and
$\calZ$ are also fixed.}
\end{defn}

\begin{defn}\label{top:real:def}
\emph{Given a partially truncated tetrahedron $\Delta$
we define its \emph{topological realization}
as the space $\Delta\!^*$ 
obtained by removing from $\Delta$ the ideal vertices, the length-0
edges, and small open stars of the non-ideal vertices. 
We will call \emph{lateral hexagon}
and \emph{truncation triangle} 
the intersection of $\Delta\!^*$ respectively with a face of $\Delta$ and with
the link in $\Delta$ of a non-ideal vertex.
The edges of the truncation triangles,
which also belong to the lateral hexagons, will be called
\emph{boundary edges}. The other edges of the lateral hexagons
will be called \emph{internal edges}.}
\end{defn}

Note that, if $\Delta$ has length-$0$ edges, some vertices 
of a truncation triangle may be missing.
Similarly, if $\Delta$ has ideal vertices or 
length-0 edges, a lateral hexagon of $\Delta\!^*$
may not quite be a hexagon, because some of its
(closed) edges may be missing. Note however that 
two consecutive edges cannot both be missing.

\begin{defn}\label{geom:real:def}
\emph{Given a partially truncated tetrahedron $\Delta$
we call \emph{geometric realization} of $\Delta$ an embedding
of $\Delta\!^*$ in $\matH^3$ such that:
\begin{enumerate}
\item   The truncation triangles are geodesic triangles, with ideal
vertices corresponding to missing vertices;
\item The lateral hexagons are geodesic polygons, with ideal vertices 
corresponding to missing edges;
\item Truncation triangles and lateral 
hexagons lie at right angles to each other.
\end{enumerate}}
\end{defn}

\begin{figure}
\begin{center}
\input{geomreal.pstex_t}
\caption{\small{A geometric tetrahedron
with one ideal vertex and one length-0 edge.}}\label{geomreal:fig}
\end{center}\end{figure}
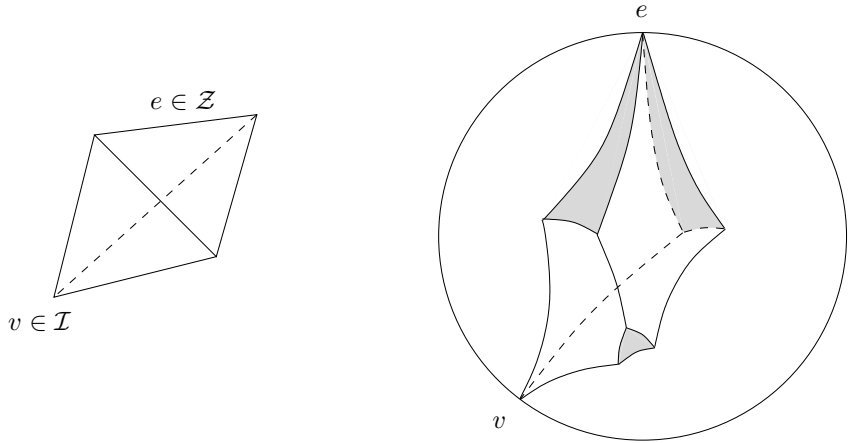

\noindent An example of geometric realization is shown in Fig.~\ref{geomreal:fig}, where
truncation triangles are shadowed.

\begin{rem}\label{geom:of:tetra:rem}
\emph{If $\Delta\!^*$ is a geometric realization of $\Delta$ and $v$ is an ideal
vertex of $\Delta$ then a neighbourhood of $v$ intersected with $\Delta\!^*$ is automatically
isometric in the half-space model 
$\Hhalf^3=\matC\times(0,\infty)$ 
of $\matH^3$ to $W\times[t_0,\infty)$, where $W\subset\matC$ is a
triangle and $t_0>0$. Similarly, if $e$ is a length-0 edge then
a neighbourhood of $e$ intersected with $\Delta\!^*$ is isometric to 
$[-1,1]\times[-b,b]\times[t_0,\infty)$. Here the triangles
$\{\pm 1\}\times [-b,b]\times[t_0,\infty)$ are contained in the truncation
triangles, the triangles $[-1,1]\times\{\pm b\}\times[t_0,\infty)$ 
are contained in the lateral hexagons, and the closure in $\Delta$ of
every triangle $\{x\}\times[-b,b]\times[t_0,\infty)$ is
obtained by adding only one point of $e$, so that the 
segment $[-1,1]\times\{0\}\times\{\infty\}$ can be viewed as a subset of $e$.}
\end{rem}

\paragraph{Triangulations}
In the language introduced above, the classical
notion of \emph{ideal triangulation} of a compact 3-manifold with boundary is a 
realization of the interior of the manifold as a gluing of some $\Delta\!^*$'s, where
the corresponding $\Delta$'s have all ideal vertices (and hence no length-0 edge)
and the gluing is induced by a simplicial pairing of the faces of the $\Delta$'s.
We can now easily extend this notion to the situation we
are interested in.

\begin{defn}\label{top:tria:def}
\emph{Let $\Nbar$ be a compact orientable manifold and let
$\calA\subset\partial\Nbar$ be a family of disjoint annuli
not lying on the toric components of $\partial\Nbar$.
Let $N$ be obtained from $\Nbar$ by removing 
$\calA$ and the toric components of $\partial\Nbar$.
We define a \emph{partially truncated triangulation} of $N$ to be 
a realization of $N$ as a gluing of some $\Delta\!^*$'s along
a pairing of the lateral hexagons induced by a simplicial
pairing of the faces of the $\Delta$'s.}
\end{defn}

\begin{rem}\label{trunc:tria:first:rem}
\begin{enumerate}
\item[\emph{1.}] \emph{Under the pairing of the faces of the $\Delta$'s, ideal vertices are
matched to each other. Similarly, length-0 edges are matched to each other.}
\item[\emph{2.}] \emph{The truncation triangles of the $\Delta\!^*$'s give a triangulation of
$\partial N$ with some genuine and some ideal vertices.}
\item[\emph{3.}] \emph{The links of the ideal vertices of the $\Delta$'s give a triangulation
of the toric components of $\partial\Nbar$.}
\item[\emph{4.}] \emph{The links of the 
length-0 edges of the $\Delta$'s give a decomposition
into rectangles of the annuli in $\calA$. On each rectangle 
$[-1,1]\times[-b,b]$ only
the two opposite edges $[-1,1]\times\{\pm b\}$
that lie on lateral hexagons get glued to other 
rectangles, while the two opposite edges $\{\pm 1\}\times[-b,b]$
that lie on truncation triangles
contribute to the boundary of $\calA$.}
\end{enumerate}
\end{rem}

\begin{defn}\label{geom:tria:def}
\emph{Let $N$ as above be endowed with a hyperbolic structure.
A partially truncated triangulation of $N$ is called \emph{geometric}
if, for each tetrahedron $\Delta$ of the triangulation, the pull-back
to $\Delta\!^*$ of the Riemannian metric of $N$ defines a
geometric realization of $\Delta$.
Equivalently, the hyperbolic structure of $N$ should be obtained by gluing
geometric realizations of the $\Delta$'s 
along isometries of their lateral hexagons.}
\end{defn}

In Section~\ref{kojima:section} we will carefully describe Kojima's 
result~\cite{kojima}
according to which every hyperbolic $N$
as above has a \emph{canonical} decomposition into partially 
truncated \emph{polyhedra}, rather than tetrahedra. Just as 
it happens with the Epstein-Penner decomposition~\cite{epstein-penner} 
of non-compact manifolds with empty boundary, 
in the vast majority of cases 
the Kojima decomposition actually consists of tetrahedra, or at least 
can be subdivided into a geometric partially truncated triangulation.

\paragraph{Manifolds with arcs}
Our aim is to employ partially truncated triangulations to construct
and understand hyperbolic manifolds with boundary, just as ideal
triangulations are employed in the cusped case without boundary. 
One disadvantage of partially truncated triangulations when compared to ideal ones
is that the length-0 edges break the symmetry of the tetrahedron, so the
situation may appear to be less flexible. It is a useful and remarkable fact
that a partially truncated triangulation of a given manifold actually
corresponds to a genuine ideal triangulation of another manifold, 
as we will now explain.

\begin{defn}\label{Nprimo:def}
\emph{Given a manifold $N$ that compactifies to an
$\Nbar$ by adding some tori and a family
$\calA$ of annuli, we define $N'$ as the quotient of $\Nbar$ in which
every annulus $[-1,1]\times S^1\in\calA$ is collapsed to an arc $[-1,1]\times\{*\}$.
Note that $N'$ is also a compact
manifold and $[-1,1]\times\{*\}$ is an arc properly embedded in $N'$, as one
readily checks by visualizing a neighbourhood of $[-1,1]\times S^1$ as 
$[-1,1]\times\{z\in\matC: 1\leqslant |z| < 2\}$.
We denote by $\alpha_N$ the family of all the arcs $[-1,1]\times\{*\}$ in $N'$.}
\end{defn}

\begin{defn}\label{rel:tria:def}
\emph{If $M$ is compact and $\beta$ is a family of disjoint properly
embedded arcs in $M$, we call \emph{ideal triangulation} of the pair $(M,\beta)$ 
an ideal triangulation of $M$ that contains as edges 
all the arcs in $\beta$.}
\end{defn}

\begin{prop}\label{tria:of:N':prop}
Partially truncated triangulations of $N$ bijectively correspond
to ideal triangulations of $(N',\alpha_N)$.
\end{prop}

\dimostraz
A partially truncated tetrahedron $\Delta\!^*$ is turned into
an ideal one by removing the truncation triangles and adding
the length-0 edges minus their ends. On the manifold
this corresponds to removing the boundary and
collapsing each rectangle $[-1,1]\times[-b,b]$ as
in Remark~\ref{trunc:tria:first:rem}(4) to
$[-1,1]\times\{*\}$, see also Remark~\ref{geom:of:tetra:rem}.
So we get precisely the interior of $N'$
with the arcs in $\alpha_N$ being the length-0 edges.

An ideal triangulation of $(N',\alpha_N)$ is turned into a 
partially truncated triangulation of $N$ by declaring
to be length-0 the edges in $\alpha_N$ and to be ideal 
the vertices on the tori of $\partial N'$  on which there 
are no ends of arcs in $\alpha_N$. 
\finedimo

Having seen how partially truncated triangulations relate to ideal
ones, it is natural to ask whether the Matveev-Piergallini 
calculus~\cite{matveev,piergallini} for ideal triangulations 
generalizes to the case of manifolds with arcs. 
Recall that the fundamental move of this calculus
is the two-to-three move (shown below in Fig.~\ref{dual:mapimove:fig})
which destroys a triangle and the two tetrahedra 
incident to it, and creates one edge and three tetrahedra incident to this edge.
Of course a positive two-to-three move can always be applied to an
ideal triangulation of $(M,\beta)$, while the inverse three-to-two move
can be applied as long as the edge it destroys does not lie in $\beta$.
The next result is due to Amendola 
(see also~\cite{BaBe} and~\cite{TuVi}). 
Since it is not strictly speaking necessary for the present paper,
we omit its proof.

\begin{teo}\label{rel:alpha:teo}
Let two ideal triangulations of the same $(M,\beta)$ be given.
Assume both triangulations contain at least two tetrahedra. Then they
are related to each other
by a finite combination of two-to-three moves and 
three-to-two moves that do not destroy the edges of $\beta$.
\end{teo}

\section[Moduli and equations for partially truncated 
tetrahedra]{Moduli and equations for\\ partially truncated tetrahedra}
\label{mod:eqns:section}
In this section we introduce moduli for the geometric realizations of
partially truncated tetrahedra, and we describe the equations ensuring that
a gluing of geometric tetrahedra gives rise to a consistent and 
complete hyperbolic structure with
geodesic boundary. The idea here is to start with a 
topological triangulation of a certain
manifold with boundary, and
try to construct its geometric structure, if any, by choosing the geometric
shape of the tetrahedra in the triangulation.
We devote the initial part of the section to putting this idea in context
and providing motivations, and only then we turn to moduli and equations.

\begin{figure}
\begin{center}
\input{newtronco.pstex_t}
\caption{\small{Notation for edges and dihedral angles 
of a truncated tetrahedron. The geometric figure on the right uses the projective model of
$\matH^3$, as explained below in 
Section~\ref{kojima:section}.}}\label{tetra:notation:fig}
\end{center}\end{figure}
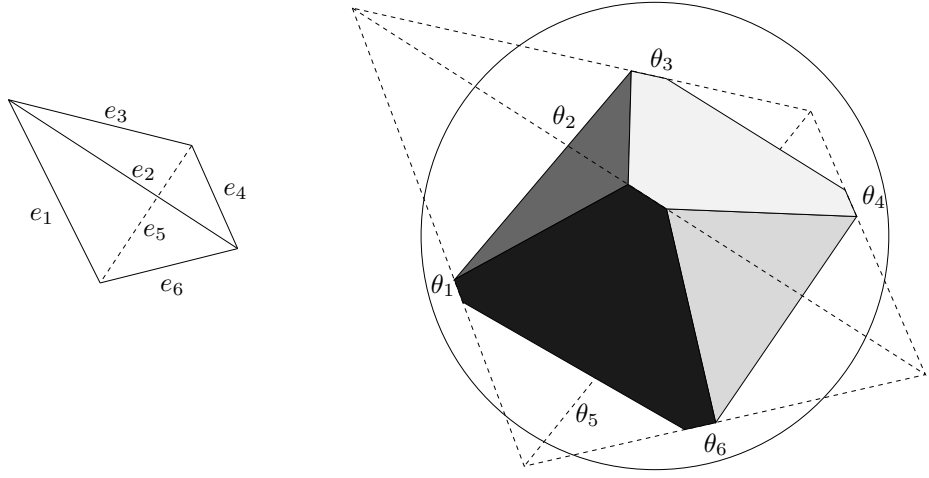

Also in this section we single out the very basic
points on which the reader could first concentrate. Moduli are introduced
in Theorem~\ref{moduli:teo}, and consistency equations in Theorem~\ref{consistency:teo},
with notation coming from formulae~(\ref{boundary:edge:length:formula}) 
to~(\ref{small:z:def:formula}) and Fig.~\ref{tetra:notation:fig}.
Completeness equations are shown
to be essentially the same as in the cusped case, and informally discussed after 
Remark~\ref{after:consistency:rem}.
In Theorem~\ref{uniqueness:teo} we also show that a solution, if any, is unique.

\paragraph{Hyperbolization with boundary}
The results of the previous section show that to build a census
of hyperbolic 3-manifolds with geodesic boundary we should first list,
according to some natural ordering,
all the pairs $(\Nbar,\calA)$ where $\Nbar$ is compact and $\calA\subset\Nbar$
is a family of disjoint annuli not lying on the toric components
of $\partial\overline{N}$.
For each such pair we should then consider
the manifold $N$
obtained by removing from $\partial\Nbar$ all the toric components
and the annuli of $\calA$, and discard the pair if the conditions of
Propositions~\ref{neg:chi:prop} or~\ref{atoroidal:prop} are violated.
For each remaining pair we should test the corresponding $N$
for hyperbolicity.

Looking more closely at the strategy just described, one sees that
it is very easy to describe an algorithm that lists, \emph{with repetitions},
all the pairs $(\Nbar,\calA)$ such that the corresponding $N$ 
satisfies the condition
of Proposition~\ref{neg:chi:prop}.
The conditions of Proposition~\ref{atoroidal:prop} are harder to check
but still manageable, at least theoretically, by means of the technology of
normal surfaces. Now, if a pair $(\Nbar,\calA)$ survives the topological tests,
we see that the double $D(N)$ of the corresponding $N$ has ends of the form
$T\times[0,\infty)$ and compactifies to
an irreducible, atoroidal, and anannular 3-manifold.
Assuming either that $N$ is non-compact or that $\partial N$ is non-empty,
we see that $D(N)$ is Haken, so Thurston's 
hyperbolization theorem~\cite{thurston:bams} shows that 
$D(N)$ is finite-volume hyperbolic,
and the involution of $D(N)$ that fixes $\partial N$ and
interchanges $N$ with its mirror copy can be realized by an 
isometry~\cite{fuji2}, so $N$ also is hyperbolic.
However, this theoretical proof of existence of 
the hyperbolic structure is not satisfactory under
at least two respects. First, it does not allow to compute the geometric 
invariants of $N$, such as the volume. Second,
it leaves unsettled the issue of removing duplicates from the list of manifolds.

The alternative strategy based on triangulations which we will now describe
overcomes both the drawbacks of the topological approach 
just pointed out. It should be noted,
however, that \emph{a priori} there could exist hyperbolic manifolds that
cannot be triangulated geometrically. These manifolds would be missed by 
our search.

\paragraph{Enumeration strategy}
To employ triangulations, we switch from the $(\Nbar,\calA)$ to the
$(N',\alpha_N)$ compactification of the candidate hyperbolic $N$. So, 
our first step is to list all pairs $(\calT,\alpha)$
where $\calT$ is an ideal triangulation of some 
compact orientable 3-manifold $N'$ with boundary, 
and $\alpha$ is a set of edges of $\calT$, also viewed as a set of properly
embedded arcs in $N'$. Of course there are infinitely many such $(\calT,\alpha)$'s,
so in practice one always deals with some finite ``initial'' segment of the list.
A pair $(\calT,\alpha)$ is immediately
discarded if $\partial N'$ contains spheres on which there are two 
or fewer ends of the arcs in $\alpha$. If $(\calT',\alpha)$ is not discarded,
we define $\Nbar$ as $N'$ minus an open tubular neighbourhood for each arc 
in $\alpha$, and $\calA$ as the family of annuli that bound the removed tubes.
Now we can define $N$ as $\Nbar$ minus
$\calA$ and the boundary tori, and $N$ 
automatically satisfies the condition of 
Proposition~\ref{neg:chi:prop}.

A pair $(\calT,\alpha)$ gives rise to a partially truncated triangulation
of the corresponding $N$ by declaring to be length-$0$
the edges in $\alpha$, and to be ideal the vertices corresponding to the
toric components of $\partial\Nbar$ on which there is no end of any arc in $\alpha$.
We will prove in the rest of this section
that there exists an algorithm to answer the
question whether \emph{can the tetrahedra of $\calT$ be 
geometrically realized in $\matH^3$ so to define a complete hyperbolic structure on $N$}.
If the answer is affirmative then we add $N$ to our census, if not we pass to the
next $(\calT,\alpha)$. 

\begin{rem}\label{change:tria:rem}
\emph{To make the search more effective, a slight modification of the
method just described can be employed. Namely, when 
the shape of the elements of a triangulation $\calT$ of a certain $(N',\alpha)$
cannot be chosen to give a complete structure on the corresponding $N$,
it is often convenient, before giving up, to try with other
triangulations of the same $(N',\alpha)$. It typically happens, at least in
the non-compact empty-boundary case dealt with by SnapPea~\cite{snappea},
that eventually a triangulation is found 
that either is geometric or suggests
which of the topological restrictions of Proposition~\ref{atoroidal:prop}
is violated.}
\end{rem}

The outcome of the strategy just outlined is a list of 
hyperbolic manifolds, each with a certain geometric triangulation.
However this list contains repetitions, that we can remove if
we can recognize manifolds.
Concentrating on those which have non-empty boundary,
we note that each of them has a unique well-defined Kojima decomposition.
This decomposition can now be viewed as the \emph{name} of the manifold,
because two such decompositions can be checked to be equal or not by comparing
the geometric shape of the polyhedra and the combinatorics of the gluings.
The recognition issue is then reduced to the issue of constructing
the Kojima decomposition starting from an arbitrary geometric triangulation.
This is the theme we will concentrate on starting from the next section.

\paragraph{Moduli}
We will now show that the dihedral angles at the non-$0$-length edges can
be used as moduli for geometric tetrahedra. 
As explained below in Remark~\ref{lengths:are:moduli:rem},
for a tetrahedron without ideal vertices, 
the lengths of the internal edges could also be employed,
but of course not in general.

\begin{teo}\label{moduli:teo}
Let $\Delta$ be a partially truncated tetrahedron and let $\Delta\!^{(1)}$ be the
set of edges of $\Delta$. The geometric realizations of $\Delta$ are
parameterized up to isometry 
by the functions $\theta:\Delta\!^{(1)}\to[0,\pi)$ such that:
\begin{itemize}
\item $\theta(e)=0$ if and only if $e$ is length-$0$;
\item For each vertex $v$ of $\Delta$, 
if $e_1,e_2,e_3$ are the edges that emanate from $v$, then
$\theta(e_1)+\theta(e_2)+\theta(e_3)$ is equal to $\pi$ for ideal $v$
and less than $\pi$ for non-ideal $v$.
\end{itemize}
The map $\theta$ corresponding to a geometric realization $\Delta\!^*$ associates to each
non-$0$-length edge $e$ the dihedral angle $\theta(e)$ of $\Delta\!^*$ along $e$.
\end{teo}

\dimostraz 
Our argument follows Fujii's~\cite{fuji}.
Let $\theta:\Delta\!^{(1)}\to[0,\pi)$ be as in the statement.
We fix notation as in Fig.~\ref{tetra:notation:fig} and set $\theta_i=\theta(e_i)$.
Our task is to show that there exists and is unique up to isometry a
geometric realization $\Delta\!^*$ of $\Delta$ with dihedral angles $\theta_i$ along
the $e_i$'s. The idea is to construct the four planes in $\matH^3$ on which the faces 
should lie, and to prove that their configuration is unique 
and determines $\Delta\!^*$ up to isometry.
The plane containing the face with edges $e_i,e_j,e_k$ will be determined
by its circle $C_{ijk}\subset\partial\matH^3$ of points at infinity.
We use the half-space model $\Hhalf^3$ with $\partial\Hhalf^3=\matE^2\cup\{\infty\}$ and
we identify a line $\ell\subset\matE^2$ to the circle $\ell\cup\{\infty\}$.

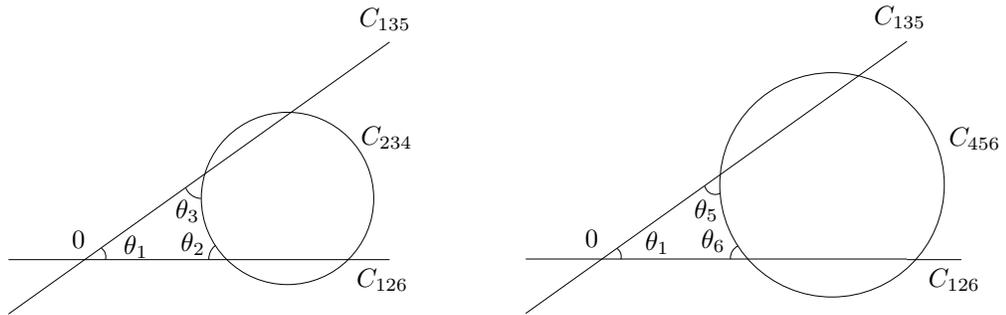
\begin{figure}
\begin{center}
\input{esiste1.pstex_t}
\caption{\small{Lines and circles at prescribed angles with each other.}}
\label{esiste1:fig}
\end{center}
\end{figure}

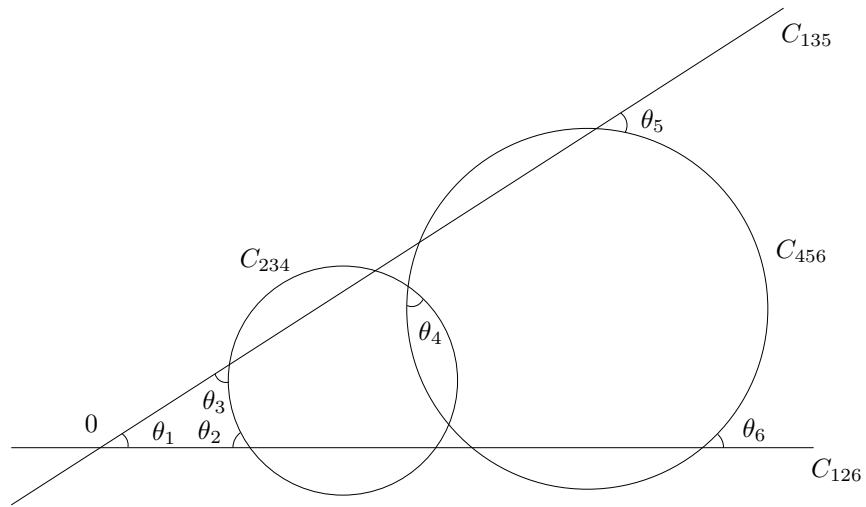
\begin{figure}
\begin{center}
\input{esiste2.pstex_t}
\caption{\small{Trace at infinity of a geometric realization.}}
\label{esiste2:fig}
\end{center}
\end{figure}

We first assume that $\Delta$ has neither
ideal vertices nor length-0 edges, and we show that the 
configuration of the $C_{ijk}$'s exists and is determined by the $\theta_i$'s.
Later we will prove that the configuration determines a unique $\Delta\!^*$.
We choose $C_{126}$ and $C_{135}$ to be lines through $0\in\matE^2$ 
at angle $\theta_1$ with each other.
Conditions $\theta_1+\theta_2+\theta_3<\pi$ and 
$\theta_1+\theta_5+\theta_6<\pi$ imply quite easily that there exist circles
$C_{234}$ and $C_{456}$ as in Fig.~\ref{esiste1:fig}. 
Next, 
we modify $C_{456}$ by a dilation and use conditions
$\theta_3+\theta_4+\theta_5<\pi$ and $\theta_2+\theta_4+\theta_6<\pi$
to show that it can be placed as in Fig.~\ref{esiste2:fig}.
This shows that the configuration of the $C_{ijk}$'s exists.
Its uniqueness follows from
uniqueness up to rotation and dilation of
the configuration of lines and circles of Fig.~\ref{esiste2:fig}.
Now we show that the $C_{ijk}$'s determine $\Delta\!^*$ uniquely.
The truncation plane at the vertex $v_{123}$ (where $e_1,e_2,e_3$ have their
common end) must be orthogonal
to $C_{126},C_{135},C_{234}$. From
Fig.~\ref{esiste1:fig} we see that such a plane exists and is unique,
and similarly for the other three truncation planes.
Moreover the truncation triangles are pairwise disjoint, and
the conclusion follows.

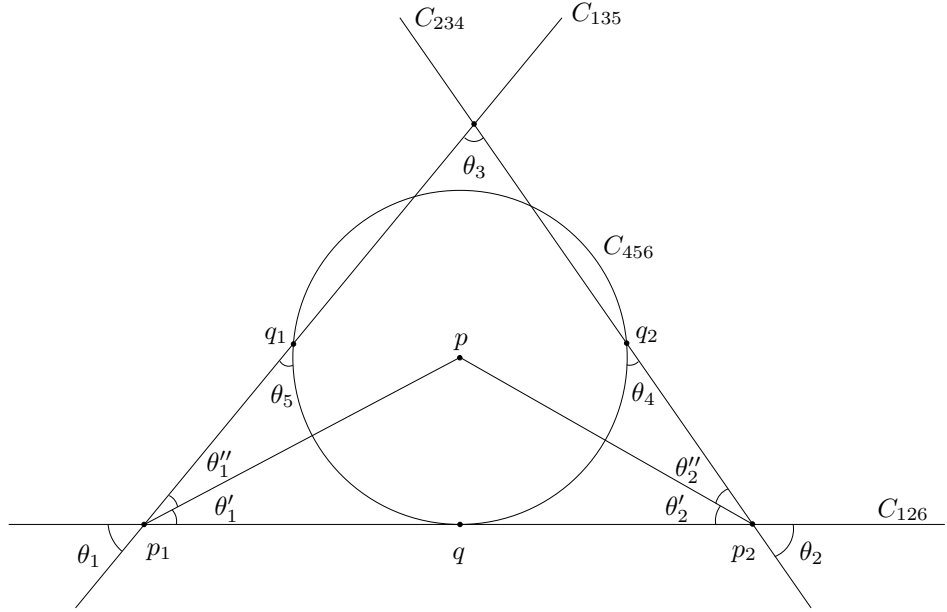
\begin{figure}
\begin{center}
\input{esiste3.pstex_t}
\caption{\small{Configuration at infinity when there is an exceptional hexagon.}}
\label{esiste3:fig}
\end{center}
\end{figure}

The same argument applies with minor changes to
the case where $\Delta$ has ideal vertices or length-$0$ edges.
See for instance Fig.~\ref{esiste3:fig} below for the configuration
to use in case ${\theta_1+\theta_2+\theta_3=\pi}$ and  
${\theta_6=0}$.\finedimo

\begin{cor}
The space of geometric realizations of $\Delta$ up to isometry is a bounded open subset
of a real Euclidean space of dimension $6-n$ where $n$ is the total number
of ideal vertices and length-$0$ edges of $\Delta$.
\end{cor}

\paragraph{From angles to lengths}
Having introduced moduli for geometric tetrahedra, our next task is to 
determine, given a triangulated orientable manifold,
which choices of moduli for the tetrahedra give a global
hyperbolic structure on the manifold. There are two obvious necessary 
conditions (which are often but not always sufficient,~\emph{e.g.}
they are not when all the vertices are ideal). Namely, 
we should have a total dihedral angle of $2\pi$ around each 
non-$0$-length edge of the manifold, and we should be able
to glue the lateral hexagons by isometries.
The first condition is directly expressed in terms of moduli.
To express the second condition recall that
the shape of a hyperbolic right-angled hexagon
is determined by the lengths of a triple of pairwise disjoint
edges. This may seem to suggest that, to ensure consistency,
one only has to compute, in terms of the dihedral angles, either the lengths
of the internal edges or the lengths of the boundary edges.
This is however false when ideal vertices and/or length-$0$ edges are involved, so
we will need to compute both. 

Let us consider a partially truncated tetrahedron $\Delta$ 
with edges $e_1,\ldots,e_6$
as above in Fig.~\ref{tetra:notation:fig}. In the rest of this paragraph we 
fix a geometric realization $\theta$ of $\Delta$ determined by dihedral angles
$\theta_i=\theta(e_i)$ for $i=1,\ldots,6$, and we denote by
$L^\theta$ the length with respect to this realization.
The boundary edges of the lateral hexagons of $\Delta$
correspond to the pairs of distinct non-opposite 
edges $\{e_i,e_j\}$, and will be denoted by $e_{ij}$.
Now $e_{ij}$ disappears towards infinity, so it has length $0$, when the 
common vertex of $e_i$ and $e_j$ is ideal, it is an infinite half-line
when one of $e_i$ or $e_j$ is $0$-length, and it is an infinite line 
when both $e_i$ and $e_j$ are $0$-length. The next result,
that is readily deduced from \cite[The Cosine Rule II, pag. 148]{beardon}
allows to compute the length of $e_{ij}$ when this length is finite. 
We refer to $e_{12}$ with notation as in Fig.~\ref{tetra:notation:fig}.

\begin{prop}\label{boundary:edge:length:prop}
If both $e_1$ and $e_2$ have non-$0$ length then
\begin{equation}\label{boundary:edge:length:formula}
\cosh L^\theta(e_{12})=\frac{\cos \theta_1 \cdot \cos \theta_2 + \cos \theta_3}
{\sin \theta_1 \cdot \sin \theta_2}.
\end{equation}
\end{prop}

Note that this result is correct (and obvious) also when
the common end of $e_1$ and $e_2$ is ideal. 
Turning to the length of an internal edge, we note that the edge 
is an infinite half-line or an infinite line when one or both its ends are ideal.
Otherwise the length is computed 
using~\cite[The Cosine Rule II, pag.~148
and Theorem 7.19.2, pag.~161]{beardon}.
To state the result of the computation we need to introduce certain functions
that will be used again below. 
With notation as in Fig.~\ref{tetra:notation:fig}, and
defining $v_{ijk}$ as the vertex from which the edges
$e_i,e_j,e_k$ emanate, we set:
\begin{eqnarray}
\begin{array}{rcl}
c^\theta(e_1)\!\!&=&\!\!\cos \theta_1\cdot \left( \cos \theta_3\cdot \cos\theta_6
        +\cos \theta_2\cdot \cos\theta_5\right)\\
        & & + \cos \theta_2\cdot \cos\theta_6
        +\cos \theta_3\cdot \cos\theta_5 + \cos\theta_4\cdot \sin^2 \theta_1;
\end{array}\label{ctheta}\\ 
d^\theta(v_{123})\,=\,2\cos \theta_1\cdot \cos\theta_2\cdot \cos\theta_3+ \cos^2 \theta_1
+\cos^2 \theta_2+\cos^2 \theta_3 -1.\label{dtheta}
\end{eqnarray}

\begin{lemma}\label{ideal:vertex:check:lemma}
$d^\theta(v_{123})=0$ if and only if
the vertex $v_{123}$ is ideal.
\end{lemma}

\begin{prop}\label{internal:edge:length:prop}
If $v_{123}$ and $v_{156}$ are both non-ideal then
\begin{equation}\label{internal:edge:length:formula}
\cosh L^\theta(e_1)={c^\theta(e_1)}\,\Big/\,{\sqrt{d^\theta(v_{123})
\cdot d^\theta(v_{156})}}.
\end{equation}
\end{prop}

The next fact will be proved in Section~\ref{computing:section}
using results from Sections~\ref{kojima:section} and~\ref{height:section}.

\begin{prop}\label{lengths:are:moduli}
If $\Delta$ has no ideal vertices then a geometric realization
of $\Delta\!^*$ is determined up to isometry by the lengths of its
internal edges.
\end{prop}

\begin{rem}\label{lengths:are:moduli:rem}
\emph{The previous proposition implies that, when there are no ideal vertices,
one could employ the lengths of the internal edges as moduli.
Besides the loss of generality, this choice is however inadvisable
because of the following drawbacks:
\begin{itemize}
\item In terms of lengths, the restriction that the three dihedral 
angles at each vertex should sum up to less than $\pi$
gets replaced by somewhat more complicated relations. Namely, one should
express boundary lengths in terms of internal lengths, and
then for each vertex impose the triangular inequalities for the three boundary edges
at that vertex;
\item Dihedral angles are needed in any case to ensure consistency,
and to express angles in terms of lengths one should invert
formula~(\ref{internal:edge:length:formula}),
which does not appear to be completely straight-forward.
\end{itemize}}\end{rem}

\paragraph{Exceptional hexagons}
Recall that we are looking for the conditions to ensure
that a gluing of lateral hexagons can be realized by isometries.
By default gluings match ideal vertices to each other and length-0
edges to each other, because these notions are part of the 
initial topological information about a triangulation.

\begin{figure}
\begin{center}
\input{quadri.pstex_t}
\caption{\small{A hexagon with a boundary edge and the
opposite internal edge being length-$0$.}}\label{special:hexagon:fig}
\end{center}\end{figure}
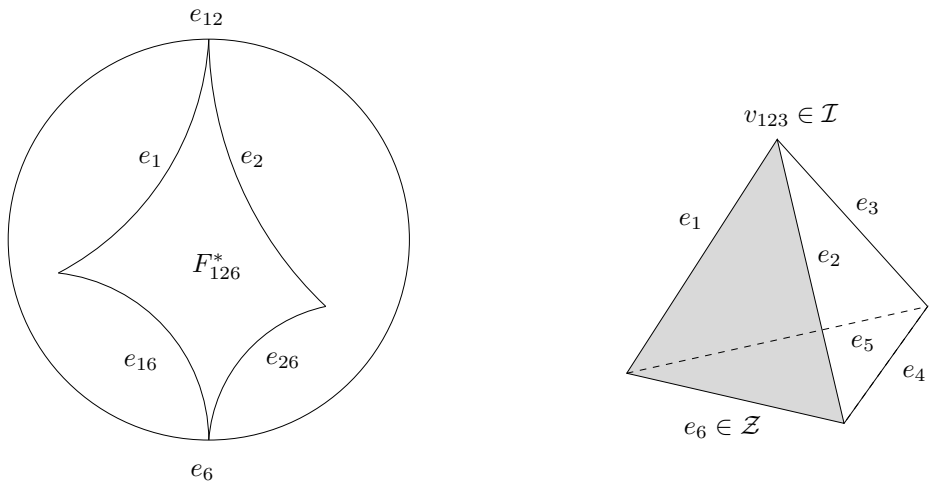

When a pairing glues two \emph{compact} hexagons, \emph{i.e.}~when there
are no ideal vertices or length-0 edges involved,
to make sure that the gluing is an isometric one
we may equivalently require the lengths of the internal edges or
those of the boundary edges to match under the gluing. 
If we actually require \emph{all} lengths to match,
we are guaranteed that the gluing is isometric also
for non-compact hexagons, except in the very special case
where a boundary edge disappears into an ideal vertex,
and the opposite internal edge is length-$0$.
In this case, which is illustrated in Fig.~\ref{special:hexagon:fig},
there is no length at all to match, because two edges have length $0$,
and the other four are infinite half-lines. But 
two hexagons as in Fig.~\ref{special:hexagon:fig} need not
be isometric to each other, as the next discussion shows.

To parameterize the special hexagons we need now to be slightly more
careful about orientation than we have been so far. 
Namely, we choose on the tetrahedra an orientation compatible
with a global orientation of the manifold. As a result
also the lateral hexagons and the truncation triangles
have a fixed orientation, and the gluing maps
reverse the orientation of the hexagons.

So, let us consider an exceptional hexagon $F^*_{126}$ 
as in Fig.~\ref{special:hexagon:fig}, with the boundary edge
$e_{12}$ lying at an ideal vertex, and the opposite internal
edge $e_6$ being length-$0$. Recall that the hexagon is
oriented and embedded in $\matH^3$ by $\theta$.
We consider now the horospheres $O_1$ and
$O_2$ centred at $e_{12}$ and passing through the 
non-ideal ends of $e_1$ and $e_2$ respectively.
We define $\sigma^\theta(F_{126})$ to be
$\pm{\rm dist}(O_1,O_2)$, the sign being positive
if $e_2,e_{12},e_1$
are arranged positively on $\partial F^*_{126}$
and $O_1$ is contained in the horoball bounded
by $O_2$, or if $e_2,e_{12},e_1$ are arranged negatively
on $\partial F^*_{126}$ and $O_2$ is contained in the horoball
bounded by $O_1$, and 
negative otherwise. This definition easily implies the following:

\begin{prop}\label{sigma:works:prop}
Let $F$ and $F'$ be paired exceptional lateral hexagons.
Their pairing can be realized by an isometry if and only
if $\sigma^\theta(F)+\sigma^\theta(F')=0$.
\end{prop}

\begin{prop}\label{compute:sigma:prop}
Let $F_{126}$ be the exceptional hexagon of 
Fig.~\ref{special:hexagon:fig}, oriented so that
$e_2,e_{12},e_1$ are positively arranged on $\partial F^*_{126}$. Then
\begin{equation}\label{compute:sigma:formula}
\sigma^\theta(F_{126})=
\log\frac{\sin \theta_2}{\cos \theta_2+\cos \theta 4}-
\log\frac{\sin \theta_1}{\cos \theta_1+\cos \theta_5}.
\end{equation}
\end{prop}

\dimostraz
We realize $\Delta\!^*$ in $\Hhalf^3$ so that
$v_{123}=\infty$ and denote by $C_{ijk}$ the circle
at infinity of the plane that contains $F^*_{ijk}$.
Let $C_{126}\cap C_{135}=\{p_1,\infty\}$ and 
$C_{126}\cap C_{234}=\{p_2,\infty\}$. Condition $\theta_6=0$ means that
$C_{456}$ is tangent to $C_{126}$. The configuration is then as in 
Fig.~\ref{esiste3:fig}, where we also introduce more notation needed for the proof.
The 
truncation plane relative to $v_{156}$ is now the Euclidean half-sphere
of radius $L(qp_1)$ centred at $p_1$, so the finite end of $e_1$ has coordinates
$(p_1,L(qp_1))$. Similarly the finite end of $e_2$ is $(p_2,L(qp_2))$, 
so $\sigma^\theta(F_{126})=\log L(qp_1)-\log L(qp_2)=\log L(qp_1)/L(qp_2)$.
Of course $L(qp_1)/L(qp_2)=\tan\theta'_2/\tan\theta'_1$. 

Now $\angle(pq_1p_1)=\theta_5+\pi/2$, whence $\sin \angle(pq_1p_1)=\cos\theta_5$,
and the sine theorem yields $L(pq_1)/L(pp_1)=
\sin\theta''_1/\cos\theta_5$.
But $L(pq_1)=L(pq)$, so 
$L(pq_1)/L(pp_1)=\sin\theta'_1$. Equaling the two expressions
of $L(pq_1)/L(pp_1)$ we  get $\sin\theta''_1=\cos\theta_5\cdot\sin\theta'_1$.
On the other hand $\theta''_1=\theta_1-\theta'_1$, so
$\sin\theta''_1=\sin\theta_1\cdot\cos\theta'_1-
\cos\theta_1\cdot\sin\theta'_1$.
Equaling the two expressions of $\sin\theta''_1$ and dividing by $\cos\theta'_1$
we get $\tan\theta'_1=\sin\theta_1/(\cos\theta_1+\cos\theta_5)$.
Similarly $\tan\theta'_2=\sin\theta_2/(\cos\theta_2+\cos\theta_4)$
and the conclusion follows.
\finedimo

\paragraph{Matching around edges}
We have discussed so far the conditions under which the
hyperbolic structure of the geometric tetrahedra matches across
lateral hexagons. As already mentioned, there is another obvious condition
we must impose if we want 
the structure to extend also along the
internal edges of the triangulation.
Namely, let us define for an edge $e$ 
\begin{equation}\label{alpha:def:formula}
\alpha^\theta(e)=\sum\{\theta(e'):\ e' \hbox{ is glued to } e\}.
\end{equation}
Then $\alpha^\theta(e)$ should be $2\pi$ for all $e$.
This condition is actually sufficient
when there are no ideal vertices, but not in general.
The point is that when the geometric tetrahedra are arranged one after 
each other around a non-$0$-length edge $e$, 
the first face of the first tetrahedron
and the second face of the last tetrahedron may overlap without coinciding.
Namely, the isometry which pairs these two faces may be a translation along
$e$ instead of being the identity. Of course the isometry has to be the identity
unless both ends of $e$ are ideal.

We recall now that a horospherical cross-section near an ideal
vertex $v$ of a geometric tetrahedron is a Euclidean triangle
well-defined up to similarity. The tetrahedron being
oriented, this triangle is also oriented, so, once a vertex
of the triangle is fixed, its similarity structure is
determined by a complex 
parameter in  the upper half-plane $\pi_+$.
Choosing a vertex of the triangle amounts to choosing an
edge $e$ ending at $v$, so we have a 
well-defined modulus $z^\theta(e,v)\in\pi_+$ whenever $v$ is ideal
and $e$ ends at $v$. We also define
\begin{equation}\label{big:Z:def:formula}
Z^\theta(e,v)=\prod\{z^\theta(e',v'):
\ (e',v')  \hbox{ is glued to } (e,v)\}.
\end{equation}
The next result is proved just as in the 
purely ideal case (see~\cite{BePe:libro}):

\begin{prop}\label{edge:match:prop}
Assume the structure defined by $\theta$ matches across lateral hexa\-gons,
and let $e$ be an internal edge with both ideal ends.
If $v$ is any one of these ends,
the structure matches 
across $e$ if and only if $\alpha^\theta(e)=2\pi$ 
and $Z^\theta(e,v)=1$.
\end{prop}

\begin{prop}\label{mod:compute:prop}
For an oriented tetrahedron as in Fig.~\ref{tetra:notation:fig}, assume
$v_{123}$ is ideal and $e_1,e_2,e_3$ are positively arranged around $v_{123}$.
Then:
\begin{equation}\label{small:z:def:formula}
z^\theta(e_1,v_{123})=
\frac{\sin \theta_2}{\sin \theta_3}\cdot e^{i\theta_1}.
\end{equation}
\end{prop}

\paragraph{Consistency equations}
We have preferred above to introduce bit after bit
our conditions for the geometric tetrahedra to define
a global structure on the manifold, but now
we collect the relevant information in one precise statement.

\begin{teo}\label{consistency:teo}
Consider an orientable manifold $N$ obtained from a compact $\Nbar$ by removing
all the tori and some annuli contained in $\partial\Nbar$.
Suppose that all the    components of $\partial N$ 
have negative Euler characteristic. Fix a partially truncated 
triangulation of $N$, with the tetrahedra oriented so that the gluings reverse
the induced orientations. Let $\theta$ be a geometric realization
of the tetrahedra in the triangulation. Then 
$\theta$ defines on $N$ a hyperbolic structure
with  geodesic boundary if and only if the following 
conditions hold:
\begin{enumerate}
\item $L^\theta(e)=L^\theta(e')$ for
all pairs $(e,e')$ of matching internal edges;
\item $L^\theta(e)=L^\theta(e')$ for
all pairs $(e,e')$ of matching boundary edges;
\item $\sigma^\theta(F)+\sigma^\theta(F')=0$
for all pairs $(F,F')$ of matching exceptional hexagons;
\item $\alpha^\theta(e)=2\pi$ for all edges $e$;
\item $Z^\theta(e,v)=1$ 
for all edges $e$ with both ideal ends and for both ends $v$  of $e$.
\end{enumerate}
\end{teo}

\begin{rem}\label{after:consistency:rem}
\emph{By Propositions~\ref{boundary:edge:length:prop},~\ref{internal:edge:length:prop},~\ref{compute:sigma:prop}, 
and~\ref{mod:compute:prop}, all the above conditions
can be expressed as analytic equations in terms of the dihedral angles $\theta$. Moreover:
\begin{itemize}
\item If $N$ has no toric end, 
\emph{i.e.}~if $\partial\Nbar$ contains
no tori, then condition (2) is a
consequence of (1), and (3) and (5) are empty. So, to
ensure hyperbolicity of $N$, one may impose (1) and (4) only;
\item If $N$ has no annular end, \emph{i.e.}~if no annuli are removed
from $\partial\Nbar$, then condition (1) is a consequence of (2),
and (3) is empty;
\item When conditions (1), (2), and (3) are in force,
condition (5) may be equivalently imposed at either end of the edge $e$.
\end{itemize}}
\end{rem}

\paragraph{Completeness equations}
The discussion of completeness is very easy. If $N$ is hyperbolic, $N$ is
complete if and only if its double $D(N)$ is. So we should ensure the toric
ends of $D(N)$ to be complete, \emph{i.e.}~their bases to 
have an induced a Euclidean
structure, rather than just a similarity structure. If an end of
$D(N)$ comes from a toric cusp of $N$, completeness is imposed as usual
by requiring the holonomy of the similarity structure on the torus to consist
of translations. If an end of $D(N)$ is the double of an annular end of $N$
then it is actually always complete. To see this, we use
Remarks~\ref{geom:of:tetra:rem} and~\ref{trunc:tria:first:rem},
which imply that the annulus at the basis of the end of $N$ is tiled by a 
(cyclic) row of Euclidean rectangles. So the annulus is Euclidean and its
boundary circles have the same length, whence the double is a Euclidean torus.

We remind the reader that completeness of an end $T\times[0,\infty)$
can be turned into a pair of equations in terms of the moduli of the
triangular horospheric cross-sections, and hence in terms
of the dihedral angles by means of Proposition~\ref{mod:compute:prop}. To do this, 
one first chooses as a basis of $H_1(T;\matZ)$ 
a pair of loops which are simplicial with respect to
the triangulation of $T$ induced by the triangulation of $N$.
Then one notes that the dilation component of the holonomy of a simplicial 
loop $\gamma$ is the product of all moduli $\gamma$ leaves to its left, 
multiplied by $-1$ if $\gamma$ has an odd number of 
vertices~---see 
\emph{e.g.}~\cite{BePe:libro}.

\paragraph{Uniqueness} 
A crucial fact for computational purposes is that a solution of the hyperbolicity
equations, if any, is unique. Before showing this we spell out the rigidity theorem
already mentioned above.

\begin{teo}\label{rigidity:teo}
On any given manifold there exists, up to isometry, at most one
finite-volume and complete hyperbolic structure with geodesic boundary.
\end{teo}

\dimostraz
If $N$ is finite-volume hyperbolic with boundary then $D(N)$ is finite-volume hyperbolic
without boundary, so its structure is unique, and we only need to show that
the embedding of $\partial N$ in $D(N)$ cannot be homotoped away from itself.
Note that $\partial N$, being geodesic, is automatically incompressible.
Uniqueness of embedding for $\partial N$ now easily 
follows by considering non-trivial simple 
loops on $\partial N$ and recalling that every non-trivial 
free-homotopy class of a loop in $D(N)$ has a unique geodesic representative.
\finedimo

\begin{teo}\label{uniqueness:teo}
Given a triangulated manifold $N$ as in Theorem~\ref{consistency:teo},
there exists at most one choice of the geometric realization $\theta$
that turns $N$ into a complete manifold with geodesic boundary.
\end{teo}

\dimostraz
Assume $\theta_0$ and $\theta_1$ yield complete structures on $N$.
By the rigidity theorem
these structures are actually the same, so we can view $\theta_0$ and $\theta_1$
as defining isotopic geometric triangulations of one hyperbolic $N$.
Of course two geometric triangulations are identical if they have the same
edges, so we are left to show that if $(e_t)_{t\in[0,1]}$ is an isotopy
of properly embedded segments, half-lines, or lines, and $e_0,e_1$ are geodesic,
then $e_0=e_1$. Consider first the case of a segment.
If we double $N$ we get for all $t$ 
a closed loop $D(e_t)$ in $D(N)$. The free-homotopy class of $D(e_t)$ is
of course independent of $t$, and it must be non-trivial,
otherwise $e_0$ would lift in $\matH^3$ to a geodesic segment with both ends on
a component of the lifting of $\partial N$, so $e_0$ would actually be contained
in $\partial N$. So all the $D(e_t)$ lift to infinite open lines in $\matH^3$,
and these lines have two well-defined ends. These ends are fixed points
of hyperbolic isometries from a discrete group, and they evolve continuously
along the isotopy. This implies that the ends are actually 
independent of $t$, whence $e_0=e_1$.

The same argument applies when one or both the ends of the edges $e_t$
tend to infinity along a cusp, except that hyperbolic fixed points get replaced
by parabolic fixed points.
\finedimo

\begin{rem}\label{certainly:no:sol:rem}
\emph{As a by-product of the previous argument we deduce that 
if a triangulation contains a boundary-parallel non-$0$-length
edge, then the triangulation is never geometric.}\end{rem}

\section{The Kojima decomposition}\label{kojima:section}

In this section we describe the canonical decomposition due to Kojima~\cite{kojima}
of a hyperbolic 3-manifold with non-empty geodesic boundary, recalling
several details because
we will be using them throughout the rest of the paper.
We omit all proofs addressing the reader to~\cite{kojima}.
All the notation introduced in this section is employed 
extensively later on, so there is basically nothing the reader could 
skip here.

\paragraph{Minkowsky space}
Kojima's construction takes place in 4-dimensional Minkowsky space,
so we start by fixing some notation about it.
We denote by $\Minkos$ the space $\mathbb{R}^4$ with coordinates $x_0,x_1,x_2,x_3$ 
endowed with the Lorentzian
inner product ${\langle x,y\rangle=-x_0 y_0+x_1 y_1+x_2 y_2+x_3 y_3}$. 
We set 
\begin{eqnarray*}
\hm&=&\{x\in\Minkos:\ \langle x,x \rangle=-1,\ x_0>0\},\\
\hp&=&\{x\in\Minkos:\ \langle x,x \rangle=1\},\\
\lp&=&\{x\in\Minkos:\ \langle x,x \rangle=0,\ x_0>0\}.
\end{eqnarray*}
We recall that $\hm$ is the upper sheet of the two-sheeted hyperboloid, and that
$\langle\,\cdot\,,\,\cdot\,\rangle$ restricts to a Riemannian metric
on $\hm$. With this metric,
$\hm$ is the so-called \emph{hyperboloid model} $\Hhyp^3$ of hyperbolic space.
The one-sheeted hyperboloid $\hp$ turns out to have a bijective
correspondence with the set of hyperbolic half-spaces in $\Hhyp^3$. Given
$w\in\hp$, the corresponding half-space, called the \emph{dual} of $w$, is given by
$$\{v\in\hm:\ \langle v,w\rangle \leqslant 0\}.$$
Similarly, the cone $\lp$ of future-oriented light-like vectors of $\Minkos$
corresponds to the set of horospheres in $\Hhyp^3$. The horosphere dual to
$u\in\lp$ is given by 
$$\{v\in\hm:\ \langle v,u\rangle=-1\}.$$ 

\paragraph{Projective model and truncated polyhedra}
Let 
${\pi:\Minkos\setminus\{0\}\to \mathbb{P}(\Minkos)}$
be the canonical projection.
We set ${\pr=\{x\in\Minkos:\ x_0=1\}}$
and note that $\pr$ can be viewed as a subset of $\mathbb{P}(\Minkos)$ via $\pi$. Moreover
$\pr$ is isometric to Euclidean 3-space $\matE^3$, and $\pi$ restricts to
a bijection between $\hm$ and the unit ball of $\pr$.
Giving this ball the metric that turns $\pi$ into an isometry with
$\Hhyp^3$, we get the projective model $\Hproj^3$ of hyperbolic space.
This model is particularly suitable for describing partially truncated tetrahedra, 
as we will now show.

Take in $\pr$ a tetrahedron $\Delta$ 
with distinct vertices outside $\Hproj^3$ or on its boundary, and call \emph{ultra-ideal}
the vertices not lying on $\partial\Hproj^3$ (those on
$\partial\Hproj^3$ are called \emph{ideal} as usual).
Suppose that the interior of every
edge of $\Delta$ meets $\Hproj^3$ or $\partial\Hproj^3$ at least in one point. 
Each ultra-ideal vertex $v$ of $\Delta$ determines a hyperbolic 
plane $H_v$ in $\Hproj^3$ obtained as the intersection of $\Hproj^3$ with the
Lorentzian orthogonal to $v$. Let us now define $\Delta\!^*$ as 
the intersection of $\Delta$ with $\Hproj^3$ truncated by the planes $H_v$
(see Fig.~\ref{tetra:notation:fig} above).
It is easy to see that $\Delta\!^*$ is
a geometric realization of the partially truncated tetrahedron
$\Delta$ in which the ideal vertices are those on
$\partial \Hproj^3$, and the length-$0$ edges 
are those tangent to $\partial \Hproj^3$.
This implies in particular that our notation $\Delta\!^*$ is consistent.

Before proceeding it is worth noting that the hyperbolic plane
$H_v$ described above can also be constructed
by elementary Euclidean geometry on $\pr$. Namely, 
if we take the cone in $\pr\cong\matE^3$ 
with vertex $v$ and tangent to $\partial\Hproj^3\cong\matS^2$, then 
$H_v$ is the Euclidean disc bounded by the circle where the cone intersects 
$\partial\Hproj^3$.

The above definition of $\Delta\!^*$ of course makes sense also for 
convex polyhedra $\Delta$ more complicated than tetrahedra, 
provided $\Delta$ only has ideal and ultra-ideal
vertices, and the interior of every edge of $\Delta$ meets
the closure of $\Hproj^3$. Partially truncated polyhedra of this sort
are the blocks of the Kojima decomposition described in the rest
of this section.

\paragraph{Convex hull of lifted boundary components}
Let $N$ be hyperbolic 
with non-empty geodesic boundary. 
Identifying the universal cover of the double $D(N)$ of $N$ with $\matH^3$,
we can realize the universal cover of $N$ itself as a closed convex
region $\nt$ of $\matH^3$ bounded by a locally finite countable family 
$\mathcal{S}$ of pairwise disjoint planes.
The group of deck transformations of the covering $\nt\to N$, denoted 
henceforth by $\Gamma$, is the stabilizer
of $\nt$ in the group of deck transformations of $\matH^3\to D(N)$.
In the rest of the section we will always use the $\Hproj^3$ model of
$\matH^3$.

Noting that $\nt$ lies on a definite side of each plane $S\in\calS$, we 
consider now the vector of 
$\hp$ dual to the half-plane that contains
$\nt$ and is bounded by $S$. We denote by 
$\calB$ the family of all these duals, and we define
$\calC\subset\Minkos$ as the closure of
$\mathrm{Conv}(\calB)$, where $\mathrm{Conv}(X)$ denotes from now on
the convex hull of a set $X\subset\Minkos$.
Kojima has shown that $\calB$ is a discrete 
subset of $\Minkos$, that $0\notin\calC$ and that 
${\pi(\calC)\supset \Hproj^3}$. In particular $\calC$ has non-empty interior,
and it is $\Gamma$-invariant by construction.
The idea is now to construct a $\Gamma$-equivariant tessellation of 
$\nt$ that projects to a decomposition of $N$ by intersecting $\nt$ with
the projections to $\Hproj^3$ of the 3-faces of $\partial\calC$.
However, it turns out that not all faces should be projected, and
that some faces have non-trivial stabilizer in $\Gamma$, so they 
must be subdivided. To explain the matter in detail we begin with the
following:

\begin{defn}
\emph{Let $X$ be a subset of $\Minkos$ such 
that $0\notin X$. A point $x\in X$ is called 
\emph{almost-visible} with respect to $X$ if the segment $[0,x]$
meets $X$ in $x$ only. The point is called 
\emph{visible} if it is almost-visible
and $\pi(x)\in \Hproj^3$.}
\end{defn}
  
Of course only the faces of $\partial\calC$ containing visible points 
should contribute to the tessellation of $\nt$. These faces are 
also called \emph{visible}.
It turns out that there are two quite different sorts of visible faces, 
called respectively \emph{elliptic} and \emph{parabolic} depending
on whether the restriction of $\langle\,\cdot\,,\,\cdot\,\rangle$ 
to the hyperplane on which
the face lies is positive-definite or positive-semi-definite.

\paragraph{Cut locus and elliptic faces}
The first type of visible faces of $\partial\calC$ 
correspond to the vertices of the
cut-locus of $\nt$ relative to $\bnt$, which we now 
define for an arbitrary manifold $M$.

\begin{defn}
\emph{Let $M$ be hyperbolic with non-empty geodesic boundary.
We define the \emph{cut-locus} $\mathrm{Cut}(M,\partial M)$
of $M$ relative to $\partial M$ as the set of points
of $M$ that admit at least two different shortest paths to $\partial M$.
A point is called a \emph{vertex} 
of the cut-locus if it admits four different shortest paths to $\partial M$
whose initial tangent vectors span the tangent space to $M$ at the point
as an affine space.}
\end{defn}

The next result is implicit in Kojima's work~\cite{kojima}.
A proof is readily deduced from 
Proposition~\ref{safe:height:prop} shown below (using discreteness, which
is easy to establish).

\begin{prop}\label{cutlocus:prop}
$\mathrm{Cut}(N,\partial N)$ has finitely many vertices. A point of 
$\nt$ is a vertex of $\mathrm{Cut}(\nt,\bnt)$ if and only
it projects in $N$ to a vertex of $\mathrm{Cut}(N,\bn)$.
\end{prop}

For every vertex $v$ of $\mathrm{Cut}(\nt,\bnt)$ 
we define now $\calB(v)\subset\calB$
as the set of dual vectors to the hyperplanes in $\calS$
having shortest distance from $v$.
Since in $\matH^3$ there is a unique shortest path joining 
a given point to a given plane, each $\calB(v)$ contains
at least four vectors. The next result describes the visible elliptic
faces of $\partial\calC$.

\begin{prop}\label{ell:face:prop}
For every vertex $v$ of $\mathrm{Cut}(\nt,\bnt)$ 
there exists a unique elliptic hyperplane $E(v)\subset\Minkos$ such that
${E(v)\cap\calB=\calB(v)}$. Moreover $E(v)$ is a support hyperplane for $\calC$, and:
\begin{enumerate}
\item The points in
$E(v)\cap\calC$ are almost-visible with respect to $\calC$;
\item $E(v)\cap\calC$ is a $3$-dimensional compact polyhedron whose 
set of vertices is $\calB(v)$;
\item The stabilizer of $E(v)\cap\calC$ in $\Gamma$ is trivial.
\end{enumerate}
\end{prop}

\paragraph{Parabolic faces and subdivision}
Besides those corresponding to vertices of the cut-locus, 
$\partial\calC$ has visible
faces coming from toric cusps of $N$.
Let us denote by $\nt_{\infty}\subset\partial\Hproj^3$ the set of 
points at infinity of $\nt$. A point
$q\in\nt_{\infty}$ is said to \emph{generate a toric cusp} 
if it is fixed under a $\matZ\oplus\matZ$ subgroup $\Gamma_q$ 
of parabolic elements of $\Gamma$,      
\emph{i.e.}~if there is 
a horoball centred at $q$ that projects to a toric cusp of $N$. 

For $\tilde{q}$ in $\lp$ and $t\in\matR$ we consider now the affine parabolic hyperplane
$$F(\tilde{q},t)=\{v\in\Minkos:\ \langle v,\tilde{q}\rangle=-t\}.$$

\begin{prop}\label{parab:face:prop}
Let $q\in\nt_\infty$ generate a toric cusp of $N$, and
take $\tilde{q}\in\lp$ such that $\pi(\tilde{q})=q$.
Then there exists a unique $t(\tilde{q})\in\mathbb{R}$ such that 
$F(\tilde{q},t(\tilde{q}))$
is a support hyperplane for $\calC$. 
Moreover $F(\tilde{q},t(\tilde{q}))$ depends only
on $q$, not on $\tilde{q}$, and,
setting ${F(q)=F(\tilde{q},t(\tilde{q}))}$, we have that:
\begin{enumerate}
\item\label{pos:tq:point} $t(\tilde{q})>0$ and 
$\calC\subseteq\bigcup_{t\geqslant t(\tilde{q})} F(\tilde{q},t)$;
\item The points in $F(q)\cap\calC$ are almost-visible with respect to $\calC$;
\item $F(q)\cap\calB$ is infinite and
$F(q)\cap\calC$ is a non-compact $3$-dimensional polyhedron whose set of vertices is
$F(q)\cap\calB$;
\item The $2$-dimensional faces of $F(q)\cap\calC$ are compact;
\item\label{non:triv:stab:point} The stabilizer in $\Gamma$ of 
$F(q)\cap\calC$ coincides with the $\matZ\oplus\matZ$ stabilizer $\Gamma_q$
of $q$.
\end{enumerate}
\end{prop}      

\noindent Point~(\ref{non:triv:stab:point}) of this
proposition shows that 
$F(q)\cap\calC$ must be subdivided before projecting to $\Hproj^3$
and intersecting with $\nt$. Continuing with the same notation,
for every 2-dimensional face $W$
of $F(q)\cap\calC$ we define now $W^{(\tilde{q})}$ as the 
cone based on $W$ with vertex in $\tilde{q}$. 
Note that $W^{(\tilde{q})}$ depends on $\tilde{q}$, and it
meets $\calC$ in $W$ only. However one easily sees that 
$W^{(q)}=\pi(W^{(\tilde{q})})$ actually
does not depend on $\tilde{q}$.
Moreover the family of all $W^{(q)}$'s, as $W$ varies
in the 2-faces of $F(q)\cap\calC$, gives a 
$\Gamma_q$-equivariant
tessellation of $\{\tilde{q}\}\cup\pi(F(q)\cap\calC)$ in which every 
polyhedron has trivial stabilizer.

\paragraph{Canonical decomposition}
We begin with the following fact: 

\begin{prop}\label{all:faces:prop}
The faces of $\calC$ described in 
Propositions~\ref{ell:face:prop} and~\ref{parab:face:prop} 
contain all the visible points of $\calC$.
\end{prop}

In addition to this, one can easily show that if $\Delta$ is a visible face of 
$\calC$ then the (partial) truncation $\pi(\Delta)^*$,
defined earlier in this section,
of $\pi(\Delta)\subset\pr$ is obtained by intersection with $\nt$.

We are now ready to summarize the construction. Let us denote by
$\calV$ the set of vertices of ${\rm Cut}(\nt,\bnt)$, and
by $Q$ the family of all points of $\nt_\infty$ that
generate toric cusps of $N$. For each $q\in Q$ we fix an arbitrary
$\tilde{q}\in\lp$ such that $\pi(\tilde{q})=q$, and we denote by 
$\widetilde{Q}$ the family of all such $\tilde{q}$'s.
We define $\kojiMinkos(\widetilde{Q})$ to be the family of all visible
elliptic faces $E(v)\cap\calC$ and all faces $W^{(\tilde{q})}$ obtained
from the visible parabolic faces.
Here $v$ varies in $\calV$, $\tilde{q}$ varies in $\widetilde{Q}$, and
$W$ varies in the 2-dimensional faces of
$F(\pi(\tilde{q}))\cap\calC$. We also denote by $\kojiproj$ the family of polyhedra
obtained by projecting to $\pr$ the elements of $\kojiMinkos(\widetilde{Q})$, and
by $\kojintstar$ the family obtained by intersecting with $\nt$ 
(or, equivalently, truncating) the elements of $\kojiproj$.
We know that indeed $\kojiproj$ and $\kojintstar$ are independent of 
$\widetilde{Q}$, and we have:

\begin{teo}\label{Koji:main:teo}
$\kojintstar$ is a $\Gamma$-equivariant decomposition of $\nt$
into partially truncated polyhedra. Every element of $\kojintstar$
has trivial stabilizer in $\Gamma$, so $\kojintstar$
projects to a canonical and finite decomposition $\kojiNstar$ of $N$
into partially truncated polyhedra.
\end{teo}

\begin{rem}
\emph{By construction, each polyhedron in $\kojiproj$ has at most one
ideal vertex. All other vertices are ultra-ideal.}
\end{rem}

\section{Choice of heights}\label{height:section}
The general strategy to decide whether a given decomposition
of a hyperbolic $N$ is the canonical one is 
as in~\cite{weeks:tilt}, \emph{i.e.}~to lift the decomposition
first to $\Hproj^3$ and then to $\Minkos$, and to make sure that
the resulting polyhedra bound a convex set. In the setting of the construction
described in the previous section, the lifting of a 
polyhedron to $\Minkos$ can be performed directly when $N$
does not have toric cusps, because one only has to select which points of
$\calB$ are the vertices of the lifting, and
$\calB$ depends on $N$ only. When there are toric cusps, 
however, we have an arbitrariness for the liftings of the ideal
vertices, corresponding to the arbitrariness of the choice of $\widetilde{Q}$.
In addition, for some choices of $\widetilde{Q}$ it can happen that
$\kojiMinkos(\widetilde{Q})$ does \emph{not} bound a convex set. The case
of toric cusps must therefore be discussed with some care.

The reader willing to catch just the main points of the section 
and proceed to the sequel could
devote attention only to Definition~\ref{Qtilde:def},
formula~(\ref{C'(O):def:formula}), 
Definition~\ref{convex:angle:def}, and 
Propositions~\ref{new:conditions:for:conv:prop},~\ref{riconoscimento:prop},~\ref{safe:height:prop}, and~\ref{verysafe:height:prop}.

\paragraph{Horospherical cross-sections}
We first show how to 
reduce the choice of $\widetilde{Q}$ to a choice that is 
more intrinsic, but still arbitrary at this stage. 
We fix $N$ hyperbolic with $\partial N\ne\emptyset$ for the rest of the section.

\begin{defn}\label{Qtilde:def}
\emph{Let $N$ be hyperbolic.
A family $\calO$ of disjoint tori in $N$ that lift to horospheres in $\nt$
and bound disjoint cusps of $N$ will be called a
\emph{horospherical cross-section} of $N$.
Given such an $\calO$,
we define $\widetilde{Q}(\calO)$ as the lifting of $Q$
where a point $q$ is lifted to the only $\tilde{q}\in\lp$
such that $\pi(\tilde{q})=q$ and the dual to $\tilde{q}$ projects in
$N$ to a component of $\calO$.
We also define $\kojiMinkos(\calO)$ as $\kojiMinkos(\widetilde{Q}(\calO))$.}
\end{defn}

\begin{rem}
\emph{The cross-section $\calO$ is determined
by a sufficiently small positive number assigned to each toric cusp,
namely the volume of the region between the cross-section at the cusp 
and infinity. Insisting all volumes to be equal to each other, one may also
determine $\calO$ by a single number. This number is naturally interpreted as
the (inverse of the) \emph{height} at which the cross-section should be taken ---the
smaller the volume, the higher the cross-section.}
\end{rem}

\paragraph{Augmented convex hull}
Before proceeding, we need to define a set of 
which $\kojiMinkos(\calO)$ gives the visible boundary. Namely, we set:
\begin{equation}\label{C'(O):def:formula}
\calC'(\calO)=\calC\cup\bigcup\Big\{
\mathrm{Conv}\Big(\{\tilde{q}\}\cup (F(q)\cap\calC)\Big):\ 
q\in Q,\ \tilde{q}=\widetilde{Q}(\calO)\cap\pi^{-1}(q)\Big\}.
\end{equation}
The fact that indeed the visible boundary of $\calC'(\calO)$ is $\kojiMinkos(\calO)$
easily follows from point~(\ref{pos:tq:point}) of
Proposition~\ref{parab:face:prop}.

\begin{lemma}\label{calC:primo:is:cone:lem}
Suppose $x$ is an almost-visible point of $\calC'(\calO)$ and 
$\pi(x)\in\nt$.
Then $\alpha\cdot x\in\calC'(\calO)$ for every $\alpha\geqslant 1$. 
\end{lemma}

\dimostraz
The statement for $\calC'(\calO)$ easily follows from the corresponding
statement for $\calC$. In addition,
since $\calC$ is closed, it is sufficient to prove the statement
when $\pi(x)$ lies in the interior of $\nt$. Convexity of $\nt$ then implies that
$\pi(x)$ lies in the interior of the convex hull of a finite number
of points of $\nt_\infty$. Moreover we know from~\cite{kojima} 
that ${\nt_{\infty}=\overline{\pi(\calB)}
\cap\partial\Hproj^3}$, and it is easy to deduce that
$\pi(x)$ lies in the convex hull
of infinitely many finite and pairwise disjoint 
subsets of $\pi(\calB)$. But $\calB$ is discrete 
in $\Minkos$, so there is a sequence
$\alpha_n$ of reals diverging to $+\infty$ 
such that $\alpha_n\cdot x\in\calC$ for all $n$.
Convexity of $\calC$ now implies the desired conclusion.
\finedimo

\begin{defn}\label{convex:angle:def}
\emph{Let $P_1,P_2$ be finite 3-dimensional convex
polyhedra in $\Minkos\setminus\{0\}$ 
which project injectively to $\matP(\Minkos)$.
Assume $P_1$ and $P_2$ share a 2-face $F$, and
for $i=1,2$ let $H_i$ be the half-hyperplane in $\Minkos$ 
such that $H_i\supset P_i$ and $\partial H_i\supset F$.
We say $P_1$ and $P_2$ form a \emph{convex angle} at $F$
if the connected component not containing $0$
of $\Minkos\setminus (H_1\cup H_2)$  
is convex. The angle 
is called \emph{strictly convex} if it is convex and $H_1\cup H_2$ 
is not a hyperplane.}
\end{defn}

\begin{lemma}\label{conditions:for:conv:lem}
The following facts are equivalent:
\begin{enumerate}
\item
$\calC'(\calO)\cap\pi^{-1}(\nt)$ is convex; 
\item The $3$-faces of $\calC'(\calO)$ meeting $\pi^{-1}(\nt)$ 
form convex angles with each other.
\end{enumerate}
\end{lemma}

\dimostraz
Implication (1) $\Rightarrow$ (2) is clear. 
To show the opposite implication, we take $x,y\in\calC'(\calO)$ such that
$\pi(x),\pi(y)\in\nt$. Using the fact that angles are convex,
we want to show that $[x,y]\subset\calC'(\calO)$. Since $\calC'(\calO)$ is closed
we can also slightly perturb $x$ and $y$, so we assume that $\pi(x)\ne\pi(y)$
and $[\pi(x),\pi(y)]$ does not meet the 1-skeleton of the decomposition
$\kojintstar$ of $\nt$. So ${[\pi(x),\pi(y)]}$ meets finitely many $2$-faces,
and the picture on $\pi^{-1}([x,y])$ is as in Fig.~\ref{convessa:fig}
\begin{figure}[t!]
\begin{center}
\input{new_convessa.pstex_t}
\caption{\small{Convexity of $\calC'(\calO)\cap\pi^{-1}(\nt)$.}}
\label{convessa:fig}
\end{center}
\end{figure}
by the assumption on convexity of angles. The conclusion now follows from
Lemma~\ref{calC:primo:is:cone:lem}.
\finedimo

The previous result easily implies that if 
$\calC'(\calO)\cap\pi^{-1}(\nt)$ is convex for some choice of ``heights'' $\calO$,
then it is also convex for any ``higher'' choice.
We will now show that such a choice is possible. To
this end we will need the following easy result (points (1) and (2) of which are
actually
only used in Section~\ref{tilt:section} below):

\begin{lemma}\label{distances:lem} 
\begin{enumerate}
\item Let $H_1$ and $H_2$ be disjoint planes in $\Hproj^3$. 
Choose disjoint half-spaces bounded by these planes and let $w_1,w_2\in\hp$ be the 
duals to these half-spaces. Then
$\cosh d(H_1,H_2)=-\langle w_1,w_2\rangle;$
\item Let $W_1$ and $W_2$ be half-spaces in $\Hproj^3$ such that 
${\partial W_1\cap \partial W_2\neq\emptyset}$ 
and let $\theta$ be the dihedral angle they determine; let
${w_1,w_2\in\hp}$ be the duals to $W_1$ and $W_2$. Then ${\cos
\theta=-\langle w_1,w_2\rangle}$;
\item Let $O$ and $H$ be respectively
a horosphere and a plane in $\Hproj^3$. Assume $O\cap H=\emptyset$,
let $u\in\lp$ be the dual to $O$ and 
$w\in\hp$ be the dual to the subspace bounded by $H$
that contains $O$. Then
$\exp {d(O,H)}=-\langle u,w\rangle;$
\item Let $O_1$ and $O_2$ be disjoint horospheres with different centres 
in $\Hproj^3$, and let $u_1,u_2\in\lp$ be their dual vectors. Then
$\exp {d(O_1,O_2)}=
-\frac{1}{2}\langle u_1,u_2\rangle.$
\end{enumerate}
\end{lemma}

\begin{prop}\label{new:conditions:for:conv:prop}
Let $N$ be hyperbolic with universal covering $\Hproj^3\supset\nt\to N$.
Let $\calO$ be a horospherical cross-section for $N$ 
such that:
\begin{enumerate}
\item For any distinct components
$O_1$ and $O_2$ of the lifting of $\calO$ we have
$$\exp d(O_1,\bnt)+\exp d(O_2,\bnt)<2 \cdot\exp d(O_1,O_2);$$
\item The toric cusps in $N$ determined by $\calO$ do not
contain vertices of $\mathrm{Cut}(N,\bn)$,
and for any such vertex $u$ we have 
$\sinh d(u,\bn)<\exp d(u,\calO)$.
\end{enumerate}
Then the visible $3$-faces of $\calC'(\calO)$
form strictly convex angles with each other. In particular, 
$\calC'(\calO)\cap\pi^{-1}(\nt)$ is convex.
\end{prop}

\dimostraz
Recall that we have two types of faces, the elliptic ones
$E(v)\cap\calC$ where $v$ is a vertex of $\mathrm{Cut}(\nt,\bnt)$, 
and the $W^{(\tilde{q})}$'s
where $\tilde{q}\in\widetilde{Q}(\calO)$ and 
$W$ is a 2-face of the parabolic face $F(\pi(\tilde{q}))\cap\calC$.
By construction elliptic faces form strictly convex angles with each other, 
and the same happens for faces $W_1^{(\tilde{q})}$ and $W_2^{(\tilde{q})}$
relative to the same $\tilde{q}$. So we have two cases to discuss.

We begin with the case of two faces of the form 
$W_1^{(\tilde{q}_1)}$ and $W_2^{(\tilde{q}_2)}$,
which of course can only have a common 2-face when $W_1=W_2=:W$.
Knowing that the horospheres $O_1$ and $O_2$
dual to $\tilde{q}_1$ and $\tilde{q}_2$
are disjoint, it is not hard to see that there exists a real number $a>1$ and 
an isometry of $\Minkos$ that carries $\tilde{q}_1$ and $\tilde{q}_2$
to the points $(a,a,0,0)$ and $(a,-a,0,0)$ respectively.
Using Lemma~\ref{distances:lem}~(4) it is also easy to see that
$a$ is intrinsically interpreted as the square root of $\exp d(O_1,O_2)$.

Of course convexity and all the relevant quantities are preserved under
isometry, so we can just assume
$\tilde{q}_1=(a,a,0,0)$ and $\tilde{q}_2=(a,-a,0,0)$.
Recall now that for $i=1,2$ we have in $\Minkos$ the hyperplanes
$F(\pi(\tilde{q}_i))=F(\tilde{q}_i,t(\tilde{q}_i))$, and
the face $W$ at which we must
prove convexity lies in the intersection of these hyperplanes.
Moreover Lemma~\ref{distances:lem}~(3) implies that 
$t(\tilde{q}_i)$ has the intrinsic meaning of 
$\exp d(O_i,\bnt)$. Now the 2-plane on which
the face $W$ lies is the following one:
$$\Bigl\{x\in\Minkos:\ 
x_0=\big(
t(\tilde{q}_1)+t(\tilde{q}_2)\big)\,/\, 2a,
\ x_1=\big(
-t(\tilde{q}_1)+t(\tilde{q}_2)\big)\,/\,2a
\Bigr\}.$$
Therefore the angle at $W$ is strictly convex if and only if
$(1/2a)\cdot\left(t(\tilde{q}_1)+t(\tilde{q}_2)\right)<a,$
and this inequality holds by  the first assumption 
of the statement and
the intrinsic interpretation of the $t(\tilde{q}_i)$'s and $a$.

Turning to the angle between a face $W^{(\tilde{q})}$
and a face $E(v)\cap\calC$, we denote by $O$ the horosphere dual 
to $\tilde{q}$, and note that the horoball bounded by $O$ cannot
contain $v$ by assumption. Using this fact it is not hard to show
that up to isometry in $\Minkos$ we can assume that $v=(1,0,0,0)$
and $\tilde{q}=(b,b,0,0)$ for some $b>1$. Again using the above
lemma, one computes $b$ to be $\exp d(v,O)$. Moreover
$E(v)$ is the hyperplane of equation $x_0=a$, with 
$a=\sinh d(v,\bnt)$. Now the angle at $W$ is strictly convex
if and only if $a<b$, and the conclusion follows.
\finedimo

\begin{cor}\label{esisteconv:cor}
$\calO$ can be chosen so that both the above conditions hold.
\end{cor}

\dimostraz
Finiteness of the number of vertices of $\mathrm{Cut}(N,\bn)$ 
readily implies that the second condition of the proposition is satisfied
for some suitably high $\calO$. Now let
${\alpha=\max\{d(O,\bn):\ O\in\calO\}}$, and redefine $\calO$ by pushing up 
each horospherical cross-section at distance $\alpha$ from its previous position.
A straight-forward computation then implies that also the first
condition of the proposition is fulfilled.
\finedimo

\paragraph{Recognizing the canonical decomposition}
Assume now that we have some decomposition $\decoNstar$ of $N$
into partially truncated polyhedra, and let us ask ourselves whether
$\decoNstar$ coincides with Kojima's canonical $\kojiNstar$.
We fix as above the universal cover $\nt\to N$ with $\nt\subset\Hproj^3$
and a horospherical cross-section $\calO$ of $N$. Now we can first lift
$\decoNstar$ to a decomposition $\decontstar$ of $\nt$, then
we can consider the associated non-truncated family of polyhedra
$\decoproj$ in $\pr$, and finally we can lift $\decoproj$ to a
$\decoMinkos(\calO)$ in $\Minkos$ by lifting each ultra-ideal vertex 
to its representative in $\calB$ and each ideal vertex $q\in Q$ to 
$\widetilde{Q}(\calO)\cap\pi^{-1}(q)$.

\begin{prop}\label{riconoscimento:prop}
Let $\decoNstar$ be a geometric decomposition of a hyperbolic $N$ and let 
$\calO$ be a horospherical cross-section for $N$ as in Proposition~\ref{new:conditions:for:conv:prop}.
Let $\decoMinkos(\calO)$ be the lifting of $\decoNstar$ 
to $\Minkos$ just described.
Then:
\begin{itemize}
\item $\decoNstar$ is the canonical 
decomposition $\kojiNstar$ if and only if
the $3$-faces of $\decoMinkos(\calO)$ form strictly convex angles with each other;
\item $\decoNstar$ is a subdivision of the canonical 
decomposition $\kojiNstar$ if and only if
the $3$-faces of $\decoMinkos(\calO)$ form convex angles with each other.
\end{itemize}
\end{prop}

\dimostraz
The ``only if'' assertions are obvious. To prove the ``if''
assertions, assume the angles are convex. Define
$Y\subset\Minkos$ as the union of the polyhedra of $\decoMinkos(\calO)$
and set $X=\{\alpha\cdot y:\ \alpha\geqslant 1,\ y\in Y\}$. 
Using convexity of angles, 
the same argument
given for Lemma~\ref{conditions:for:conv:lem} shows 
that $X$ is convex. Now $X\supset Y\supset \calB\cup\widetilde{Q}(\calO)$,
so $X\supset\calC'(\calO)$ by convexity. 
Moreover, the polyhedra of $\decoMinkos(\calO)$
have vertices in $\calB\cup\widetilde{Q}(\calO)$,
so $Y\subset\calC'(\calO)$.
Since $Y$ is the visible boundary of $X$, we can conclude
that it is also the visible boundary of $\calC'(\calO)$, so
$\kojiMinkos(\calO)$ is the natural decomposition of $Y$ into faces.
Our given decomposition $\decoMinkos(\calO)$ of $Y$ then coincides with
$\kojiMinkos(\calO)$ when all its angles are strictly convex.
If there are flat angles, however, $\decoMinkos(\calO)$ is only
a subdivision of $\kojiMinkos(\calO)$.
\finedimo

We have now some remarks about how to
apply Proposition~\ref{riconoscimento:prop} in practice.
To begin, note that to ensure (strict) convexity at all the 
infinitely many $2$-faces of $\decoMinkos(\calO)$, it is actually
sufficient to check it for one lifting of each $2$-face of $\decoN$,
and there are finitely many of such faces. However, two serious issues remain.
First, we need an effective method to check convexity, and we will 
provide one in Section~\ref{tilt:section}.
Second, we need a way to determine $\calO$ using the geometry of $\decoNstar$ only. 
A partial step in this direction is discussed in
the rest of this section, and the conclusion is given
(when $\decoNstar$ is actually a triangulation) in Section~\ref{computing:section}.

\paragraph{Intrinsic computation of height}
Now we show how to find a horospherical cross-section
as in Proposition~\ref{new:conditions:for:conv:prop}
in terms of the geometry of $N$. 

\begin{prop}\label{safe:height:prop}
Let $N$ be hyperbolic and realize the universal cover $\nt\to N$ with $\nt\subset\Hhalf^3$ 
so that $\infty$ generates a toric cusp of $N$. 
Let $\Gamma_\infty$ be the $\matZ\oplus\matZ$ stabilizer
of $\infty$ in the group of deck transformations of $\nt\to N$. 
Let $r_1>r_2$ be the first and second largest radii of 
the components of $\bnt$, viewed as Euclidean half-spheres in $\matC\times(0,\infty)$.
Let $d$ be the diameter of a fundamental domain for the action of $\Gamma_\infty$
on $\matC\times\{0\}$, and define
$k(r_1,r_2,d)=\sqrt{3}\cdot(r_1^2+d^2/4)/(r_1-r_2).$
Consider the 
horosphere $\widetilde{O}=\matC\times\{k(r_1,r_2,d)\}\subset\Hhalf^3$ and let 
$O$ be the projection of $\widetilde{O}$ in $N$. Then:
\begin{enumerate}
\item\label{toric:section:point}
$O$ is an embedded toric cross-section of the cusp, and it is disjoint from
any other arbitrarily chosen embedded toric cross-section at any of the other cusps;
\item\label{no:vertices:point}
The cusp bounded by $O$ does not contain vertices of $\mathrm{Cut}(N,\bn)$;
\item\label{far:from:verts:point}
If $u$ is a vertex of ${\mathrm{Cut}(N,\bn)}$ we have 
$\sinh d(u,\bn)<\exp d(u,O)$.
\end{enumerate}\end{prop}

\dimostraz
The constants $r_1,r_2,d$ are fixed, so we set $k=k(r_1,r_2,d)$.
Since $k>r_1$, we see that $\widetilde{O}\cap\partial\nt=\emptyset$.
So both assertions of point~(\ref{toric:section:point})
follow from the following claim:
\emph{If $\widetilde{O}'$ is a horosphere centred at some point  
of ${\partial\Hhalf^3\setminus\{\infty\}}$ 
and ${\widetilde{O}'\cap\partial\nt=\emptyset}$ then 
$\widetilde{O}'\cap\widetilde{O}=\emptyset$.} To prove this claim, let 
$\widetilde{O}'$ be a Euclidean sphere of radius $x$ tangent to
$\matC\times\{0\}$ at a point $z$. We must show that $2x<k$.
Now within distance $d/2$ from $z$ there exists the centre $w$ of a component 
of $\partial\nt$ of Euclidean radius $r_1$. Knowing that 
$\widetilde{O}'$ and this component are disjoint, we deduce that ${x<(d^2/4-
r_1^2)/(2\cdot r_1)}$ whence
the conclusion at once.

To prove points~(\ref{no:vertices:point}) and~(\ref{far:from:verts:point})
we set $h=\sqrt{(r_1^2+d^2/4)/(1-r_2/r_1)}$
and claim the following:
\emph{${\matC\times[h,\infty)}$ does not contain any vertex
of ${\mathrm{Cut}(\nt,\bnt)}$.} 
Using Proposition~\ref{cutlocus:prop} and the easy fact that
${k\geqslant h}$, our claim readily implies 
point~(\ref{no:vertices:point}), and it will be used below for 
point~(\ref{far:from:verts:point}).

To prove the claim, let 
$v$ a vertex of ${\mathrm{Cut}(\nt,\bnt)}$. We first show that the components 
of $\bnt$ nearest to $v$ cannot all have the same Euclidean radius.
If this were the case, turning to the setting of Kojima's construction
in $\Minkos$ and using Lemma~\ref{distances:lem}~(3), we would deduce that
$\calB(v)$ is contained in a parabolic affine hyperplane of $\Minkos$
(Lorentz orthogonal to the dual of $\widetilde{O}$ in $\lp$), against
Proposition~\ref{ell:face:prop}. 

\begin{figure}
\begin{center}
\input{altezza.pstex_t}
\caption{\small{Notation for the proof of 
Proposition~\ref{safe:height:prop} (2,3).}}
\label{altezza:fig}
\end{center}
\end{figure}
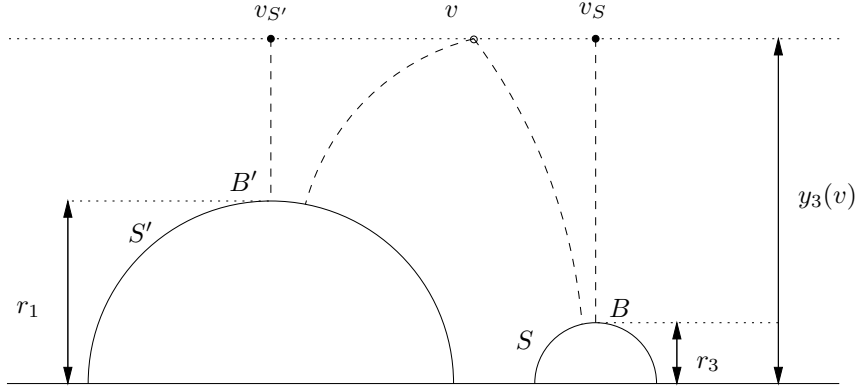

Since the components 
of $\bnt$ nearest to $v$ do not all have the same radius, one of them, say $S$,
has some radius ${r_3\leqslant r_2}$.
By definition of $d$ there exists
another component $S'$ of $\bnt$ with radius
$r_1$ such that, with notation as in Fig.~\ref{altezza:fig},
the Euclidean 
distance between $v$ and $v_{S'}$ is at most $d/2$. 
Recall now that
if $y$ is the $(0,\infty)$-coordinate on $\Hhalf^3=\matC\times(0,\infty)$
and $\|\cdot\|$ is the Euclidean norm,
the hyperbolic distance between ${p,q\in\Hhalf^3}$ satisfies
\begin{equation}\label{dist:form}
\cosh d(p,q)=1+{\|p-q\|^2}\,\big/\,{2y(p)y(q)}.
\end{equation}
Now 
$d(v_S,B)=d(v_S,S)\leqslant d(v,S)$, and we know that 
$d(v,S)\leqslant d(v,S')\leqslant d(v,B')$, so $d(v_S,B)\leqslant d(v,B')$.
Using~(\ref{dist:form}) we easily deduce that
\[
y(v)^2
\leqslant\frac{r_1\cdot r_3+d^2/4\cdot r_3/r_1-r_3^2}{1-r_3/r_1}
\leqslant\frac{r_1^2+d^2/4}{1-r_3/r_1}
\leqslant\frac{r_1^2+d^2/4}{1-r_2/r_1}=h^2.
\] 
Our claim, and hence point~(\ref{no:vertices:point}), are proved. For
point~(\ref{far:from:verts:point}) it is
sufficient to show that for any vertex $v$ of 
$\mathrm{Cut}(\nt,\bnt)$ we have ${\sinh d(v,\bnt)<
\exp d(v,\widetilde{O})}$.
Choose again $S'$ as in Fig.~\ref{altezza:fig}.
The claim just proved and the easy fact that $h>r_1$ show that
$\|v-B'\|^2<h^2+d^2/4$. But $h>d/2$, so $\|v-B'\|^2<2h^2$.
This inequality and~(\ref{dist:form}) easily imply that
$$\sinh^2 d(v,\bnt)\leqslant\sinh^2 d(v,B')
< {\big(h^4+2\cdot h^2\cdot y(v)\cdot r_1\big)}\,\big/\,{\big(y(v)^2\cdot r_1^2\big)}.$$
As noted above we have $r_1<h$. Moreover $y(v)<h$ by the claim shown, so 
$$\sinh^2 d(v,\bnt)<
{(3\cdot h^4)}\,\big/{(y(v)^2 \cdot r_1^2)}=
{k^2}\,\big/\,{y(v)^2}=\exp^2 d(v,\widetilde{\calO}),$$
and the proof is complete. 
\finedimo

\begin{prop}\label{verysafe:height:prop}
Let a hyperbolic $N$ have $n$ toric cusps, and for $i=1,\ldots,n$ choose a realization
$\nt_i$ of $\nt$ in $\Hhalf^3$ so that
$\infty$ generates the $i$-th cusp.
Let $r_1^{(i)},r_2^{(i)},d^{(i)}$ be
the constants relative to $\nt_i$ as in Proposition~\ref{safe:height:prop}. Set
$$\lambda=\max\Big\{k\big(r_1^{(i)},r_2^{(i)},d^{(i)}\big)\big/r_1^{(i)}:\ i=1,\ldots,n\Big\}.$$
Define $O_i$ as the projection of 
$\matC\times\{\lambda\cdot k(r_1^{(i)},r_2^{(i)},d^{(i)})\}$
from $\nt_i$ to $N$. Then $\{O_i\}_{i=1}^n$ is a horospherical cross-section
as in Proposition~\ref{new:conditions:for:conv:prop}.
\end{prop}

\dimostraz
If $O_i$ is first defined as the projection of 
$\matC\times\{k(r_1^{(i)},r_2^{(i)},d^{(i)})\}$, then $d(O_i,\partial N)$
equals $\log(k(r_1^{(i)},r_2^{(i)},d^{(i)})/r_1^{(i)})$,
and the conclusion follows from the same argument given for 
Corollary~\ref{esisteconv:cor}.
\finedimo

\section{The tilt formula}\label{tilt:section}
Proposition~\ref{riconoscimento:prop} shows that, to determine
whether a geometric decomposition of a hyperbolic manifold is Kojima's
canonical one, we must lift the decomposition to Minkowski 4-space and then
check convexity of all the angles at the 2-faces of the lifting.
In this section we provide the explicit formula that allows to check convexity.
This formula, already
obtained by Ushijima~\cite{ushijima} in more implicit terms,
extends Weeks' tilt formula~\cite{weeks:tilt}. Our main contribution
here is the computations of tilts in terms of moduli.
The statements of Proposition~\ref{tiltconv:prop}, 
Remark~\ref{horo:encode:rem}, Theorem~\ref{tiltformula:teo}, 
and Propositions~\ref{D:for:ideal:prop}
and~\ref{formula:prop}, with notation as in 
equation~(\ref{gtheta:def:formula}), 
may already be sufficient to proceed to the next section.

\paragraph{Tilts and convexity}
Let $\widehat{\Delta}$ be a tetrahedron in $\Minkos$ that projects
in $\pr$ to a tetrahedron 
$\Delta$ with vertices outside $\Hproj^3$ or on its boundary,
and edges meeting $\partial\Hproj^3$ or tangent to it.
Let $\widehat{F}$ be a face of $\widehat{\Delta}$ with image $F$ in $\Delta$.
Let $H$ be the unique half-space in $\Hproj^3$ 
such that $H\supset \Delta\cap\Hproj^3$ and $\partial H\supset F\cap\Hproj^3$.
Let $m\in\hp$ be the dual to $H$, and let $p\in\Minkos$ be the unique
vector such that ${\langle p,x\rangle=-1}$ for every $x\in\widehat{\Delta}$. 
We define the \emph{tilt} of $\widehat{\Delta}$ relative to 
$\widehat{F}$ as the real number $\langle m,p\rangle$.
The next result shows how tilts relate to convexity. For a proof
see~\cite{ushijima} or~\cite{weeks:tilt}.

\begin{prop}\label{tiltconv:prop}
Let $\widehat{\Delta}$ and $\widehat{\Delta}'$ 
be tetrahedra in $\Minkos$ sharing a $2$-face $\widehat{F}$. 
Assume that ${\widehat{\Delta}\cup\widehat{\Delta}'}$ 
projects injectively to $\pr$, and let $t$ and $t'$
be the tilts of $\widehat{\Delta}$ and $\widehat{\Delta}'$ relative to $\widehat{F}$.
Then the angle formed by $\widehat{\Delta}$ and $\widehat{\Delta}'$ at $\widehat{F}$ is
convex (respectively, 
strictly convex) if and only if ${t+t'\leqslant 0}$ (respectively,
${t+t'<0}$). 
\end{prop}

\paragraph{From moduli to tilts} Our task is now to
compute the tilts of the lifting
of a partially truncated tetrahedron from the intrinsic 
geometry of the tetrahedron itself. Recall that the lifting
of a non-ideal vertex $u$ is uniquely determined 
by the requirement that it should belong to $\hp$.
However, when $u$ is ideal, to get uniqueness we must choose
a horosphere at $u$ and lift $u$ to the dual in $\lp$ to this horosphere.
We begin by fixing some notation and 
recalling how horospheres are encoded~\cite{weeks:tilt}.

Let $\Delta$ be an abstract partially truncated tetrahedron.
Fix $\theta:\Delta\!^{(1)}\to[0,\pi)$ as in Theorem~\ref{moduli:teo}
and denote by $\Delta\!^{\theta,*}$ 
the corresponding geometric realization (up to isometry) of $\Delta\!^*$ in $\Hproj^3$.
Let $\Delta\!^\theta$ be the associated 
tetrahedron with ideal and ultra-ideal vertices in $\pr$.

\begin{rem}\label{horo:encode:rem}
\emph{If $u$ is an ideal vertex of $\Delta\!^\theta$, the set of horospheres at $u$
is parameterized by the positive reals, with a horosphere $O$ corresponding
to $r>0$ if $r$ is the radius of the smallest Euclidean disc on $O$ containing
$O\cap\Delta\!^{\theta}$.}
\end{rem}

So we fix a function $r:\calI\to(0,\infty)$,
where $\calI$ is the set of ideal vertices of $\Delta$,
and we denote by $\calO^{\theta,r}$
the associated family of horospheres at the ideal vertices of $\Delta\!^\theta$.
Now $\calO^{\theta,r}$ determines a unique lifting 
$\widehat{\Delta}^{\theta,r}\subset\Minkos$ 
of $\Delta\!^\theta$, and
$(\Delta\!^\theta,\widehat{\Delta}^{\theta,r})$
is well-defined up to isometry of pairs. 
Denoting by $u_1,\ldots,u_4$ the vertices of $\Delta$  and
by $F_i$ the face opposite
to $u_i$, we can then define $t_i^{\theta,r}$ as the tilt of
$\widehat{\Delta}^{\theta,r}$ relative to $\widehat{F}_i^{\theta,r}$. 
The tilt is unchanged under isometry, so
$t_i^{\theta,r}$ is indeed well-defined. 

To compute $t_i^{\theta,r}$ explicitly in terms of $\theta$ and $r$
we must introduce certain positive real numbers $D_i^{\theta,r}$.  
To this end, we denote by $H_i^\theta$ the plane in $\Hproj^3$
containing $F_i^{\theta,*}$. When $u_i$
is non-ideal we denote by $T_i^\theta$ the truncation plane for $\Delta\!^{\theta,*}$ 
corresponding to $u_i$, and when $u_i$ is ideal we denote by $O_i^{\theta,r}$ the 
horosphere at $u_i$ determined by $\theta$ and $r$. 
Then we define $D_i^{\theta,r}$ as $D(H_i^\theta,T_i^\theta)$ or 
$D(H_i^\theta,O_i^{\theta,r})$ depending on the type of $u_i$, 
where $D$ is the function defined as follows:
\begin{itemize} 
\item If $H$ and $T$ are planes in $\matH^3$ and $H\cap T=\emptyset$, we set
$D(H,T)=\cosh d(H,T)$ where $d$ is the usual distance in $\matH^3$;
if $H\cap T\neq \emptyset$ we set $D(H,T)=\cos \angle(H,T)$ where
$\angle(H,T)\in(0,\pi/2]$ is the angle formed by $H$ and $T$;
\item If $H$ is a plane in $\matH^3$ and $O$ is a horosphere not centred at
a point of $H$, we set $D(H,O)=\exp \pm \ell$, where 
$\ell$ is the length of the unique geodesic arc that joins $H$ to $O$ and is
orthogonal to both, with negative sign taken when $H\cap O\ne\emptyset$.
\end{itemize}

\begin{rem}
\emph{The choice of $\angle(H,T)$ in $(0,\pi/2]$ may look artificial 
at first sight, and indeed one could extend the definition to half-spaces rather
than planes, and choose $\angle(H,T)$ in $(0,\pi)$. However it is easy to show
that, in our situation, the dihedral angle at $H$ and $T$ that contains
$\Delta\!^{\theta,*}$ is always the acute one, so the definition of $D$ would remain
the same.}
\end{rem}

The following result was proved in~\cite{ushijima}:

\begin{teo}\label{tiltformula:teo}
Let $\theta$ and $r$ determine the geometry and a family of horospheres at the ideal
vertices of a partially truncated tetrahedron $\Delta$ with vertices
$u_1,\ldots,u_4$. Let $\widehat{\Delta}^{\theta,r}$ be the 
tetrahedron in $\Minkos$ determined by $\theta$ and $r$.
Let $t_i^{\theta,r}$ be the tilt of 
$\widehat{\Delta}^{\theta,r}$ relative to the face opposite to $u_i$.
Set $\theta_{ij}=\theta([u_i,u_j])$ and 
let $D_i^{\theta,r}$ be the number just introduced. Then:
$$\left(\begin{array}{c}
t_1^{\theta,r}\\ t_2^{\theta,r}\\ t_3^{\theta,r}\\ t_4^{\theta,r} 
\end{array} \right) = \left( \begin{array}{cccc}
1 & -\cos \theta_{34} & -\cos \theta_{24} & -\cos \theta_{23}\\
-\cos \theta_{34} & 1 & -\cos \theta_{14} & -\cos \theta_{13}\\
-\cos \theta_{24} & -\cos \theta_{14} & 1 & -\cos \theta_{12}\\
-\cos \theta_{23} & -\cos \theta_{13} & -\cos \theta_{12}  & 1
\end{array}\right) \left(\begin{array}{c}
1/D_1^{\theta,r}\\
1/D_2^{\theta,r}\\
1/D_3^{\theta,r}\\
1/D_4^{\theta,r} \end{array} \right).$$
\end{teo}

To make the calculation of tilts explicit, we are left to compute
the $D_i^{\theta,r}$'s. For ideal $u_i$ we denote $r(u_i)$ by $r_i$,
and the computation is easy:

\begin{prop}\label{D:for:ideal:prop}
If $u_1$ is ideal then:
$$D_1^{\theta,r}=\frac{1}{2r_1}\cdot
\frac{\sin \theta_{12}\cos \theta_{34}+ 
\sin \theta_{13}\cos \theta_{24}+\sin \theta_{14}
\cos\theta_{23}}{\sin \theta_{12}\sin 
\theta_{13}\sin\theta_{14}}.$$
\end{prop}
\dimostraz
We realize $\Delta\!^{\theta,*}$ in $\Hhalf^3$ setting $u_1=\infty$ 
and we denote by $C_i$ the trace at 
infinity of the plane that contains 
$F_i^{\theta,*}$. If $a$ is the circumradius of the Euclidean triangle
determined by $C_2,C_3,C_4$ and $a'$ is the Euclidean radius 
of $C_1$, then it is readily shown that ${D_1^{\theta,r}=a/(r_1\cdot a')}$.
Now the conclusion follows from the computation of $a/a'$ starting 
from the moduli, which involves only tools of elementary Euclidean 
\begin{figure}\begin{center}
\input{ratio_radii.pstex_t}
\caption{\small{The ratio of the radii of $C$ and $C_1$ is function
of the $\theta_{ij}$'s.}}\label{radii:ratio:fig}
\end{center}\end{figure}
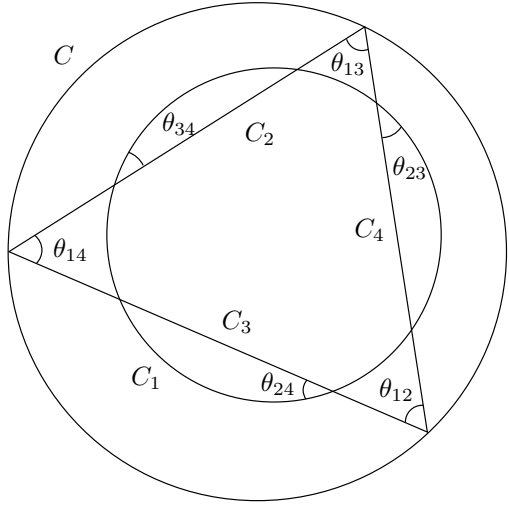
geometry, see Fig.~\ref{radii:ratio:fig}.
\finedimo

For non-ideal $u_i$ we denote
$D_i^{\theta,r}$ by $D_i^\theta$, because it is independent of $r$.
To compute it we need to introduce the following 
constant $g^\theta>0$:
\begin{equation}\label{gtheta:def:formula}
\begin{array}{rcl}
g^\theta =
-1\,+&\!\!\!\sum\!\!\!&\Big\{ \cos^2\theta(e):\ {e\in\Delta\!^{(1)}}\Big\}\\
+\, 2&\!\!\!\sum\!\!\!& 
\Big\{\cos\theta(e') \cos\theta(e'') \cos\theta(e'''):\ \
v\in\Delta\!^{(0)},\ v=e'\cap\ e''\cap e'''\Big\} \\
+\, 2&\!\!\!\sum\!\!\!&
\Big\{
\frac{\prod\{\cos\theta(e):\ {e\in\Delta^{(1)}}\}}
{\cos\theta(e')cos\theta(e'')}:\ \ 
\{e',e''\}\subset\Delta\!^{(1)},\ e'\cap e''=\emptyset\Big\}\\
-&\!\!\!\sum\!\!\!&
\Big\{\cos^2\theta(e') \cos^2\theta(e''):\ \ 
\{e',e''\}\subset\Delta\!^{(1)},\ e'\cap e''=\emptyset\Big\}.\\
\end{array}
\end{equation}
Of course $g^\theta$ is well-defined. Moreover:

\begin{prop}\label{formula:prop}
If $u_i$ is non-ideal and $d^\theta$ is as in formula~(\ref{dtheta}) of Section~\ref{mod:eqns:section} then
$$D_i^\theta=\sqrt{g^\theta/d^{\theta}(u_i)}.$$
\end{prop}

The proof of this result will be divided in several lemmas.
In the course of our argument we will need to use twice the
following explicit formula for the hyperbolic distance 
in $\Hhalf^2=\{z\in\mathbb{C}:\ \Im(z)>0\}$: for $x\in\matR$, $\rho>0$ and
$0<\alpha,\beta<\pi$ we have
\begin{equation}\label{h2half:distance:eqn}
d\big(x+\rho{\rm e}^{i\alpha},x+\rho{\rm e}^{i\beta}\big)=
\big|\log \tan(\alpha/2)-\log \tan(\beta/2)\big|.
\end{equation}

\begin{lemma}\label{con2:lem}
Let $c,a>0$ with $|c-a|<1$. In $\Hhalf^2$ let $\gamma_1$ and 
$\gamma_2$ be the
geodesics with ends at $\pm1$ and at $c\pm a$ respectively.
\begin{enumerate}
\item
If $\gamma_1$ and $\gamma_2$ are disjoint then
$\cosh d(\gamma_1,\gamma_2)=(1+a^2-c^2)\,/\,2a$.
\item
If $\gamma_1$ and $\gamma_2$ intersect at $z$ and
$\alpha=\angle(c-a,z,-1)$ is the 
angle they form then
$\cos \alpha=(1+a^2-c^2)\,/\,2a$. 
\end{enumerate}
\end{lemma}

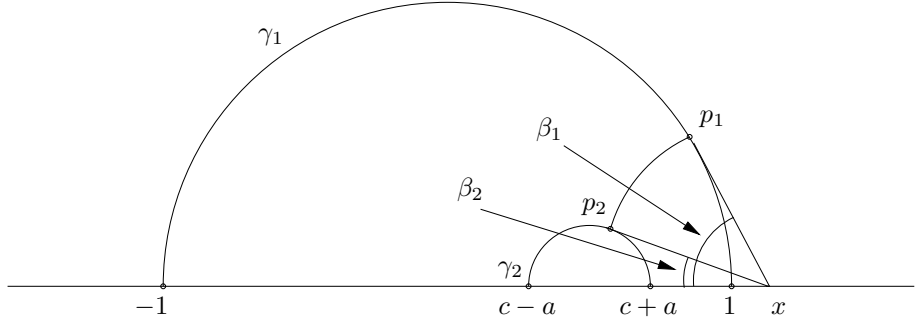
\begin{figure}
\begin{center}
\input{con1.pstex_t}
\caption{\small{The geodesic arc from $p_1$ to $p_2$ is the shortest
path between $\gamma_1$ and $\gamma_2$.}}
\label{con1:fig}
\end{center}
\end{figure}

\dimostraz
For point~(1),
let $p_1,p_2,\beta_1,\beta_2,x$ be as in Fig.~\ref{con1:fig}.
Since ${\sin\beta_1=1/x}$ and $\sin\beta_2=a/(x-c)$,
using~(\ref{h2half:distance:eqn})
we easily get $d(\gamma_1,\gamma_2)=d(p_1,p_2)=
\mathrm{arccosh\ } x-\mathrm{arccosh\ } ((x-c)/a).$
It follows that
\begin{equation}\label{espr1:formula}
\cosh d(\gamma_1,\gamma_2)=\Big(
x(x-c)-\sqrt{(x^2-1)((x-c)^2-a^2)}\Big)\,\Big/\, a.
\end{equation}
On the other hand, 
imposing that $|x-p_1|^2=|x-p_2|^2$ we get that
$x^2-1=(x-c)^2-a^2$.
Using this relation in the right-hand side
of~(\ref{espr1:formula}) we easily get the claimed equality.
A very similar argument proves point~(2).
\finedimo

\begin{figure}
\begin{center}
\input{tetra.pstex_t}
\caption{\small{Notations for the proof of Proposition~\ref{formula:prop}.}}
\label{tetra:fig}
\end{center}
\end{figure}
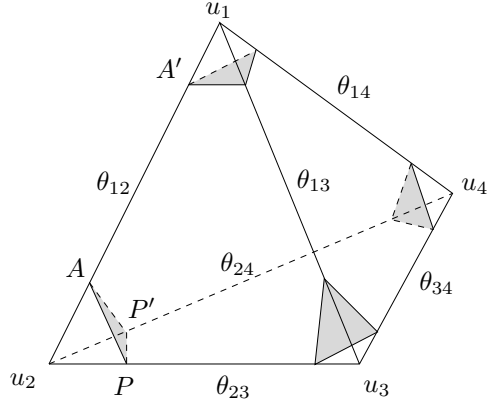

Now assume that $\Delta$ has neither ideal vertices nor 
length-0 edges, \emph{i.e.} that $\Delta\!^*$ is compact, and let
$A,A',P,P'$ be the points of $\Delta\!^*$ shown in Fig.~\ref{tetra:fig}.
Let $\Delta\!^{\theta,*}$ be a geometric realization of $\Delta$  
in ${\Hhalf^3=\matC\times(0,\infty)}$  
such that ${A= (0,\exp (-L^{\theta}
([u_1,u_2])))}$ and ${A'= (0,1)}$.
Let the hemisphere $H_1^\theta$ containing $F_1^{\theta,*}$ have 
Euclidean radius $R$ and centre $C=(z,0)$.
Note that the truncation planes $T_1^\theta$ and $T_2^\theta$ 
relative to $u_1$ and $u_2$ 
are hemispheres
centred at $(0,0)$ with Euclidean radii $1$ and 
${\exp (-L^{\theta} ([u_1,u_2]))}$ respectively.   

\begin{prop}\label{con3:prop}
$D_1^{\theta}=\frac{1}{R}\cdot\exp (-L^{\theta}([u_1,u_2]))
\cdot\sinh (L^{\theta}([u_1,u_2])).$
\end{prop}

\dimostraz
Since $H_1^\theta\perp T_2$, we have $|z|^2=R^2+\exp (-2L^{\theta}([u_1,u_2]))$.
Lemma~\ref{con2:lem} with $c=|z|$ and $a=R$ yields $D_1^{\theta}=\left(1-(\exp 
-2L^{\theta}([u_1,u_2]))\right)\,\big/\,2R$.
\finedimo

Now let $\ell=d(A,P)$ and $\ell'=d(A',P')$.
In the sequel we shall use the following equalities, which are readily 
deduced from~\cite[The Cosine Rule II, pag. 148]{beardon}: 
\begin{equation}\label{relk:formula}
\tanh \ell=\frac{\sqrt{d^{\theta}(u_2)}}
{\cos\theta_{12}\cdot\cos \theta_{23}+\cos \theta_{24}},
\qquad 
\tanh \ell'=\frac{\sqrt{d^{\theta}(u_2)}}
{\cos\theta_{12}\cdot\cos \theta_{24}+\cos \theta_{23}}.
\end{equation}

\begin{prop}\label{con4:prop}
$R^2=\exp (-2L^{\theta}([u_1,u_2]))\cdot
\sin^2 \theta_{12}\,\big/\,d^{\theta}(u_2).$
\end{prop}

\dimostraz
Using equality~(\ref{h2half:distance:eqn}) it 
is easily seen that 
\begin{eqnarray*}
P&=&\exp (-L^{\theta}([u_1,u_2]))\,\big/\,\cosh \ell\,\cdot\,
\big(\sinh \ell,1\big),
\\
P'&=&\exp (-L^{\theta}([u_1,u_2]))\,\big/\,\cosh \ell'\,\cdot\,
\big(\sinh \ell'\cos\theta_{12}
+i\sinh \ell'\sin\theta_{12},
1\big).
\end{eqnarray*}
Set $z=x+iy$.
Since $0P\perp CP$ and $0P'\perp CP'$, we have
\begin{equation}\label{centro}
\left\{
\begin{array}{ll} 
x\cdot\tanh \ell=\exp
(-L^{\theta}([u_1,u_2])),\\
x\cdot\tanh \ell'\cdot\cos\theta_{12}+ y\cdot\tanh \ell'
\cdot\sin \theta_{12}=\exp (-L^{\theta}([u_1,u_2])).\\
\end{array}
\right.
\end{equation}
The desired equality is now readily proved by
solving equations~(\ref{centro}) 
with respect to $x$ and $y$,
using~(\ref{relk:formula}),
and recalling that $R^2=x^2+y^2-(\exp(-2L^{\theta}([u_1,u_2])))$.
\finedimo

We can now prove
Proposition~\ref{formula:prop} for compact $\Delta\!^*$.
Equation~(\ref{internal:edge:length:formula}) yields
\[
\cosh L^{\theta}([u_1,u_2])=c^{\theta}([u_1, u_2])\,\Big/\,
\sqrt{d^{\theta}(u_1)d^{\theta}(u_2)},
\]
where $c^\theta([u_1, u_2])$ is defined by equation~(\ref{ctheta}).
By Propositions~\ref{con3:prop} and~\ref{con4:prop} we deduce
\begin{equation}\label{edue:formula} 
\left(D_1^{\theta}\right)^2=
\frac{\cosh^2 L^{\theta}([u_1,u_2])-1}{R^2\cdot\exp 2L^{\theta}
([u_1,u_2])}
=\frac{\left(c^{\theta}([u_1, u_2])\right)^2-d^{\theta}(u_1)d^{\theta}(u_2)}
{d^{\theta}(u_1)\sin^2\theta_{12}}.
\end{equation}
A long but straight-forward computation shows that
the right-hand side of   
equation~(\ref{edue:formula}) is in fact equal 
to ${g^{\theta}/d^{\theta}(u_1)}$.
This proves
Proposition~\ref{formula:prop} when $\Delta$ is a truncated tetrahedron
with no ideal vertices and no length-0 edges.
In the general case we 
can approximate a geometric realization  
$\Delta\!^{\theta,*}$ of any partially truncated
tetrahedron $\Delta$ with geometric realizations 
of compact truncated tetrahedra.
Using Proposition~\ref{formula:prop} in the compact case 
and a standard continuity argument 
we then deduce that the proposition
holds in general.

\section{Computing the canonical triangulation}\label{computing:section}
In this section we show how to compute the canonical decomposition $\kojiNstar$
of a hyperbolic 3-manifold $N$ starting from an arbitrary geometric 
triangulation $\triaNstar$ of $N$. This is achieved by a step-by-step 
modification of $\triaNstar$ until a triangulation is reached whose 
lifting to $\Minkos$ has only convex angles. 
According to Proposition~\ref{riconoscimento:prop},
if all the angles are actually strictly convex then $\kojiNstar=\triaNstar$, otherwise
$\kojiNstar$ is obtained from $\triaNstar$ by removing the 2-faces at which 
the lifting has flat angles. We warn the reader that, just as in~\cite{weeks:tilt},
the process of 
transforming $\triaNstar$ into $\kojiNstar$ may a priori get stuck at some
point, so we are not entitled to call it an algorithm in a strict sense.
On the other hand, in the next section we will show that if the process does
not get stuck then it converges in finite time.
The essential points of this section are the initial paragraph about
topological and geometric moves,  Theorem~\ref{self:convex:teo} and
the outline of the algorithm described in the last
two paragraphs.

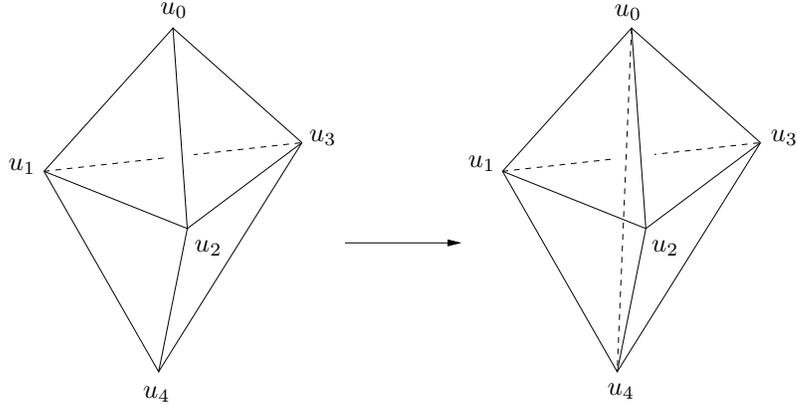
\begin{figure} 
\begin{center}
\input{mapimove.pstex_t}
\caption{\small{The two-to-three move.}}\label{dual:mapimove:fig}
\end{center}\end{figure}

\paragraph{Topological and geometric moves}
The fundamental move of the Matveev-Piergallini calculus for 
topological ideal triangulations
is the two-to-three move, already mentioned in Section~\ref{trunc:tria:section} 
and shown in Fig.~\ref{dual:mapimove:fig}. To fix notation, let us say
that the move replaces two distinct tetrahedra $T_{0123}$ and 
$T_{1234}$ sharing a face $F_{123}$ with three distinct tetrahedra
$T_{0124}$, $T_{0134}$, and $T_{0234}$ sharing an edge $e_{04}$.
This move cannot always be performed in a geometric setting,
but when it can we call the initial pair of tetrahedra an 
\emph{admissible} one. More precisely:

\begin{defn}\label{adm:tetra:def}
\emph{A triple $(T_{0123},T_{1234},F_{123})$ consisting of two
Euclidean tetrahedra and their common face, embedded in $\pr\cong\matE^3$ as in
Fig.~\ref{dual:mapimove:fig}-left, is called \emph{admissible} if
$e_{04}$ meets the interior of $F_{123}$.}
\end{defn}

Recall now that we are considering a partially truncated
triangulation $\triaNstar$ of a hyperbolic $N$.
We denote by $\triaN$ the corresponding abstract triangulation,
we fix the universal covering $\Hproj^3\supset\nt\to N$
and note that $\triaN$ determines a triangulation $\triaproj$ 
contained in $\pr$. In the sequel we will often lift tetrahedra
from $\triaN$ to $\triaproj$: the reader is invited to check that
all our considerations are independent of the lifting chosen.

\begin{defn}\label{adm:face:def}
\emph{If $F$ is a 2-face of $\triaN$ and the two tetrahedra
$\Delta$ and $\Delta'$ incident to $F$ are distinct, we call $F$
\emph{admissible} if the lifting to $\triaproj\subset\pr\cong\matE^3$
of the triple $(\Delta,\Delta',F)$
is a Euclidean admissible triple.}
\end{defn}

\begin{rem}\label{adm:then:move:rem}
\emph{If $F$ is an admissible $2$-face of $\triaN$ then
the two-to-three move that destroys $F$ yields a new
geometric partially truncated triangulation of $N$.}
\end{rem}

Turning to the inverse (three-to-two) move, we show that it is
always geometric:

\begin{lemma}\label{three:to:two:OK:lem}
Assume in $\triaN$ there are precisely three distinct tetrahedra sharing a 
non-$0$-length edge. Lift the edge and the tetrahedra to $\triaproj$, 
with notation as in Fig.~\ref{dual:mapimove:fig}-right.
Then the triple $(T_{0123},T_{1234},F_{123})$ is admissible, and the
three-to-two move in $\triaN$ that destroys $e_{04}$ gives rise to a new geometric
triangulation of $N$.
\end{lemma}

\dimostraz
Since $T_{0124},T_{0234},T_{1234}$ are cyclically arranged around $e_{04}$,
the line $r_{04}$ through $u_0$ and $u_4$ meets the interior of $F_{123}$.
To show that $r_{04}\cap F_{123}$ is actually a point of the interior
of $e_{04}$, we must show that $u_0$ and $u_4$ cannot 
lie on opposite sides of the plane which contains $F_{123}$. 
If this were the case, using again the fact that 
$r_{04}\cap F_{123}\ne\emptyset$, we would deduce that 
$u_0\in T_{1234}$ up to interchanging $u_0$ and $u_4$. 
From the Euclidean point of view, $u_0,\ldots,u_4$ lie outside
the unit ball $\Hproj^3$ or on its boundary, and all the edges 
$e_{ij}$ meet the ball
or are tangent to it. 
So $T_{1234}\setminus\Hproj^3$ is the union of (at most) four regions $W_i$, where
$W_i$ is star-shaped with respect to $u_i$. Since $u_0$ belongs to one of the
$W_i$'s, the corresponding edge
$e_{0i}$ does not meet $\Hproj^3$. A contradiction.
\finedimo

\paragraph{Effectiveness of moves}
To transform a geometric triangulation into Kojima's canonical decomposition, we will apply
both the two-to-three and the three-to-two moves, trying to 
remove concave angles from the lifting to $\Minkos$. The next result shows that,
when we remove a concave angle, the new ones that we create are not concave. Note 
however that some of the ``old'' convex angles may become concave.

\begin{lemma}\label{casi:lem}
Let $(T_{0123},T_{1234},F_{123})$ be an admissible triple in $\pr$, with
notation as in Fig.~\ref{dual:mapimove:fig}-left. For $i=0,\ldots,4$ let 
$\widehat{u}_i$ be a lifting of $u_i$ to $\Minkos$. Assume that 
$\widehat{u}_0,\ldots,\widehat{u}_4$ are affinely independent in $\Minkos$
and their convex hull $X$ does not contain $0$. Let $Y$
be the set of almost-visible points of $X$. Then
one and only one of the following possibilities occurs:
\begin{itemize}
\item $Y=\widehat{T}_{0124}\cup \widehat{T}_{0134}\cup \widehat{T}_{0234}$,
the angles at $\widehat{F}_{014}$, 
$\widehat{F}_{024}$,  and $\widehat{F}_{034}$, are strictly convex 
and the angle at $\widehat{F}_{123}$ is strictly concave.
\item $Y=\widehat{T}_{0123}\cup\widehat{T}_{1234}$,the angle at
$\widehat{F}_{123}$ is strictly convex and the angles at
$\widehat{F}_{014}$, $\widehat{F}_{024}$ and $\widehat{F}_{034}$
are strictly concave. 
\end{itemize}
\end{lemma}

\dimostraz
Let 
$Y'=\widehat{T}_{0124}\cup \widehat{T}_{0134}\cup \widehat{T}_{0234}$
and 
$Y''=\widehat{T}_{0123}\cup\widehat{T}_{1234}$.
Admissibility of ${(T_{0123},T_{1234},F_{123})}$ 
readily implies that either $Y=Y'$ or $Y=Y''$. More precisely,
let ${c=e_{04}\cap F_{123}}$, 
${c'=\pi^{-1}(c)\cap \widehat{e}_{04}}$,
and ${c''=\pi^{-1} (c)\cap \widehat{F}_{123}}$, and define
$\lambda>0$ so that ${c''=\lambda\cdot c'}$.
Affine independence of $\widehat{u}_0,\ldots,\widehat{u}_4$ shows that $\lambda\ne1$.
Now ${Y=Y'}$ when $\lambda>1$,
and ${Y=Y''}$ when $\lambda<1$.

On the other hand, choosing coordinates on $\Minkos$ such that
${\widehat{T}_{0123}\subset\{x_0=1\}}$, we see that
the angle at $\widehat{F}_{123}$ 
is convex precisely when
${x_0(\widehat{u}_4)>1}$, namely when $\lambda<1$.
Similarly, with coordinates 
such that ${\widehat{T}_{0124}\subset\{x_0=1\}}$,
the angle at $\widehat{F}_{024}$
is convex precisely when $x_0(\widehat{u}_3)>1$, namely when
$\lambda>1$. This concludes the proof.
\finedimo

\paragraph{Self-adjacent tetrahedra}
Now let $\calO$ be a horospherical cross-section for $N$,
and consider the corresponding
lifting $\triaMinkos(\calO)$ of $\triaproj$ to $\Minkos$.
In $\triaMinkos(\calO)$
it is always possible to apply a two-to-three
move to a 2-face. Similarly, we can always apply a three-to-two move to an edge with
three neighbouring tetrahedra. The same moves may however be impossible
in $\triaN$, when the involved tetrahedra are not distinct. The next
result shows that in this case we actually do not need to worry about convexity of angles.
In other words, when we are tempted to make a move (because of concavity), then
we are guaranteed that the move is \emph{topologically} possible.
Recall however that a topological two-to-three move may not be geometric.

\begin{teo}\label{self:convex:teo}
Let $\triaNstar$ be a geometric triangulation of a hyperbolic $N$.
Take a universal cover
$\Hproj^3\supset\nt\to N$, let $\calO$ be a horospherical
cross-section for $N$, and let $\triaMinkos(\calO)$ be the associated lifting
of $\triaN$ to $\Minkos$. Then:
\begin{enumerate}
\item If $\widehat{F}$ is a $2$-face of $\triaMinkos(\calO)$ 
and the two tetrahedra
incident to $\widehat{F}$ have the same image in $\triaN$, then the angle at 
$\widehat{F}$ is strictly convex;
\item If $\widehat{e}$ is an edge of $\triaMinkos(\calO)$ 
with three incident 
tetrahedra $\widehat{\Delta}_i$ and three incident $2$-faces $\widehat{F}_i$, for
$i=1,2,3$, and 
the images of the $\widehat{\Delta}_i$'s in $\triaN$ are not pairwise distinct,
then the angles at the 
$\widehat{F}_i$'s are strictly convex.
\end{enumerate}
\end{teo}

The proof of this result makes a crucial use
of Proposition~\ref{lengths:are:moduli}, announced
but not shown in Section~\ref{mod:eqns:section}.
So we show it now before proceeding.

\dimo{lengths:are:moduli}
We use $\Hproj^3\subset\pr\subset\Minkos$ 
and take a tetrahedron $\Delta$ with vertices 
$u_i$ outside $\Hproj^3\cup\partial\Hproj^3$ and
edges $e_{ij}$ that intersect $\Hproj^3$ or are tangent to it.
Let ${\ell_{ij}\in[0,\infty)}$ be the length of $e_{ij}\cap\Delta\!^*$.
We must show that if $\Delta'$ is another such tetrahedron
and $\ell'_{ij}=\ell_{ij}$ then there exists an isometry between
$\Delta\!^*$ and $\Delta\!^{\prime\,*}$.

Let $\widehat{u}_i$ be the only positive multiple of $u_i$ that
lies in $\hp$. Since $\ell_{ij}$ is the distance between the truncation
planes for $\Delta$ relative to $u_i$ and $u_j$, and these planes bound the 
half-spaces dual to $\widehat{u}_i$ and $\widehat{u}_j$,
Lemma~\ref{distances:lem}~(1)
shows 
that $\langle \widehat{u}_i,\widehat{u}_j\rangle=-\cosh \ell_{ij}$.
This implies that $\langle \widehat{u}'_i,\widehat{u}'_j\rangle=
\langle \widehat{u}_i,\widehat{u}_j\rangle$ for all $i,j$ (including $i=j$).
Now it is easy to see that $(\widehat{u}_i)_{i=1}^4$ and 
$(\widehat{u}'_i)_{i=1}^4$ are bases of $\Minkos$. Then there
exists an isometry $\varphi$ of $\Minkos$ such that 
$\varphi(\widehat{u}_i)=\widehat{u}'_i$ for all $i$, and the conclusion follows.
\finedimo

\dimo{self:convex:teo}
For point~(1), let $\widehat{\Delta}_1$ and $\widehat{\Delta}_2$ be the tetrahedra of
$\triaMinkos(\calO)$ incident to $\widehat{F}$ and let
$\Delta$ be their common image in $\triaproj$. 
Let $\theta:\Delta\!^{(1)}\to[0,\pi)$ give the dihedral angles of $\Delta$
and let ${r:\calI\to (0,\infty)}$ be the map determined by $\calO$
as explained in Remark~\ref{horo:encode:rem}, 
where $\calI$ is the set of ideal vertices of $\Delta$.
Let $u_1,\ldots,u_4$ be the vertices of $\Delta$ and let $F_i$ be the face
opposite to $u_i$, with notation such that 
${(\widehat{\Delta}_i,\widehat{F})}$ projects to $(\Delta,F_i)$ for $i=1,2$.
According to Proposition~\ref{tiltconv:prop}, to prove 
strict convexity 
at $\widehat{F}$ we have to check that ${t_1^{\theta,r}+t_2^{\theta,r}<0}$.
To do so we will need to discuss various possibilities for the
geometry of $\Delta$. Recall first that 
the combinatorial data defining $\calT$ determine a
simplicial isomorphism ${\varphi:F_1\to F_2}$ , and 
$\varphi$ induces
an isometry ${\varphi^*:F_1^{\theta,*}\to F_2^{\theta,*}}$.

We consider first the case where $\theta_{34}=0$. In this case 
Theorem~\ref{tiltformula:teo} gives
\begin{eqnarray*}
t^{\theta,r}_1&=&1/D_1^{\theta,r}-1/D_2^{\theta,r}-\cos \theta_{24}/D_3^{\theta,r}
-\cos \theta_{23}/D_4^{\theta,r},\\
t^{\theta,r}_2&=&-1/D_1^{\theta,r}+1/D_2^{\theta,r}-\cos \theta_{14}/D_3^{\theta,r}
-\cos \theta_{13}/D_4^{\theta,r}.
\end{eqnarray*}
Now $\cos \theta_{24}+\cos\theta_{14}>0$ and 
$\cos\theta_{23}+\cos\theta_{13}>0$, so
\[
t^{\theta,r}_1+t^{\theta,r}_2
=-\Big( (\cos \theta_{24}+\cos\theta_{14})/D_3^{\theta,r}+
(\cos\theta_{23}+\cos\theta_{13})/D_4^{\theta,r}\Big)<0.
\]

The case $\theta_{34}=0$ is settled, so we will assume
henceforth that $\theta_{34}\ne0$.
This implies that 
${\varphi(u_2)\ne u_1}$, otherwise the total dihedral angle in $N$ along
the image of $e_{34}$ would reduce to $\theta_{34}$, but $\theta_{34}<2\pi$.
Orientability of $N$ then implies 
that, up to interchanging $u_3$ and $u_4$, we have 
\begin{equation}\label{cond1:eqn}
\varphi(u_2)=u_4,\ \varphi(u_3)=u_1,\ \varphi(u_4)=u_3.
\end{equation}

These conditions easily imply that either all the $u_i$'s are ideal
or none of them is. Assume first they are all ideal.
Then the dihedral angles along opposite edges are
the same. We set $\alpha=\theta_{12}=\theta_{34}$, $\beta=\theta_{13}=\theta_{24}$,
$\gamma=\theta_{14}=\theta_{23}$, and note further that $\alpha+\beta+\gamma=\pi$.
Setting $r_i=r(u_i)$ and using
Proposition~\ref{D:for:ideal:prop} we see that
$r_i=1/D_i^{\theta,r}$. Recall now that the length of an edge 
of a Euclidean triangle is twice the circumradius
times the sine of the opposite angle. 
Since $\varphi^*$ matches the triangular cross-sections determined by the $r_i$'s
at the vertices of $\Delta\!^{\theta,*}$, we have:
\[
r_2\sin\alpha=r_4\sin\beta,\quad r_3\sin\beta=r_1\sin\alpha,\quad r_4\sin\gamma=
r_3\sin\gamma,
\]
whence $r_2=r_1$ and 
$r_4=r_3=r_1\cdot{\sin\alpha}\,/\,{\sin\beta}$.
Using Theorem~\ref{tiltformula:teo} we then get
$$t^{\theta,r}_1=t^{\theta,r}_2=r_1\cdot\big(1-\cos\alpha-(\cos\beta+\cos\gamma)\cdot
{\sin\alpha}\,/\,{\sin\beta}\big).$$
Relation ${\alpha+\beta+\gamma=\pi}$ now implies that
$t^{\theta,r}_i=r_1\cdot
\sin\gamma\cdot(\cos\alpha -1)\,/\,{\sin\beta}< 0$.

The only case left to settle to prove point (1) is when $\theta_{34}\ne 0$
and no $u_i$ is ideal. Using~(\ref{cond1:eqn}) and
the fact that $\varphi^*$ is an isometry we see that
\[
L^{\theta}(e_{23})=L^{\theta}(e_{14}),\qquad
L^{\theta}(e_{24})=L^{\theta}(e_{34})=L^{\theta}(e_{13}).
\]
Proposition~\ref{lengths:are:moduli} just proved and these relations now
imply that there exists an isometry of $\Delta\!^{\theta,*}$ 
onto itself that interchanges $u_1$ with $u_2$ and $u_3$ with $u_4$.
Then, with notation as in 
Theorem~\ref{tiltformula:teo}, we have
\[
D_1^{\theta}=D_2^{\theta},\ D_3^{\theta}=D_4^{\theta},
\ \theta_{23}=\theta_{14},\ \theta_{13}=\theta_{24}.
\]
Now we define 
$x=\cos \theta_{24}=\cos \theta_{13},\ y=\cos \theta_{23}=\cos \theta_{14},
\ z=\cos \theta_{34},\ w=\cos \theta_{12}$.
From Theorem~\ref{tiltformula:teo} 
and Proposition~\ref{formula:prop}
we deduce that
\begin{equation}\label{primotilt:formula}
t_1^\theta=t_2^\theta={\Big((1-z)\sqrt{d^{\theta}(u_1)}-(x+y)
\sqrt{d^{\theta}(u_3)}\Big)}\,\Big/\,\sqrt{g^{\theta}}.
\end{equation}
We claim that
the following formula holds:
\begin{equation}\label{w:formula}
w=(x+zy)\sqrt{(1-x^2)\,\big/\,(1-z^2)}-xy.
\end{equation}
We first show that~(\ref{w:formula}) implies
$t_i^\theta<0$. Later we will establish~(\ref{w:formula}).
Note that
$d^{\theta}(u_1)=x^2+y^2+w^2+2xyw-1$ and 
$d^{\theta}(u_3)=x^2+y^2+z^2+2xyz-1$. Using~(\ref{w:formula}) we then get
$d^{\theta}(u_1)=d^{\theta}(u_3) \cdot \frac{1-x^2}{1-z^2}.$
By equation~(\ref{primotilt:formula}) it follows that
\[
t^\theta_i=\sqrt{d^{\theta}(u_3)}\cdot
\Big((1-z)\sqrt{(1-x^2)\,\big/\,(1-z^2)}-(x+y)\Big)\,\Big/\, \sqrt{g^\theta}.
\]
Since $\theta_{14}+\theta_{24}+\theta_{34}<\pi$
we deduce that
\begin{equation}\label{condang:formula}
x+y> 1-z\geqslant 0.
\end{equation}
Then
$t^\theta_i<0 \Leftrightarrow  x+y>(1-z)\sqrt{\frac{1-x^2}{1-z^2}}
\Leftrightarrow (x+y)^2>\Big( (1-z)\sqrt{\frac{1-x^2}{1-z^2}}\Big)^2.$
After some computations we deduce that
\[ 
t^\theta_i<0
\Leftrightarrow (x^2+y^2+z^2+2xyz-1)+(x+y)^2+(1-z)(z-y)^2>0.
\]
Since $x^2+y^2+z^2+2xyz-1=d^\theta(u_3)>0$, 
it is now sufficient to show that
\[
(x+y)^2+(1-z)(z-y^2)\geqslant 0,
\]
and this inequality follows quite easily from~(\ref{condang:formula}). 

We are left to establish~(\ref{w:formula}). To do so we distinguish two
cases, according to whether $\theta_{14}$ vanishes or not. 
We first assume $\theta_{14}\ne 0$. 
Let $T_2$ and $T_4$ be the 
truncation triangles relative to $u_2$ and $u_4$ respectively
and let $\ell,\ell'$ be the boundary edges of $\Delta\!^{\theta,*}$ defined
by $\ell=F_1^{\theta,*}\cap T_2$, $\ell'=F_2^{\theta,*}\cap T_4$.   
Now $\theta_{34}\neq 0$ and 
$\theta_{14}\neq 0$, so $\ell'$ has
finite length. Moreover, we know that 
$\varphi^*$ identifies $\ell$ to $\ell'$.
Equation~(\ref{w:formula}) 
is now obtained by equaling the lengths of $\ell$ and $\ell'$ via
Proposition~\ref{boundary:edge:length:prop}. 
We turn at last to the case when 
$\theta_{14}=0$. 
Imposing $L^{\theta}(e_{24})=L^{\theta}(e_{34})$
via Proposition~\ref{internal:edge:length:prop} we get
after some computations 
$(x+w)^2(1-z^2)=(x+z)^2(1-x^2)$.
Since $y=1$, this relation is in fact equivalent to
equation~(\ref{w:formula}), and we are done.
 
Now we show point~(2). 
As a consequence of~(1), we have strict convexity along at least one of the
$\widehat{F}_i$'s. Lemma~\ref{casi:lem} then implies strict convexity at all the 
$\widehat{F}_i$'s.
\finedimo

\paragraph{Outline of the algorithm}
The input of our process is a geometric triangulation $\triaNstar$
of a hyperbolic 3-manifold $N$. More precisely, we start with
a topological partially truncated
triangulation $\triaN$ of $N$ in the sense of Proposition~\ref{tria:of:N':prop}, and 
a solution of the system
of consistency and completeness equations for $\triaN$, as 
discussed in Section~\ref{mod:eqns:section}.
Then we perform the following steps:
\begin{enumerate}
\item We choose a horospherical cross-section $\calO$
as in Proposition~\ref{new:conditions:for:conv:prop};
\item We pick a 2-face $F$ of $\triaN$ such that the two tetrahedra of $\triaN$
incident to $F$ are distinct. We lift $F$ and its two incident tetrahedra
to $\triaMinkos(\calO)$. Using Proposition~\ref{tiltconv:prop}, 
we check whether the angle at the lifted face is strictly concave. If it is, 
we move to step 3. If it is not, we move to another 2-face. 
If all faces are visited and no concave angle is found, 
$\triaNstar$ is the output 
(because it is $\kojiNstar$ or a subdivision of it,
by Proposition~\ref{riconoscimento:prop} and Theorem~\ref{self:convex:teo});
\item If $F$ is admissible, we change $\triaN$ by performing
the geometric two-to-three move that kills $F$, 
and we go back to step 2. If $F$ is non-admissible, we 
check whether one of the non-$0$-length edges of $F$ is shared by precisely three
tetrahedra of $\triaN$. If it is, we change $\triaN$ by 
applying the geometric three-to-two move that kills this edge, and we go back to step 2.
If it is not, we do not change $\triaN$ but we go back to step 2 moving
to a different concave face. If all concave faces are visited and no move can be applied
to any of them, we give up.
\end{enumerate}
In the next section we will show
that the process, if it does not get stuck during step 3, 
outputs the canonical decomposition in finite time. 
Steps 2 and 3 are of course directly implementable, whereas step 1 requires
a careful discussion, to which the rest of the present section is devoted.

\paragraph{Algorithmic choice of horospherical cross-sections}
According to Proposition~\ref{verysafe:height:prop}, 
to determine a horospherical cross-section as in 
Proposition~\ref{new:conditions:for:conv:prop}
we must find for each cusp a realization of $\nt$ in $\Hhalf^3$ so that
the cusp is generated by $\infty$, and compute the corresponding
$d,r_1,r_2$ of Proposition~\ref{safe:height:prop}. 
We recall that the datum to use is a geometric triangulation $\triaN$
of $N$.

Let us concentrate on a cusp $C$ and 
fix a tetrahedron $\Delta_0^*\in\triaNstar$ with a certain ideal vertex $v_0$ 
asymptotic to $C$. We take a realization $\tilde\Delta_0^*$ in $\Hhalf^3$
such that $v_0$ gets identified to $\infty$. Here and in the sequel the realizations
we consider are of course all 
compatible with the geometric structure given on the
tetrahedra. Choosing a horospherical cross-section
at $C$ now amounts to choosing a positive real number, namely the height at which 
the lifted cross-section should intersect $\tilde\Delta_0$. In the course
of our argument, starting from $\tilde\Delta_0^*$, we will be successively gluing
new tetrahedra to free faces of tetrahedra we already have, as dictated by the
combinatorics and the geometry of $\triaNstar$. We warn the reader that it is not possible
to predict \emph{a priori} how many different copies of each tetrahedron
of $\triaNstar$ will need to be glued, but the process is guaranteed to be finite anyway,
as we will carefully explain.

\smallskip\noindent\textbf{Step 1.A.} \emph{We take one copy of each tetrahedron $\Delta\!^*$ of 
$\triaNstar$ for each vertex $v$ of $\Delta\!^*$ asymptotic to $C$,
and, starting from $\tilde\Delta_0^*$,
we do gluings along free vertical faces until each $(\Delta\!^*,v)$ has been
realized once in $\Hhalf^3$ with $v=\infty$.}

After Step 1.A we have a certain finite family $\calF_1$ 
of partially truncated tetrahedra in $\Hhalf^3$, all having 
$\infty$ as a vertex, and we can compute the following:
\begin{itemize}
\item $\rho=\max\{\rho(\tilde\Delta\!^*):\ \tilde\Delta\!^*\in\calF_1\}$,
where $\rho(\tilde\Delta\!^*)$ is the Euclidean radius of the
half-sphere that contains the face of $\tilde\Delta\!^*$ opposite to $\infty$;
\item $r=\max\bigcup\{r(\tilde\Delta\!^*):\ \tilde\Delta\!^*\in\calF_1\}$,
where $r(\tilde\Delta\!^*)$ is the set of Euclidean radii of 
the half-spheres that contain the truncation triangles of $\tilde\Delta\!^*$.
We define $r$ to be $-\infty$ if all the tetrahedra of $\calF_1$ are ideal,
and we note for later purpose that the definition of $r(\tilde\Delta\!^*)$
makes sense also if $\tilde\Delta\!^*$ does not have $\infty$ as a vertex; 
\item The intersection $\Omega$ of the horizontal plane at height $z=\max\{\rho,r\}$
with the union of the tetrahedra in $\calF_1$;
\item The first number we need to determine, \emph{i.e.}~the diameter $d$ of $\Omega$ with
respect to the ordinary Euclidean metric on $\matC\times\{z\}$.
\end{itemize}

\smallskip\noindent\textbf{Step 1.B.} \emph{Starting from $\calF_1$, 
we perform gluings along free non-vertical
faces, adding new truncated tetrahedra, until we get a family $\calF_2$ such that 
$\bigcup\{r(\tilde\Delta\!^*):\ \tilde\Delta\!^*\in\calF_2\}$ contains at least two distinct 
values $r'_1>r'_2$.}

The way to realize Step 1.B algorithmically is as follows. 
We list the free non-vertical faces of $\calF_1$, we perform the gluings
along these faces getting a family $\calF_1'$, and we check whether $\calF_1'$ already
works. If it does not, we proceed similarly with $\calF_1'$ instead of $\calF_1$, until
the desired $\calF_2$ is reached. Of course this procedure only has to be iterated a finite
number of times, even if the number of iterations is hard to predict \emph{a priori}.

\smallskip\noindent\textbf{Step 1.C.} \emph{Starting from $\calF_2$, we perform gluings 
along free non-vertical
faces, adding new truncated tetrahedra, until we get a family $\calF_3$ such that 
any further tetrahedron glued to 
$\calF_3$ along a non-vertical face would lie entirely outside
$\Omega\times[r'_2,\infty)$.}

Of course this step is also a finite one, even if its length
is not easily predictable. 
Note also that $\Omega$ has finite diameter
$d$, so we could replace $\Omega$ by an easier set, like a disc or a square.
The choice of $\calF_3$ guarantees that its union contains
$\nt\cap(\Omega\times[r'_2,\infty))$, so
the two other constants $r_1$ and $r_2$ we need to determine are now
the first and second largest elements  of
$\bigcup\{r(\tilde\Delta\!^*):\ \tilde\Delta\!^*\in\calF_3\}$.

\section{Finiteness}\label{finiteness:section}
This section is entirely devoted to the proof
that the algorithm to transform a geometric triangulation into Kojima's canonical 
decomposition, if it does not get stuck, converges in finite time.
This fact was already announced above and is
accurately stated as follows:

\begin{teo}\label{finiteness:teo}
Let $N$ be hyperbolic with non-empty boundary,
let $\calO$ be a horospherical cross-section for $N$ as in 
Proposition~\ref{new:conditions:for:conv:prop}, 
and let $\triaN$ be a geometric triangulation of $N$.
Then there exists an integer $\nu=\nu(N,\calO,\calT)$
such that the following holds:
Assume $\{\triaN^i\}_{i=0}^j$ is a sequence of
geometric triangulations of $N$ starting at $\triaN^0=\triaN$, and 
for all $i=0,\ldots,j-1$ we have:
\begin{itemize}
\item there is a $2$-face $\widehat{F}^i$
of $\triaMinkos^i(\calO)$ along which 
$\triaMinkos^i(\calO)$ has a strictly concave angle;
\item 
$\triaN^{i+1}$ is obtained from $\triaN^i$ by a two-to-three or a three-to-two
move killing $\widehat{F}^i$.
\end{itemize} 
Then $j\leqslant \nu$.
\end{teo}

For the proof we fix as above the universal cover $\matH^3\supset\nt\to N$ 
with deck transformation group $\Gamma<{\mathrm{Isom}(\matH^3)}$, and
we denote by $\widetilde{\calO}$ the lifting of $\calO$.
We start with a series of lemmas, the first of which is taken from~\cite{kojima}. 

\begin{lemma}\label{cuspbordo:lem}
A point $q\in\nt_\infty$ generates an annular cusp of $N$
if and only if it belongs to the circle at infinity
of two distinct components of
$\bnt$.
\end{lemma}
 
\begin{lemma}\label{penultimo:lemma}
Let $S_0$ be a component of $\bnt$ 
with stabilizer $\Gamma_0$ in $\Gamma$. If $c>0$ and
\[
A(S_0,c)=\big\{S\subset\bnt:
\ S \textrm{\ is\ a\ component\ of}\ \bnt,\ d(S,S_0)\leqslant c\big\},
\]
then $\Gamma_0$ leaves $A(S_0,c)$ invariant and 
$\# (A(S_0,c)/\Gamma_0)<\infty$.
\end{lemma}

\dimostraz
The first assertion is obvious. For $S\in A(S_0,c)$ we have either $d(S,S_0)>0$
or $d(S,S_0)=0$. Correspondingly we have a splitting
$A(S_0,c)=A^+(S_0,c)\sqcup A^0(S_0)$, which of course is $\Gamma_0$-equivariant.
Recall now from~\cite{kojima} that
$S_0/\Gamma_0=S_0/\Gamma$
is a complete finite-area hyperbolic surface.
Lemma~\ref{cuspbordo:lem} readily implies that there is a bijection between
$A^0(S_0)/\Gamma_0$ and the set of cusps of $S_0/\Gamma_0$, so
$A^0(S_0)/\Gamma_0$ is finite. 

We are left to show that 
$A^+(S_0,c)/\Gamma_0$ is finite. To this end note first that $\nt/\Gamma_0$
is complete hyperbolic with geodesic boundary (but probably infinite volume),
and its boundary components constitute a locally finite family.
Now let $q$ be a point of $\partial S_0$ 
that generates a cusp
of $S_0/\Gamma_0$, \emph{i.e.}~a point in $\partial S_0\subset \nt_{\infty}$ 
that generates an annular
cusp of $N$. Using Lemma~\ref{cuspbordo:lem} and realizing $\nt$ in $\Hhalf^3$
with $q=\infty$, it is easily proved that there exists a horoball $B_q$ centred
at $q$ such that, if $S$ is a component of $\bnt$ and $d(S,S_0)>0$,
then ${d(S,S_0)=d(S,S_0\setminus O_q)}$. Repeating this argument for all the
finitely many cusps of $S_0/\Gamma_0$ we deduce that for some $\varepsilon>0$
our quotient $A^+(S_0,c)/\Gamma_0$ naturally corresponds to the set
of boundary components of $\nt/\Gamma_0$ whose distance from the
$\varepsilon$-thick part of $S_0/\Gamma_0$ is positive and bounded by $c$.
Compactness of the $\varepsilon$-thick 
part of $S_0/\Gamma_0$ and local finiteness of the components of
$\partial(\nt/\Gamma_0)$ then imply the conclusion.
\finedimo

\begin{lemma}\label{ancorauno:lemma}
Let $O_0$ be a component of $\widetilde{\calO}$ 
with stabilizer $\Gamma_0$ in $\Gamma$. If $c>0$ and
\begin{eqnarray*}
A^{(1)}(O_0,c)&=&\big\{O\subset \widetilde{\calO}:\ O\ \textrm{is\ a\ component\ of}\ 
\widetilde{\calO},\ d(O,O_0)\leqslant c\big\},\\
A^{(2)}(O_0,c)&=&
\big\{S\subset\bnt:\ S\textrm{\ is\ a\ component\ of}\ \bnt,
\  d(S,O_0)\leqslant c\big\},
\end{eqnarray*}
then $\Gamma_0$ leaves $A^{(j)}(O_0,c)$ 
invariant and $\# (A^{(j)}(S_0,c)/\Gamma_0)<\infty$ for both $j=1,2$.
\end{lemma}

\dimostraz
Realize $\nt$ in $\Hhalf^3=\matC\times(0,\infty)$ so that $O_0$ is centred at $\infty$.
Now the components of $\widetilde{O}$ are Euclidean spheres, and 
$A^{(1)}(O_0,c)$ consists of those whose radius is bounded from below by
a certain constant. Compactness of $\matC/\Gamma_0$ easily implies 
finiteness of $A^{(j)}(S_0,c)/\Gamma_0$ for $j=1$. A very similar
argument is employed for $j=2$.
\finedimo

Now we set $\calX_{\triaN}(\calO)=\{\alpha\cdot x:\ x\in\triaMinkos(\calO),
\ \alpha\geqslant 1\}$ and for real $c$ we define
$\mathcal{L}_{c}=\{v\in\Minkos:\ x_0(v)>0,\ \langle v,v\rangle< c\}.$

\begin{lemma}\label{ultimo:lemma}
There exists $c<0$ such that $\calX_{\triaN}(\calO)\supset\mathcal{L}_{c}$.
\end{lemma}

\dimostraz
Since $\Gamma$ acts isometrically on $\triaMinkos(\calO)$, the Lorentzian 
norm induces a continuous map
$\triaMinkos(\calO)\to\mathbb{R}$.
The domain of this map is homeomorphic to a compactification of $N$, so the 
map has a minimum $c$, which of course is negative. 
Knowing that the projection of $\triaMinkos(\calO)$ to $\pr$ contains
the unit ball $\Hproj^3$ we easily deduce the conclusion.
\finedimo

\begin{prop}\label{finitezza:prop}
There exist
only a finite number of $\Gamma$-inequivalent segments with
ends in ${\calB\cup \widetilde{Q}(\calO)}$ whose midpoint does not belong to
$\calX_{\triaN}(\calO)$.
\end{prop}

\dimostraz
Let $p_1,p_2\in\calB\cup \widetilde{Q}(\calO)$ be distinct and 
assume that the midpoint of $[p_1,p_2]$ 
does not belong to $\calX_{\triaN}(\calO)$. 
Choosing $c$ as in Lemma~\ref{ultimo:lemma}
we deduce that the Lorentzian norm of $(p_1+p_2)/2$ is at least $c$.
We will now consider three cases depending 
on whether $p_1$ and $p_2$ belong to 
$\calB$ or to $\widetilde{Q}(\calO)$. In all three cases we will show
that there are finitely many choices for $\{p_1,p_2\}$ up to the action
of $\Gamma$. 

\noindent\textsc{Case 1}: $p_1,p_2\in\calB$. 
Then
$c\leqslant \big\langle {(p_1+p_2)}/{2},{(p_1+p_2)}/{2}\big\rangle
=\big(1+\big\langle p_1,p_2\big\rangle\big)\big/{2}$.
If $S_i$ is the component
of $\bnt$ dual to $p_i$,
using Lemma~\ref{distances:lem}~(1) we deduce that
$\cosh d(S_1,S_2)\leqslant 1-2c$.
Finiteness of the number of components of $\bn$ 
and Lemma~\ref{penultimo:lemma} then imply the desired finiteness.

\noindent\textsc{Case 2}: $p_1,p_2\in \widetilde{Q}(\calO)$.
Then 
$c\leqslant \big\langle {(p_1+p_2)}/{2},{(p_1+p_2)}/{2}\big\rangle
=\big\langle p_1,p_2\big\rangle\big/2$.
If $O_i$ is  the horosphere dual
to $p_i$, we deduce that  
$\exp d(O_1,O_2)\leqslant -c$
by Lemma~\ref{distances:lem}~(4).
We then get the desired conclusion using
finiteness of the number of toric cusps of $N$
and Lemma~\ref{ancorauno:lemma} with $j=1$.

\noindent\textsc{Case 3}: $p_1\in \widetilde{Q}(\calO),p_2\in\calB$.
Then
$c\leqslant \big\langle {(p_1+p_2)}/{2},{(p_1+p_2)}/{2}\big\rangle
=\big(1+2\big\langle p_1,p_2\big\rangle\big)\big/{4}$.
If $O$ is the horosphere dual
to $p_1$ and $S$ is the component
of $\bnt$ dual to $p_2$, we deduce that 
$2\cdot\exp d(O,S)\leqslant 1-4c$
by Lemma~\ref{distances:lem}~(3).
We conclude using
again finiteness of the number of toric cusps of $N$
and Lemma~\ref{ancorauno:lemma} with $j=2$.
\finedimo

\dimo{finiteness:teo}
Lemma~\ref{casi:lem} implies that 
${\calX_{\triaN^i}(\calO)\subset\calX_{\triaN^{i+1}}(\calO)}$ for all $i$. Moreover
${\calX_{\triaN^{i+1}}(\calO)}$ contains at least an edge
with endpoints in $\calB\cup\widetilde{Q}(\calO)$ whose midpoint does not belong to
${\calX_{\triaN^{i}}(\calO)}$. Then we achieve the desired property
by defining $\nu$
as the number of $\Gamma$-inequivalent edges with endpoints in 
$\calB\cup\widetilde{Q}(\calO)$
whose midpoints do not belong to
$\calX_{\triaN}(\calO)$.
\finedimo

\vspace{1cm}

\noindent\hspace{6.5cm} Scuola Normale Superiore

\noindent\hspace{6.5cm} Piazza dei Cavalieri, 7

\noindent\hspace{6.5cm} 56126 Pisa, Italy

\noindent\hspace{6.5cm} frigerio@sns.it

\vspace{.2cm}

\noindent\hspace{6.5cm} Dipartimento di Matematica Applicata

\noindent\hspace{6.5cm} Universit\`a di Pisa

\noindent\hspace{6.5cm} Via Bonanno Pisano, 25/B

\noindent\hspace{6.5cm} 56126 Pisa, Italy

\noindent\hspace{6.5cm} petronio@dma.unipi.it

\end{document}

%% file: defs_hypeq_bd.tex
\usepackage{amssymb}
\usepackage{amsfonts}
\usepackage{amsmath}
\usepackage{amsthm}
\usepackage[english]{babel}
\usepackage{amscd,amsthm}
\usepackage[dvips]{graphics}

\newtheorem{lemma}{Lemma}[section] 
\newtheorem{teo}[lemma]{Theorem}
\newtheorem{rem}[lemma]{Remark} 
\newtheorem{prop}[lemma]{Proposition}
\newtheorem{cor}[lemma]{Corollary}

\newtheorem{defn}[lemma]{Definition}

\newcommand{\matR}{\ensuremath {\mathbb{R}}}

\newcommand{\matZ}{\ensuremath {\mathbb{Z}}}
\newcommand{\matC}{\ensuremath {\mathbb{C}}}
\newcommand{\matP}{\ensuremath {\mathbb{P}}}
\newcommand{\matH}{\ensuremath {\mathbb{H}}}
\newcommand{\matM}{\ensuremath {\mathbb{M}}}
\newcommand{\matE}{\ensuremath {\mathbb{E}}}
\newcommand{\matS}{\ensuremath {\mathbb{S}}}

\newcommand{\calZ}{\ensuremath {\mathcal{Z}}}
\newcommand{\calX}{\ensuremath {\mathcal{X}}}

\newcommand{\calK}{\ensuremath {\mathcal{K}}}
\newcommand{\calI}{\ensuremath {\mathcal{I}}}

\newcommand{\calS}{\ensuremath {\mathcal{S}}}

\newcommand{\calC}{\ensuremath {\mathcal{C}}}
\newcommand{\calA}{\ensuremath {\mathcal{A}}}
\newcommand{\calB}{\ensuremath {\mathcal{B}}}
\newcommand{\calD}{\ensuremath {\mathcal{D}}}
\newcommand{\calF}{\ensuremath {\mathcal{F}}}

\newcommand{\calT}{\ensuremath {\mathcal{T}}}

\newcommand{\calO}{\ensuremath {\mathcal{O}}}
\newcommand{\calV}{\ensuremath {\mathcal{V}}}

\reversemarginpar\textwidth=14cm \textheight=19cm

\font\titsc=cmcsc10 scaled 1200

\newcommand{\Nbar}{{\overline N}}
\newcommand{\Hhyp}{\matH_{\,{\rm hyp}}}
\newcommand{\Hproj}{\matH_{\,{\rm proj}}}

\newcommand{\Hhalf}{\matH_{\,{\rm half}}}

\newcommand{\nt}{\widetilde{N}}
\newcommand{\dimo}[1]{\vspace{2pt}\noindent\textit{Proof of \ref{#1}}.\ }
\newcommand{\finedimo}{{\hfill\hbox{$\square$}\vspace{2pt}}}
\newcommand{\dimostraz}{\vspace{2pt}\noindent\textit{Proof}.\ }

\newcommand{\bnt}{\partial\widetilde{N}}
\newcommand{\bn}{\partial{N}}          

\newcommand{\pr}{\Pi^3}

\newcommand{\hp}{\mathcal{H}_{+}^3}
\newcommand{\hm}{\mathcal{H}_{-}^3}
\newcommand{\Minkos}{\mathbb{M}^{3+1}}

\newcommand{\diag}{\SelectTips{cm}{} \xymatrix@1}

\newcommand{\lp}{L^3_+}

\newcommand{\triaN}{\calT}
\newcommand{\triaNstar}{\calT^*}

\newcommand{\triaproj}{\calT_\matP}
\newcommand{\triaMinkos}{\calT_\matM}

\newcommand{\decoN}{\calD}
\newcommand{\decoNstar}{\calD^*}
\newcommand{\decontstar}{\calD_{\nt}^*}
\newcommand{\decoproj}{\calD_\matP}
\newcommand{\decoMinkos}{\calD_\matM}

\newcommand{\kojiNstar}{\calK_N^*}
\newcommand{\kojintstar}{\calK_{\nt}^*}
\newcommand{\kojiproj}{\calK_\matP}
\newcommand{\kojiMinkos}{\calK_\matM}

%% file: geomreal.pstex_t
\begin{picture}(0,0)%
\includegraphics{geomreal.pstex}%
\end{picture}%
\setlength{\unitlength}{3355sp}%
\begingroup\makeatletter\ifx\SetFigFont\undefined%
\gdef\SetFigFont#1#2#3#4#5{%
  \reset@font\fontsize{#1}{#2pt}%
  \fontfamily{#3}\fontseries{#4}\fontshape{#5}%
  \selectfont}%
\fi\endgroup%
\begin{picture}(6195,3308)(871,-6376)
\put(1921,-3916){\makebox(0,0)[lb]{\smash{\SetFigFont{10}{12.0}{\rmdefault}{\mddefault}{\updefault}\special{ps: gsave 0 0 0 setrgbcolor}$e\in\calZ$\special{ps: grestore}}}}
\put(871,-5551){\makebox(0,0)[lb]{\smash{\SetFigFont{10}{12.0}{\rmdefault}{\mddefault}{\updefault}\special{ps: gsave 0 0 0 setrgbcolor}$v\in\calI$\special{ps: grestore}}}}
\put(4448,-6278){\makebox(0,0)[lb]{\smash{\SetFigFont{10}{12.0}{\rmdefault}{\mddefault}{\updefault}\special{ps: gsave 0 0 0 setrgbcolor}$v$\special{ps: grestore}}}}
\put(5506,-3233){\makebox(0,0)[lb]{\smash{\SetFigFont{10}{12.0}{\rmdefault}{\mddefault}{\updefault}\special{ps: gsave 0 0 0 setrgbcolor}$e$\special{ps: grestore}}}}
\end{picture}

%% file: newtronco.pstex_t
\begin{picture}(0,0)%
\includegraphics{newtronco.pstex}%
\end{picture}%
\setlength{\unitlength}{1895sp}%
\begingroup\makeatletter\ifx\SetFigFont\undefined%
\gdef\SetFigFont#1#2#3#4#5{%
  \reset@font\fontsize{#1}{#2pt}%
  \fontfamily{#3}\fontseries{#4}\fontshape{#5}%
  \selectfont}%
\fi\endgroup%
\begin{picture}(12024,6134)(-3011,-7163)
\put(4126,-2611){\makebox(0,0)[lb]{\smash{\SetFigFont{10}{12.0}{\rmdefault}{\mddefault}{\updefault}\special{ps: gsave 0 0 0 setrgbcolor}$\theta_2$\special{ps: grestore}}}}
\put(5401,-1861){\makebox(0,0)[lb]{\smash{\SetFigFont{10}{12.0}{\rmdefault}{\mddefault}{\updefault}\special{ps: gsave 0 0 0 setrgbcolor}$\theta_3$\special{ps: grestore}}}}
\put(8176,-3661){\makebox(0,0)[lb]{\smash{\SetFigFont{10}{12.0}{\rmdefault}{\mddefault}{\updefault}\special{ps: gsave 0 0 0 setrgbcolor}$\theta_4$\special{ps: grestore}}}}
\put(4426,-6511){\makebox(0,0)[lb]{\smash{\SetFigFont{10}{12.0}{\rmdefault}{\mddefault}{\updefault}\special{ps: gsave 0 0 0 setrgbcolor}$\theta_5$\special{ps: grestore}}}}
\put(2536,-4816){\makebox(0,0)[lb]{\smash{\SetFigFont{10}{12.0}{\rmdefault}{\mddefault}{\updefault}\special{ps: gsave 0 0 0 setrgbcolor}$\theta_1$\special{ps: grestore}}}}
\put(6106,-6901){\makebox(0,0)[lb]{\smash{\SetFigFont{10}{12.0}{\rmdefault}{\mddefault}{\updefault}\special{ps: gsave 0 0 0 setrgbcolor}$\theta_6$\special{ps: grestore}}}}
\put(-1739,-2536){\makebox(0,0)[lb]{\smash{\SetFigFont{10}{12.0}{\rmdefault}{\mddefault}{\updefault}\special{ps: gsave 0 0 0 setrgbcolor}$e_3$\special{ps: grestore}}}}
\put(-1379,-3301){\makebox(0,0)[lb]{\smash{\SetFigFont{10}{12.0}{\rmdefault}{\mddefault}{\updefault}\special{ps: gsave 0 0 0 setrgbcolor}$e_2$\special{ps: grestore}}}}
\put(-1004,-4801){\makebox(0,0)[lb]{\smash{\SetFigFont{10}{12.0}{\rmdefault}{\mddefault}{\updefault}\special{ps: gsave 0 0 0 setrgbcolor}$e_6$\special{ps: grestore}}}}
\put(-2714,-3871){\makebox(0,0)[lb]{\smash{\SetFigFont{10}{12.0}{\rmdefault}{\mddefault}{\updefault}\special{ps: gsave 0 0 0 setrgbcolor}$e_1$\special{ps: grestore}}}}
\put(-1229,-4081){\makebox(0,0)[lb]{\smash{\SetFigFont{10}{12.0}{\rmdefault}{\mddefault}{\updefault}\special{ps: gsave 0 0 0 setrgbcolor}$e_5$\special{ps: grestore}}}}
\put(-179,-3541){\makebox(0,0)[lb]{\smash{\SetFigFont{10}{12.0}{\rmdefault}{\mddefault}{\updefault}\special{ps: gsave 0 0 0 setrgbcolor}$e_4$\special{ps: grestore}}}}
\end{picture}

%% file: esiste1.pstex_t
\begin{picture}(0,0)%
\includegraphics{esiste1.pstex}%
\end{picture}%
\setlength{\unitlength}{2250sp}%
\begingroup\makeatletter\ifx\SetFigFont\undefined%
\gdef\SetFigFont#1#2#3#4#5{%
  \reset@font\fontsize{#1}{#2pt}%
  \fontfamily{#3}\fontseries{#4}\fontshape{#5}%
  \selectfont}%
\fi\endgroup%
\begin{picture}(10524,3499)(1789,-5173)
\put(3076,-4486){\makebox(0,0)[lb]{\smash{\SetFigFont{10}{12.0}{\rmdefault}{\mddefault}{\itdefault}$\theta_1$}}}
\put(5686,-3286){\makebox(0,0)[lb]{\smash{\SetFigFont{10}{12.0}{\rmdefault}{\mddefault}{\itdefault}$C_{234}$}}}
\put(5671,-1981){\makebox(0,0)[lb]{\smash{\SetFigFont{10}{12.0}{\rmdefault}{\mddefault}{\itdefault}$C_{135}$}}}
\put(5641,-4846){\makebox(0,0)[lb]{\smash{\SetFigFont{10}{12.0}{\rmdefault}{\mddefault}{\itdefault}$C_{126}$}}}
\put(11341,-1959){\makebox(0,0)[lb]{\smash{\SetFigFont{10}{12.0}{\rmdefault}{\mddefault}{\itdefault}$C_{135}$}}}
\put(12166,-3294){\makebox(0,0)[lb]{\smash{\SetFigFont{10}{12.0}{\rmdefault}{\mddefault}{\itdefault}$C_{456}$}}}
\put(11941,-4846){\makebox(0,0)[lb]{\smash{\SetFigFont{10}{12.0}{\rmdefault}{\mddefault}{\itdefault}$C_{126}$}}}
\put(3646,-4096){\makebox(0,0)[lb]{\smash{\SetFigFont{10}{12.0}{\rmdefault}{\mddefault}{\itdefault}$\theta_3$}}}
\put(3706,-4456){\makebox(0,0)[lb]{\smash{\SetFigFont{10}{12.0}{\rmdefault}{\mddefault}{\itdefault}$\theta_2$}}}
\put(9361,-4089){\makebox(0,0)[lb]{\smash{\SetFigFont{10}{12.0}{\rmdefault}{\mddefault}{\itdefault}$\theta_5$}}}
\put(9436,-4426){\makebox(0,0)[lb]{\smash{\SetFigFont{10}{12.0}{\rmdefault}{\mddefault}{\itdefault}$\theta_6$}}}
\put(8821,-4441){\makebox(0,0)[lb]{\smash{\SetFigFont{10}{12.0}{\rmdefault}{\mddefault}{\itdefault}$\theta_1$}}}
\put(2506,-4426){\makebox(0,0)[lb]{\smash{\SetFigFont{10}{12.0}{\rmdefault}{\mddefault}{\itdefault}$0$}}}
\put(8161,-4411){\makebox(0,0)[lb]{\smash{\SetFigFont{10}{12.0}{\rmdefault}{\mddefault}{\itdefault}$0$}}}
\end{picture}

%% file: esiste2.pstex_t
\begin{picture}(0,0)%
\includegraphics{esiste2.pstex}%
\end{picture}%
\setlength{\unitlength}{2368sp}%
\begingroup\makeatletter\ifx\SetFigFont\undefined%
\gdef\SetFigFont#1#2#3#4#5{%
  \reset@font\fontsize{#1}{#2pt}%
  \fontfamily{#3}\fontseries{#4}\fontshape{#5}%
  \selectfont}%
\fi\endgroup%
\begin{picture}(8424,5229)(3289,-5773)
\put(4081,-4966){\makebox(0,0)[lb]{\smash{\SetFigFont{10}{12.0}{\rmdefault}{\mddefault}{\updefault}$0$}}}
\put(9901,-1786){\makebox(0,0)[lb]{\smash{\SetFigFont{10}{12.0}{\rmdefault}{\mddefault}{\itdefault}$\theta_5$}}}
\put(11371,-886){\makebox(0,0)[lb]{\smash{\SetFigFont{10}{12.0}{\rmdefault}{\mddefault}{\itdefault}$C_{135}$}}}
\put(11317,-3145){\makebox(0,0)[lb]{\smash{\SetFigFont{10}{12.0}{\rmdefault}{\mddefault}{\itdefault}$C_{456}$}}}
\put(11686,-5446){\makebox(0,0)[lb]{\smash{\SetFigFont{10}{12.0}{\rmdefault}{\mddefault}{\itdefault}$C_{126}$}}}
\put(5701,-3241){\makebox(0,0)[lb]{\smash{\SetFigFont{10}{12.0}{\rmdefault}{\mddefault}{\itdefault}$C_{234}$}}}
\put(4786,-5056){\makebox(0,0)[lb]{\smash{\SetFigFont{10}{12.0}{\rmdefault}{\mddefault}{\itdefault}$\theta_1$}}}
\put(5266,-5041){\makebox(0,0)[lb]{\smash{\SetFigFont{10}{12.0}{\rmdefault}{\mddefault}{\itdefault}$\theta_2$}}}
\put(10966,-5026){\makebox(0,0)[lb]{\smash{\SetFigFont{10}{12.0}{\rmdefault}{\mddefault}{\itdefault}$\theta_6$}}}
\put(7576,-4006){\makebox(0,0)[lb]{\smash{\SetFigFont{10}{12.0}{\rmdefault}{\mddefault}{\itdefault}$\theta_4$}}}
\put(5311,-4726){\makebox(0,0)[lb]{\smash{\SetFigFont{10}{12.0}{\rmdefault}{\mddefault}{\itdefault}$\theta_3$}}}
\end{picture}

%% file: quadri.pstex_t
\begin{picture}(0,0)%
\includegraphics{quadri.pstex}%
\end{picture}%
\setlength{\unitlength}{2763sp}%
\begingroup\makeatletter\ifx\SetFigFont\undefined%
\gdef\SetFigFont#1#2#3#4#5{%
  \reset@font\fontsize{#1}{#2pt}%
  \fontfamily{#3}\fontseries{#4}\fontshape{#5}%
  \selectfont}%
\fi\endgroup%
\begin{picture}(8270,4290)(2993,-6226)
\put(4666,-4396){\makebox(0,0)[lb]{\smash{\SetFigFont{10}{12.0}{\rmdefault}{\mddefault}{\updefault}$F^*_{126}$}}}
\put(9031,-3976){\makebox(0,0)[lb]{\smash{\SetFigFont{10}{12.0}{\rmdefault}{\mddefault}{\updefault}$e_1$}}}
\put(10576,-5071){\makebox(0,0)[lb]{\smash{\SetFigFont{10}{12.0}{\rmdefault}{\mddefault}{\updefault}$e_5$}}}
\put(11041,-5341){\makebox(0,0)[lb]{\smash{\SetFigFont{10}{12.0}{\rmdefault}{\mddefault}{\updefault}$e_4$}}}
\put(9076,-5851){\makebox(0,0)[lb]{\smash{\SetFigFont{10}{12.0}{\rmdefault}{\mddefault}{\updefault}$e_6\in\calZ$}}}
\put(4186,-3421){\makebox(0,0)[lb]{\smash{\SetFigFont{10}{12.0}{\rmdefault}{\mddefault}{\updefault}$e_1$}}}
\put(5101,-3421){\makebox(0,0)[lb]{\smash{\SetFigFont{10}{12.0}{\rmdefault}{\mddefault}{\updefault}$e_2$}}}
\put(4051,-5236){\makebox(0,0)[lb]{\smash{\SetFigFont{10}{12.0}{\rmdefault}{\mddefault}{\updefault}$e_{16}$}}}
\put(5326,-5221){\makebox(0,0)[lb]{\smash{\SetFigFont{10}{12.0}{\rmdefault}{\mddefault}{\updefault}$e_{26}$}}}
\put(4651,-6226){\makebox(0,0)[lb]{\smash{\SetFigFont{10}{12.0}{\rmdefault}{\mddefault}{\updefault}$e_6$}}}
\put(10276,-4306){\makebox(0,0)[lb]{\smash{\SetFigFont{10}{12.0}{\rmdefault}{\mddefault}{\updefault}$e_2$}}}
\put(10621,-3841){\makebox(0,0)[lb]{\smash{\SetFigFont{10}{12.0}{\rmdefault}{\mddefault}{\updefault}$e_3$}}}
\put(9616,-3046){\makebox(0,0)[lb]{\smash{\SetFigFont{10}{12.0}{\rmdefault}{\mddefault}{\updefault}$v_{123}\in\calI$}}}
\put(4643,-2146){\makebox(0,0)[lb]{\smash{\SetFigFont{10}{12.0}{\rmdefault}{\mddefault}{\updefault}$e_{12}$}}}
\end{picture}

%% file: esiste3.pstex_t
\begin{picture}(0,0)%
\includegraphics{esiste3.pstex}%
\end{picture}%
\setlength{\unitlength}{2763sp}%
\begingroup\makeatletter\ifx\SetFigFont\undefined%
\gdef\SetFigFont#1#2#3#4#5{%
  \reset@font\fontsize{#1}{#2pt}%
  \fontfamily{#3}\fontseries{#4}\fontshape{#5}%
  \selectfont}%
\fi\endgroup%
\begin{picture}(8424,5517)(1189,-5473)
\put(1837,-5041){\makebox(0,0)[lb]{\smash{\SetFigFont{10}{12.0}{\rmdefault}{\mddefault}{\updefault}$\theta_1$}}}
\put(2461,-4966){\makebox(0,0)[lb]{\smash{\SetFigFont{10}{12.0}{\rmdefault}{\mddefault}{\updefault}$p_1$}}}
\put(5221,-3076){\makebox(0,0)[lb]{\smash{\SetFigFont{10}{12.0}{\rmdefault}{\mddefault}{\updefault}$p$}}}
\put(3515,-3027){\makebox(0,0)[lb]{\smash{\SetFigFont{10}{12.0}{\rmdefault}{\mddefault}{\updefault}$q_1$}}}
\put(5296,-1546){\makebox(0,0)[lb]{\smash{\SetFigFont{10}{12.0}{\rmdefault}{\mddefault}{\updefault}$\theta_3$}}}
\put(6271,-166){\makebox(0,0)[lb]{\smash{\SetFigFont{10}{12.0}{\rmdefault}{\mddefault}{\updefault}$C_{135}$}}}
\put(4861,-211){\makebox(0,0)[lb]{\smash{\SetFigFont{10}{12.0}{\rmdefault}{\mddefault}{\updefault}$C_{234}$}}}
\put(6556,-2266){\makebox(0,0)[lb]{\smash{\SetFigFont{10}{12.0}{\rmdefault}{\mddefault}{\updefault}$C_{456}$}}}
\put(6845,-3016){\makebox(0,0)[lb]{\smash{\SetFigFont{10}{12.0}{\rmdefault}{\mddefault}{\updefault}$q_2$}}}
\put(9016,-4621){\makebox(0,0)[lb]{\smash{\SetFigFont{10}{12.0}{\rmdefault}{\mddefault}{\updefault}$C_{126}$}}}
\put(8317,-5011){\makebox(0,0)[lb]{\smash{\SetFigFont{10}{12.0}{\rmdefault}{\mddefault}{\updefault}$\theta_2$}}}
\put(7711,-4981){\makebox(0,0)[lb]{\smash{\SetFigFont{10}{12.0}{\rmdefault}{\mddefault}{\updefault}$p_2$}}}
\put(5206,-4996){\makebox(0,0)[lb]{\smash{\SetFigFont{10}{12.0}{\rmdefault}{\mddefault}{\updefault}$q$}}}
\put(7111,-4606){\makebox(0,0)[lb]{\smash{\SetFigFont{10}{12.0}{\rmdefault}{\mddefault}{\updefault}$\theta_2'$}}}
\put(6811,-3586){\makebox(0,0)[lb]{\smash{\SetFigFont{10}{12.0}{\rmdefault}{\mddefault}{\updefault}$\theta_4$}}}
\put(3571,-3601){\makebox(0,0)[lb]{\smash{\SetFigFont{10}{12.0}{\rmdefault}{\mddefault}{\updefault}$\theta_5$}}}
\put(2995,-4210){\makebox(0,0)[lb]{\smash{\SetFigFont{10}{12.0}{\rmdefault}{\mddefault}{\updefault}$\theta_1''$}}}
\put(7201,-4261){\makebox(0,0)[lb]{\smash{\SetFigFont{10}{12.0}{\rmdefault}{\mddefault}{\updefault}$\theta_2''$}}}
\put(3067,-4606){\makebox(0,0)[lb]{\smash{\SetFigFont{10}{12.0}{\rmdefault}{\mddefault}{\updefault}$\theta_1'$}}}
\end{picture}

%% file: new_convessa.pstex_t
\begin{picture}(0,0)%
\includegraphics{new_convessa.pstex}%
\end{picture}%
\setlength{\unitlength}{2368sp}%
\begingroup\makeatletter\ifx\SetFigFont\undefined%
\gdef\SetFigFont#1#2#3#4#5{%
  \reset@font\fontsize{#1}{#2pt}%
  \fontfamily{#3}\fontseries{#4}\fontshape{#5}%
  \selectfont}%
\fi\endgroup%
\begin{picture}(5274,5026)(1489,-4936)
\put(4306,-4936){\makebox(0,0)[lb]{\smash{\SetFigFont{10}{12.0}{\rmdefault}{\mddefault}{\updefault}$0$}}}
\put(3631,-406){\makebox(0,0)[lb]{\smash{\SetFigFont{10}{12.0}{\rmdefault}{\mddefault}{\updefault}$\calC'(\calO)$}}}
\put(2409,-1500){\makebox(0,0)[lb]{\smash{\SetFigFont{10}{12.0}{\rmdefault}{\mddefault}{\updefault}$x$}}}
\put(6071,-706){\makebox(0,0)[lb]{\smash{\SetFigFont{10}{12.0}{\rmdefault}{\mddefault}{\updefault}$y$}}}
\end{picture}

%% file: altezza.pstex_t
\begin{picture}(0,0)%
\includegraphics{altezza.pstex}%
\end{picture}%
\setlength{\unitlength}{3355sp}%
\begingroup\makeatletter\ifx\SetFigFont\undefined%
\gdef\SetFigFont#1#2#3#4#5{%
  \reset@font\fontsize{#1}{#2pt}%
  \fontfamily{#3}\fontseries{#4}\fontshape{#5}%
  \selectfont}%
\fi\endgroup%
\begin{picture}(6174,2952)(1789,-5773)
\put(1876,-5221){\makebox(0,0)[lb]{\smash{\SetFigFont{10}{12.0}{\rmdefault}{\mddefault}{\updefault}$r_1$}}}
\put(5041,-3031){\makebox(0,0)[lb]{\smash{\SetFigFont{10}{12.0}{\rmdefault}{\mddefault}{\updefault}$v$}}}
\put(6031,-3016){\makebox(0,0)[lb]{\smash{\SetFigFont{10}{12.0}{\rmdefault}{\mddefault}{\updefault}$v_S$}}}
\put(3631,-3031){\makebox(0,0)[lb]{\smash{\SetFigFont{10}{12.0}{\rmdefault}{\mddefault}{\updefault}$v_{S'}$}}}
\put(6901,-5626){\makebox(0,0)[lb]{\smash{\SetFigFont{10}{12.0}{\rmdefault}{\mddefault}{\updefault}$r_3$}}}
\put(3451,-4321){\makebox(0,0)[lb]{\smash{\SetFigFont{10}{12.0}{\rmdefault}{\mddefault}{\updefault}$B'$}}}
\put(6256,-5251){\makebox(0,0)[lb]{\smash{\SetFigFont{10}{12.0}{\rmdefault}{\mddefault}{\updefault}$B$}}}
\put(5566,-5491){\makebox(0,0)[lb]{\smash{\SetFigFont{10}{12.0}{\rmdefault}{\mddefault}{\updefault}$S$}}}
\put(2701,-4696){\makebox(0,0)[lb]{\smash{\SetFigFont{10}{12.0}{\rmdefault}{\mddefault}{\updefault}$S'$}}}
\put(7651,-4411){\makebox(0,0)[lb]{\smash{\SetFigFont{10}{12.0}{\rmdefault}{\mddefault}{\updefault}$y_3(v)$}}}
\end{picture}

%% file: ratio_radii.pstex_t
\begin{picture}(0,0)%
\includegraphics{ratio_radii.pstex}%
\end{picture}%
\setlength{\unitlength}{3947sp}%
\begingroup\makeatletter\ifx\SetFigFont\undefined%
\gdef\SetFigFont#1#2#3#4#5{%
  \reset@font\fontsize{#1}{#2pt}%
  \fontfamily{#3}\fontseries{#4}\fontshape{#5}%
  \selectfont}%
\fi\endgroup%
\begin{picture}(3146,3146)(3873,-8339)
\put(5913,-5628){\makebox(0,0)[lb]{\smash{\SetFigFont{10}{12.0}{\rmdefault}{\mddefault}{\updefault}$\theta_{13}$}}}
\put(5228,-7246){\makebox(0,0)[lb]{\smash{\SetFigFont{10}{12.0}{\rmdefault}{\mddefault}{\updefault}$C_3$}}}
\put(5472,-7639){\makebox(0,0)[lb]{\smash{\SetFigFont{10}{12.0}{\rmdefault}{\mddefault}{\updefault}$\theta_{24}$}}}
\put(6217,-7678){\makebox(0,0)[lb]{\smash{\SetFigFont{10}{12.0}{\rmdefault}{\mddefault}{\updefault}$\theta_{12}$}}}
\put(4853,-6001){\makebox(0,0)[lb]{\smash{\SetFigFont{10}{12.0}{\rmdefault}{\mddefault}{\updefault}$\theta_{34}$}}}
\put(4169,-6798){\makebox(0,0)[lb]{\smash{\SetFigFont{10}{12.0}{\rmdefault}{\mddefault}{\updefault}$\theta_{14}$}}}
\put(6300,-6278){\makebox(0,0)[lb]{\smash{\SetFigFont{10}{12.0}{\rmdefault}{\mddefault}{\updefault}$\theta_{23}$}}}
\put(6060,-6668){\makebox(0,0)[lb]{\smash{\SetFigFont{10}{12.0}{\rmdefault}{\mddefault}{\updefault}$C_4$}}}
\put(4659,-7598){\makebox(0,0)[lb]{\smash{\SetFigFont{10}{12.0}{\rmdefault}{\mddefault}{\updefault}$C_1$}}}
\put(5372,-6068){\makebox(0,0)[lb]{\smash{\SetFigFont{10}{12.0}{\rmdefault}{\mddefault}{\updefault}$C_2$}}}
\put(4171,-5581){\makebox(0,0)[lb]{\smash{\SetFigFont{10}{12.0}{\rmdefault}{\mddefault}{\updefault}$C$}}}
\end{picture}

%% file: con1.pstex_t
\begin{picture}(0,0)%
\includegraphics{con1.pstex}%
\end{picture}%
\setlength{\unitlength}{3355sp}%
\begingroup\makeatletter\ifx\SetFigFont\undefined%
\gdef\SetFigFont#1#2#3#4#5{%
  \reset@font\fontsize{#1}{#2pt}%
  \fontfamily{#3}\fontseries{#4}\fontshape{#5}%
  \selectfont}%
\fi\endgroup%
\begin{picture}(6724,2309)(489,-4761)
\put(1456,-4761){\makebox(0,0)[lb]{\smash{\SetFigFont{10}{12.0}{\rmdefault}{\mddefault}{\updefault}{$-1$}%
}}}
\put(4146,-4761){\makebox(0,0)[lb]{\smash{\SetFigFont{10}{12.0}{\rmdefault}{\mddefault}{\updefault}{$c-a$}%
}}}
\put(5041,-4761){\makebox(0,0)[lb]{\smash{\SetFigFont{10}{12.0}{\rmdefault}{\mddefault}{\updefault}{$c+a$}%
}}}
\put(5806,-4761){\makebox(0,0)[lb]{\smash{\SetFigFont{10}{12.0}{\rmdefault}{\mddefault}{\updefault}{$1$}%
}}}
\put(6161,-4761){\makebox(0,0)[lb]{\smash{\SetFigFont{10}{12.0}{\rmdefault}{\mddefault}{\updefault}{$x$}%
}}}
\put(3841,-3886){\makebox(0,0)[lb]{\smash{\SetFigFont{10}{12.0}{\rmdefault}{\mddefault}{\updefault}{$\beta_2$}%
}}}
\put(4764,-3991){\makebox(0,0)[lb]{\smash{\SetFigFont{10}{12.0}{\rmdefault}{\mddefault}{\updefault}{$p_2$}%
}}}
\put(4423,-3444){\makebox(0,0)[lb]{\smash{\SetFigFont{10}{12.0}{\rmdefault}{\mddefault}{\updefault}{$\beta_1$}%
}}}
\put(5634,-3340){\makebox(0,0)[lb]{\smash{\SetFigFont{10}{12.0}{\rmdefault}{\mddefault}{\updefault}{$p_1$}%
}}}
\put(2371,-2731){\makebox(0,0)[lb]{\smash{\SetFigFont{10}{12.0}{\rmdefault}{\mddefault}{\updefault}{$\gamma_1$}%
}}}
\put(4141,-4456){\makebox(0,0)[lb]{\smash{\SetFigFont{10}{12.0}{\rmdefault}{\mddefault}{\updefault}{$\gamma_2$}%
}}}
\end{picture}

%% file: tetra.pstex_t
\begin{picture}(0,0)%
\includegraphics{tetra.pstex}%
\end{picture}%
\setlength{\unitlength}{2565sp}%
\begingroup\makeatletter\ifx\SetFigFont\undefined%
\gdef\SetFigFont#1#2#3#4#5{%
  \reset@font\fontsize{#1}{#2pt}%
  \fontfamily{#3}\fontseries{#4}\fontshape{#5}%
  \selectfont}%
\fi\endgroup%
\begin{picture}(4365,3900)(1891,-8296)
\put(2896,-8296){\makebox(0,0)[lb]{\smash{\SetFigFont{10}{12.0}{\rmdefault}{\mddefault}{\updefault}\special{ps: gsave 0 0 0 setrgbcolor}$P$\special{ps: grestore}}}}
\put(3031,-7546){\makebox(0,0)[lb]{\smash{\SetFigFont{10}{12.0}{\rmdefault}{\mddefault}{\updefault}\special{ps: gsave 0 0 0 setrgbcolor}$P'$\special{ps: grestore}}}}
\put(2431,-7156){\makebox(0,0)[lb]{\smash{\SetFigFont{10}{12.0}{\rmdefault}{\mddefault}{\updefault}\special{ps: gsave 0 0 0 setrgbcolor}$A$\special{ps: grestore}}}}
\put(3301,-5221){\makebox(0,0)[lb]{\smash{\SetFigFont{10}{12.0}{\rmdefault}{\mddefault}{\updefault}\special{ps: gsave 0 0 0 setrgbcolor}$A'$\special{ps: grestore}}}}
\put(3781,-4606){\makebox(0,0)[lb]{\smash{\SetFigFont{10}{12.0}{\rmdefault}{\mddefault}{\updefault}\special{ps: gsave 0 0 0 setrgbcolor}$u_1$\special{ps: grestore}}}}
\put(5311,-8266){\makebox(0,0)[lb]{\smash{\SetFigFont{10}{12.0}{\rmdefault}{\mddefault}{\updefault}\special{ps: gsave 0 0 0 setrgbcolor}$u_3$\special{ps: grestore}}}}
\put(6256,-6316){\makebox(0,0)[lb]{\smash{\SetFigFont{10}{12.0}{\rmdefault}{\mddefault}{\updefault}\special{ps: gsave 0 0 0 setrgbcolor}$u_4$\special{ps: grestore}}}}
\put(1891,-8191){\makebox(0,0)[lb]{\smash{\SetFigFont{10}{12.0}{\rmdefault}{\mddefault}{\updefault}\special{ps: gsave 0 0 0 setrgbcolor}$u_2$\special{ps: grestore}}}}
\put(2731,-6316){\makebox(0,0)[lb]{\smash{\SetFigFont{10}{12.0}{\rmdefault}{\mddefault}{\updefault}\special{ps: gsave 0 0 0 setrgbcolor}$\theta_{12}$\special{ps: grestore}}}}
\put(4666,-6286){\makebox(0,0)[lb]{\smash{\SetFigFont{10}{12.0}{\rmdefault}{\mddefault}{\updefault}\special{ps: gsave 0 0 0 setrgbcolor}$\theta_{13}$\special{ps: grestore}}}}
\put(5056,-5386){\makebox(0,0)[lb]{\smash{\SetFigFont{10}{12.0}{\rmdefault}{\mddefault}{\updefault}\special{ps: gsave 0 0 0 setrgbcolor}$\theta_{14}$\special{ps: grestore}}}}
\put(5851,-7276){\makebox(0,0)[lb]{\smash{\SetFigFont{10}{12.0}{\rmdefault}{\mddefault}{\updefault}\special{ps: gsave 0 0 0 setrgbcolor}$\theta_{34}$\special{ps: grestore}}}}
\put(3886,-8281){\makebox(0,0)[lb]{\smash{\SetFigFont{10}{12.0}{\rmdefault}{\mddefault}{\updefault}\special{ps: gsave 0 0 0 setrgbcolor}$\theta_{23}$\special{ps: grestore}}}}
\put(3931,-7096){\makebox(0,0)[lb]{\smash{\SetFigFont{10}{12.0}{\rmdefault}{\mddefault}{\updefault}\special{ps: gsave 0 0 0 setrgbcolor}$\theta_{24}$\special{ps: grestore}}}}
\end{picture}

%% file: mapimove.pstex_t
\begin{picture}(0,0)%
\includegraphics{mapimove.pstex}%
\end{picture}%
\setlength{\unitlength}{2368sp}%
\begingroup\makeatletter\ifx\SetFigFont\undefined%
\gdef\SetFigFont#1#2#3#4#5{%
  \reset@font\fontsize{#1}{#2pt}%
  \fontfamily{#3}\fontseries{#4}\fontshape{#5}%
  \selectfont}%
\fi\endgroup%
\begin{picture}(7980,4230)(1291,-5131)
\put(2701,-5131){\makebox(0,0)[lb]{\smash{\SetFigFont{10}{12.0}{\rmdefault}{\mddefault}{\updefault}$u_4$}}}
\put(6106,-2716){\makebox(0,0)[lb]{\smash{\SetFigFont{10}{12.0}{\rmdefault}{\mddefault}{\updefault}$u_1$}}}
\put(7636,-1126){\makebox(0,0)[lb]{\smash{\SetFigFont{10}{12.0}{\rmdefault}{\mddefault}{\updefault}$u_0$}}}
\put(7561,-5101){\makebox(0,0)[lb]{\smash{\SetFigFont{10}{12.0}{\rmdefault}{\mddefault}{\updefault}$u_4$}}}
\put(2881,-1111){\makebox(0,0)[lb]{\smash{\SetFigFont{10}{12.0}{\rmdefault}{\mddefault}{\updefault}$u_0$}}}
\put(4441,-2446){\makebox(0,0)[lb]{\smash{\SetFigFont{10}{12.0}{\rmdefault}{\mddefault}{\updefault}$u_3$}}}
\put(9271,-2446){\makebox(0,0)[lb]{\smash{\SetFigFont{10}{12.0}{\rmdefault}{\mddefault}{\updefault}$u_3$}}}
\put(1291,-2746){\makebox(0,0)[lb]{\smash{\SetFigFont{10}{12.0}{\rmdefault}{\mddefault}{\updefault}$u_1$}}}
\put(3241,-3601){\makebox(0,0)[lb]{\smash{\SetFigFont{10}{12.0}{\rmdefault}{\mddefault}{\updefault}$u_2$}}}
\put(8026,-3601){\makebox(0,0)[lb]{\smash{\SetFigFont{10}{12.0}{\rmdefault}{\mddefault}{\updefault}$u_2$}}}
\end{picture}